\documentclass[openany, amssymb, psamsfonts]{amsart}
\usepackage{mathrsfs,comment,amsmath}
\usepackage[usenames,dvipsnames]{color}
\usepackage[normalem]{ulem}
\usepackage{amssymb}
\usepackage{url}
\usepackage[all,arc,2cell]{xy}
\UseAllTwocells
\usepackage{enumerate,enumitem}
\usepackage{hyperref}  
\usepackage{braket}
\usepackage{framed}
\usepackage{graphicx}
\usepackage{bbm}
\usepackage{stmaryrd}
\usepackage{multirow}
\usepackage[table]{xcolor}
\usepackage{array}
\usepackage{caption}

\usepackage{tabularx}

\hypersetup{%
  bookmarksnumbered=true,%
  bookmarks=true,%
  colorlinks=true,%
  linkcolor=blue,%
  citecolor=blue,%
  filecolor=blue,%
  menucolor=blue,%
  urlcolor=blue,%
  pdfnewwindow=true,%
  pdfstartview=FitBH}   

\usepackage[utf8]{inputenc}
\usepackage{amsmath}
\usepackage{xcolor}
\usepackage{soul}
\usepackage{geometry}
\def\makeautorefname#1#2{\expandafter\def\csname#1autorefname\endcsname{#2}}
\def\equationautorefname~#1\null{(#1)\null}
\makeautorefname{footnote}{footnote}%
\makeautorefname{item}{item}%
\makeautorefname{figure}{Figure}%
\makeautorefname{table}{Table}%
\makeautorefname{part}{Part}%
\makeautorefname{appendix}{Appendix}%
\makeautorefname{chapter}{Chapter}%
\makeautorefname{section}{Section}%
\makeautorefname{subsection}{Section}%
\makeautorefname{subsubsection}{Section}%
\makeautorefname{theorem}{Theorem}%
\makeautorefname{thm}{Theorem}%
\makeautorefname{cor}{Corollary}%
\makeautorefname{lem}{Lemma}%
\makeautorefname{prop}{Proposition}%
\makeautorefname{pro}{Property}
\makeautorefname{conj}{Conjecture}%
\makeautorefname{defn}{Definition}%
\makeautorefname{notn}{Notation}
\makeautorefname{fact}{Fact}
\makeautorefname{notns}{Notations}
\makeautorefname{rem}{Remark}%
\makeautorefname{quest}{Question}%
\makeautorefname{exmp}{Example}%
\makeautorefname{ax}{Axiom}%
\makeautorefname{claim}{Claim}%
\makeautorefname{ass}{Assumption}%
\makeautorefname{asss}{Assumptions}%
\makeautorefname{con}{Construction}%
\makeautorefname{prob}{Problem}%
\makeautorefname{warn}{Warning}%
\makeautorefname{obs}{Observation}%
\makeautorefname{conv}{Convention}%

\newtheorem{thm}{Theorem}[section]
\newtheorem{cor}{Corollary}[section]
\newtheorem{prop}{Proposition}[section]
\newtheorem{lem}{Lemma}[section]

\theoremstyle{definition}
\newtheorem{defn}{Definition}[section]

\newtheorem{rem}{Remark}[section]

\makeatletter
\let\c@obs=\c@thm
\let\c@cor=\c@thm
\let\c@prop=\c@thm
\let\c@lem=\c@thm
\let\c@prob=\c@thm
\let\c@con=\c@thm
\let\c@conj=\c@thm
\let\c@defn=\c@thm
\let\c@fact=\c@thm
\let\c@notn=\c@thm
\let\c@notns=\c@thm
\let\c@exmp=\c@thm
\let\c@ax=\c@thm
\let\c@pro=\c@thm
\let\c@ass=\c@thm
\let\c@warn=\c@thm
\let\c@rem=\c@thm
\let\c@sch=\c@thm
\numberwithin{equation}{section}
\makeatother

\newcommand{\abs}[1]{\left|#1\right|}
\newcommand{\inorm}[1]{[\![#1]\!]}
\newcommand{\norm}[1]{\left\lVert#1\right\rVert}

\newcommand{\R}{\mathbb{R}}

\newcommand{\C}{\mathbb{C}}

\newcommand{\CC}{\mathcal{C}}

\newcommand{\SSS}{\mathcal{S}}
\newcommand{\Z}{\mathbb{Z}}

\newcommand{\bracket}[1]{\left(#1\right)}
\newcommand{\cbracket}[1]{\left\{ #1\right\}}
\newcommand{\sbracket}[1]{\left[ #1\right]}
\newcommand{\inner}[2]{\left\langle #1,#2 \right\rangle}

\newcommand{\dd}{\mathrm{d}}

\usepackage{blindtext}

\begin{document}
\setcounter{tocdepth}{1}
\title{Homogeneous Polynomial Carleson Operators along the Paraboloid}

\author{Alina Zhu}
\begin{abstract}
In this paper, we prove boundedness on $L^p(\R^{n+1})$, $1<p<\infty$, of polynomial Carleson operators along the paraboloid with phase belonging to a class $Y$ of homogeneous polynomials. For $n,d\ge 2$, let $\mathcal{P}^{n,d}$ denote the class of all real-valued homogeneous polynomials in $n$ variables of degree $d$. We may take the class $Y$ to be an arbitrary linear subspace $Y \subseteq \mathcal{P}^{n,d}$ such that $\dim Y \le M_{n , d}$, where $M_{n , d}$ is a constant depending only on $n$ and $d$.
The main novelty of our result is that we may allow a large collection of homogeneous polynomials of the same degree, in contrast to earlier results that relied on distinct degree assumptions. The heart of our proof is a $TT^*$ method, and the derivation of new estimates for oscillatory integral kernels involving interactions between homogeneous polynomial phases and the paraboloid.
\end{abstract}

\maketitle
\tableofcontents

\section{Introduction}

In this paper, we prove an \textit{a priori} bound on $L^p$, $1< p < \infty$, for maximally modulated singular Radon transforms over the paraboloid. Let $n , d\ge 1$ be fixed integers and $\R[X_1, ... ,X_n]$ be the vector space of real-valued polynomials in $n$ variables. For each $P\in \R[X_1,..., X_n]$,
define an integral operator, initially acting on Schwartz functions $f \in \mathcal{S}(\R^{n+1})$, by
\begin{align}
  R_Pf(x, t) : = \mathrm{p.v.}\int_{\R^n} e^{i P(y)}f(x - y, t - \abs{y}^2)K(y) ~\dd y~.\label{radontransformdef}
\end{align} 
Here, $K$ is a Calderón-Zygmund kernel, that is, a tempered distribution agreeing with a $\CC^1$ function $K$ on $x\neq 0$, such that $\hat{K} \in L^\infty(\R^n)$, and $\abs{\partial^\alpha K(x)} \le C \abs{x}^{-n - \abs{\alpha}}$ for all $0 \le \abs{\alpha} \le 1$ and a uniform constant $C$. 
We consider the corresponding maximally modulated operator,
\begin{align}
  f \longmapsto \sup_{P \in Y} \abs{R_Pf}~,\label{operator}
\end{align}
in which $Y \subseteq \R[X_1,...,X_n]$ is a chosen subspace. Such operators were first introduced in \cite{MR3945729}. It is reasonable to expect that the operator in \eqref{operator} satisfies an \textit{a priori} bound on $L^p(\R^{n+1})$, $1 < p < \infty$, for any subspace $Y$ of bounded degree. Let $\mathcal{P}^{n,d}\subseteq \R[X_1,...,X_n]$ denote the $\binom{d+n-1}{n-1}$ dimensional vector space of homogeneous polynomials in $n$ variables of degree $d$. We verify this conjecture for certain subspaces $Y \subseteq \mathcal{P}^{n,d}$.

\begin{thm}\label{mainresult2}
  Let $n , d\ge 2$ be fixed integers. Let $Y \subseteq \mathcal{P}^{n , d}$ be a linear subspace of dimension at most $M_{n,d}$, where
\begin{align*}
  M_{n,d} = \left\{\begin{array}{cl}
  n-1 & \mbox{if } d = 2~ \\[0.5cm] 
  \dbinom{d+ n -3}{n-1} -1 & \mbox{if $d$ is even}~, ~d \ge 4 \\[0.5cm]
  \dbinom{d+ n -2}{n-1} -1 & \mbox{if $d$ is odd}~\\
  \end{array}\right.~.
\end{align*}
In the case when $d = 2$, assume furthermore that $\abs{X}^2\notin Y$, where $|X|^2$ denotes the radial polynomial $|X|^2 = X_1^2 + ... + X_n^2$.
Then, for each $1 < p< \infty$, there exists a constant $C$ which may depend on $n , d, p$, and a chosen basis of $Y$, such that
\begin{align*}
  \norm{ \sup_{P \in Y} \abs{R_{P} f}}_{L^p(\R^{n+1})} \le C \norm{f}_{L^p(\R^{n+1})}~,~\mbox{ for all $f \in \mathcal{S}(\R^{n+1})$~.}
\end{align*}
\end{thm}

The initial work of \cite{MR3945729} and later \cite{MR4790541} could consider an arbitrary linear span $Y$ of homogeneous polynomials in $n$ variables of distinct degrees, such that
\begin{align*}
    Y \cap \mathrm{span}_\R\{X_1,..., X_n , |X|^2\} = \{0\}~.
\end{align*}
In particular, if $Y$ consists of polynomials of degree at most $d$, then $\dim(Y) \le d-1$. 
In contrast, Theorem \ref{mainresult2} treats a large collection $Y$ of homogeneous polynomials of the same degree. In the case when $Y \subseteq \mathcal{P}^{n , d}$ is an arbitrary linear subspace, the upper bound $M_{n , d}$ for $\dim(Y)$ admissible in Theorem \ref{mainresult2} is allowed to grow polynomially with respect to both $n$ and $d$.
Table \ref{tab:Mnd} below contrasts for small $n , d$ the largest permissible dimension of $Y$ in \cite{MR3945729} (i.e., $d-1$) and the largest permissible dimension of $Y$ in Theorem \ref{mainresult2} (i.e., $M_{n , d}$).
\begin{table}[ht]
\centering
\fontsize{8pt}{12pt}\selectfont
\setlength{\tabcolsep}{1.4pt}
\renewcommand{\arraystretch}{1.5}
\begin{tabular}{|
  >{\columncolor{gray!20}\centering\arraybackslash}m{1cm}|
  *{25}{>{\centering\arraybackslash}m{0.5cm}|}}
\hline
\rowcolor{gray!20}
$n$
& \multicolumn{5}{c|}{\cellcolor{gray!20}2}
& \multicolumn{5}{c|}{\cellcolor{gray!20}3}
& \multicolumn{5}{c|}{\cellcolor{gray!20}4}
& \multicolumn{5}{c|}{\cellcolor{gray!20}5}\\
\hline
\rowcolor{gray!20}
$d$
& 2 & 3 & 4 & 5 & 6
& 2 & 3 & 4 & 5 & 6
& 2 & 3 & 4 & 5 & 6
& 2 & 3 & 4 & 5 & 6 \\
\hline
$M_{n,d}$
& 1 & 2 & 2 & 4 & 4
& 2 & 5 & 5 & 14 & 14
& 3 & 9 & 9 & 34 & 34
& 4 & 14 & 14 & 69 & 69\\
\hline
$d-1$
& 1 & 2 & 3 & 4 & 5
& 1 & 2 & 3 & 4 & 5
& 1 & 2 & 3 & 4 & 5
& 1 & 2 & 3 & 4 & 5\\
\hline
\end{tabular}
\caption{Comparison of $M_{n,d}$ and $d-1$ for $2\le n\le5,$ $2\le d\le6$.}
\label{tab:Mnd}
\end{table}

To illustrate the growth of $M_{n,d}$ with respect to $n$ and $d$, we compute the leading term in its asymptotic expansion as $n , d \to \infty$.
For $k,s \in \Z_{\ge 0}$, the leading term of $\binom{k+s}{s}$ as $k \to \infty$ is given by
\begin{align*}
\binom{k+s}{s} = \frac{(k+s)(k+s-1) \cdots (k+1)}{s!} = \frac{k^s}{s!} + O_s(k^{s-1})~.
\end{align*}
Applying this to the definition of $M_{n , d}$ in Theorem \ref{mainresult2}, for fixed $n\ge 2$, the leading term as $d \to \infty$ of $M_{n,d}$ is $d^{n-1}/(n-1)!$. For fixed $d\ge 3$, as $n \to \infty$, the leading term of $M_{n,d}$ is $n^{d-1}/(d-1)!$ if $d$ is odd, and $n^{d-2}/(d-2)!$ if $d$ is even. For comparison, the leading terms of $\dim(\mathcal{P}^{n,d}) = \binom{d+n-1}{n-1}$ as $d \to \infty$ and $n \to \infty$ are $d^{n-1}/(n-1)!$ and $n^d/d!$ respectively.

We may also compute a universal lower bound for $M_{n , d}$ which is satisfied for all $n , d \ge 2$.
Since the function $\binom{d+n-i}{n-1}$ is increasing in both $n$ and $d$ for $i = 2,3$, we have for each $n \ge 2 $, $d \ge 3$ the lower bound
\begin{align*}
  \binom{d+n-i}{n-1} \ge \binom{d+2 - i}{1} ~,~ \mbox{ and }~ \binom{d+n-i}{n-1} = \binom{d+n - i}{d+ 1 - i} \ge \binom{n + 3 - i}{4-i}~.
\end{align*}
Hence, when $n = 2$, $M_{n,d}$ admits the following lower bound for all $d \ge 2$:
\begin{align*}
    M_{n,d} \ge \begin{cases}d-1 & \mbox{ if $d$ is odd} \\[0.2cm]
    \max\{d-2, 1\} & \mbox{ if $d$ is even} \\[0.2cm]
    \end{cases}~.
\end{align*}
If $n \ge 3$, we have the improved lower bound for all $d \ge 2$:
\begin{align*}
    M_{n,d} \ge \begin{cases}
    ~~\max\cbracket{d(d+1)/2 -1 ~,~ n(n+1)/2 -1} & \mbox{ if $d$ is odd} \\[0.2cm]
    ~~\max\cbracket{d(d-1)/2 -1~,~ n(n+1)/2 -1} & \mbox{ if $d$ is even , $d \ge 4$} \\[0.2cm]
    ~~n-1 & \mbox{ if $d = 2$} \\[0.2cm]
    \end{cases}~.
\end{align*}
Therefore, when $n\ge3$, our dimension bound is strictly stronger than the dimension constraint of $d-1$ given in \cite{MR3945729}.

The bound $M_{n , d}$ on $\dim(Y)$ in Theorem \ref{mainresult2} does not indicate a genuine failure of boundedness of \eqref{operator} but rather a natural limitation of our method. More specifically, this constraint encodes intrinsic properties of integration kernels obtained through the $TT^*$ method, whose decay is central to our proof. As will be explained later, we cannot demonstrate decay of the required kernels for certain subspaces $Y\subseteq \mathcal{P}^{n,d}$ of dimension exceeding $M_{n,d}$. The dependence of the constant $C$ on the choice of basis of $Y$ is also a natural limitation of our method, as the definition of the aforementioned integration kernel is also sensitive to the basis of $Y$.

\subsection{A More General Theorem} \label{corollaries}
The proof of Theorem \ref{mainresult2} relies on the observation that any linear subspace $Y \subseteq \mathcal{P}^{n , d}$ admissible in Theorem \ref{mainresult2} can be embedded into a codimension $n$ subspace $\mathcal{Q} \subseteq \mathcal{P}^{n , d}$ satisfying certain nondegeneracy conditions.
To state the nondegeneracy conditions, we use the following notation. For $c \in \R^n$, $c\cdot X$ denotes the linear polynomial $c \cdot X = c_1X_1 + ... + c_nX_n$, and $\abs{X}^2$ denotes the radial polynomial $|X|^2 =X_1^2 + ... + X_n^2$.
\begin{defn}
Let $c \in \R^n$ and let $\mathcal{Q}\subseteq \mathcal{P}^{n,d}$ be a linear subspace. We say that $\mathcal{Q}$  satisfies the \textbf{tangential nondegeneracy condition} for $c$ if
\begin{align}
  \mathcal{Q} \cap \mathrm{span}_\R \cbracket{ X_j(c \cdot X)^{d-1}~~|~~ 1 \le j \le n} = \{0\}~.\label{cond1}
\end{align}
We say that $\mathcal{Q}$ satisfies the \textbf{radial nondegeneracy condition} for $c$ if
\begin{align}
  \mathcal{Q} \cap \mathrm{span}_\R \cbracket{ (c \cdot X)^d~,~ \abs{X}^{2}(c \cdot X)^{d - 2}} = \{0\}~. \label{cond2}
\end{align}
\end{defn}
The reason for why we named conditions \eqref{cond1} and \eqref{cond2} ``tangential'' and ``radial'' is related to how certain orthogonal projection maps behave once we fix a particular inner product on $\mathcal{P}^{n , d}$; we refer to Remark \ref{conditions} for a more detailed discussion.
We now state our main result in its full generality.

\begin{thm}\label{mainresult}
Fix a dimension $n \ge 2$ and a degree $d \ge 2$. Let $\mathcal{Q}\subseteq \mathcal{P}^{n,d}$ be a codimension $n$ linear subspace and fix a basis $\mathfrak{B} $ of $\mathcal{Q}$.
Suppose there exists $c \in (\R\backslash\{0\})^n$ 
such that $\mathcal{Q}$ satisfies the tangential nondegeneracy condition \eqref{cond1} for $c$. If $d$ is even, assume further that $\mathcal{Q}$ simultaneously satisfies the radial nondegeneracy condition \eqref{cond2} for the same $c$.
Then, for each $1 < p< \infty$, the following \textit{a priori} estimate holds for each Schwartz function $f \in \mathcal{S}(\R^{n+1})$:
\begin{align*}
  \norm{ \sup_{P \in \mathcal{Q}} \abs{R_{P} f}}_{L^p(\R^{n+1})} \le C \norm{f}_{L^p(\R^{n+1})}~,
\end{align*}
with a constant $C$ that may depend on $n , d, p$, and the basis $\mathfrak{B} $.

\end{thm}
Theorem \ref{mainresult2} directly follows from Theorem \ref{mainresult} once we verify the following proposition.
\begin{prop}\label{derivationofmainthm}
Let $n , d \ge 2$ be fixed integers.
Let $Y\subseteq \mathcal{P}^{n , d}$ be an arbitrary linear subspace. Then the following holds:
\begin{itemize}
\item[(i)] If $\dim(Y) < \binom{d+ n -2}{n-1}  $, then there exists a codimension $n$ subspace $\mathcal{Q} \subseteq \mathcal{P}^{n , d}$ containing $Y$ and an open dense set $U\subseteq \R^n$ such that $\mathcal{Q}$ satisfies the tangential nondegeneracy condition \eqref{cond1} for all $c \in U$.
\item[(ii)] If $d \ge 4$ is even and $\dim(Y) < \binom{d+ n -3}{n-1} $, then there exists a codimension $n$ subspace $\mathcal{Q} \subseteq \mathcal{P}^{n , d}$ containing $Y$ and $c \in (\R \backslash \{0\})^n$ such that $\mathcal{Q}$ satisfies both the tangential nondegeneracy condition \eqref{cond1} and the radial nondegeneracy condition \eqref{cond2} for $c$.
\item[(iii)] If $d = 2$, $\dim(Y) < n$, and $|X|^2 \notin Y$, then there exists a codimension $n$ subspace $\mathcal{Q} \subseteq \mathcal{P}^{n , 2}$ containing $Y$ and $c \in (\R \backslash \{0\})^n$ such that $\mathcal{Q}$ satisfies both the tangential nondegeneracy condition \eqref{cond1} and the radial nondegeneracy condition \eqref{cond2} for $c$.
\end{itemize}
\end{prop}

The dimension upper bound $M_{n,d}$ proposed in Theorem \ref{mainresult2} is maximal in the context of embedding an arbitrary linear subspace $Y \subseteq \mathcal{P}^{n,d}$ into a codimension $n$ subspace $\mathcal{Q} \subseteq \mathcal{P}^{n , d}$ which satisfies the hypothesis of Theorem \ref{mainresult}. Problematic examples start to appear once we allow the dimension of $Y$ to exceed $M_{n,d}$. Let $n , d \ge 1$ be fixed integers and $M_{n , d}$ as defined in Theorem \ref{mainresult2}. If $l \in \mathcal{P}^{n , 1} \backslash \{0\}$ is a nonzero linear form,  then the subspace $Y = l\mathcal{P}^{n , d-1}$ has dimension $\binom{d+n-2}{n-1} > M_{n , d}$. No codimension $n$ subspace $\mathcal{Q}$ containing such a $Y$ can satisfy the tangential nondegeneracy condition required for Theorem \ref{mainresult}, since
\begin{align*}
    l\mathcal{P}^{n , d-1} \cap \mathrm{span}_\R \cbracket{ X_j(c \cdot X)^{d-1}~~|~~ 1 \le j \le n} \neq \{0\}~,~ \mbox{ for all } c \in (\R \backslash \{0 \})^n~.
\end{align*}
If $d \ge 4$ is even, then the subspace $Y = \abs{X}^2\mathcal{P}^{n , d-2}$ has dimension $\binom{d+n-3}{n-1}= M_{n,d} + 1$, and we have a similar obstruction since
\begin{align*}
    \abs{X}^2\mathcal{P}^{n , d-2} \cap \mathrm{span}_\R \cbracket{ (c \cdot X)^d~,~ \abs{X}^{2}(c \cdot X)^{d - 2}} \neq \{0\}~,~ \mbox{ for all } c \in (\R \backslash \{0\})^n~.
\end{align*}
The assumption when $d = 2$ that $|X|^2 \notin Y$ in Theorem \ref{mainresult2} comes from the observation that the radial nondegeneracy condition when $d = 2$ \eqref{cond2} requires
\begin{align*}
  \mathcal{Q} \cap \mathrm{span}_\R \cbracket{ (c \cdot X)^2~,~ \abs{X}^{2}} = \{0\}~.
\end{align*}
The condition thus fails for all $c \in (\R\backslash\{0\})^n$ once $Y$ contains $|X|^2$.

As a second illustration of a setting in which Theorem \ref{mainresult} applies, we identify a class of codimension $n$ subspaces $\mathcal{Q} \subseteq \mathcal{P}^{n ,d}$ for which the nondegeneracy conditions admit a simple linear-algebraic characterization. 
\begin{prop}\label{monomials}
Fix $n , d \ge2$ and let
\begin{align*}
\cbracket{\alpha \in \Z^{n}_{\ge 0}~~|~~ \abs{\alpha} = d} = T \sqcup S
\end{align*}
be a partition of the set of multi-indices of degree $d$ in $n$ variables into two disjoint sets, where
\begin{align*}
    \abs{T} = n~,~ \abs{S} = \dim(\mathcal{P}^{n ,d}) - n~.
\end{align*}
Define the codimension $n$ linear subspace
\begin{align*}
  \mathcal{Q} = \mathrm{span}_\R \cbracket{ X^{\alpha} ~~|~~ \alpha \in S} \subseteq \mathcal{P}^{n,d}~.
\end{align*}  
If $T$ is a linearly independent set in $\R^n$, then $\mathcal{Q}$ satisfies the tangential nondegeneracy condition for every $c \in (\R \backslash \{0\})^n$. 
If $T$ further satisfies that there is a vector $\theta = (\theta_1, ..., \theta_n) \in T $ such that $\theta_j \ge 2$ for some $1\le j \le n$, then there exists $c \in (\R \backslash \{0\})^n$ for which $\mathcal{Q}$ satisfies the radial nondegeneracy condition for $c$.
\end{prop}
Boundedness of the maximal operator in \eqref{operator} can thus be directly derived for subspaces $\mathcal{Q} \subseteq \mathcal{P}^{n,d}$ described in Proposition \ref{monomials}.
\begin{cor}[Span of Monomials] \label{monomialcor}
  Let $n , d \ge 2$ be fixed integers. Fix a partition
  \begin{align*}
      \cbracket{\alpha \in \Z^{n}_{\ge 0}~~|~~ \abs{\alpha} = d} = T \sqcup S
  \end{align*}
  such that $|T| = n$. Suppose $T$ is a linearly independent set in $\R^n$, and, if $d$ is even, suppose further that there exists $\theta \in T$ with $\theta_j \ge 2$ for some $1 \le j \le n$. Consider the codimension $n$ subspace
  \begin{align*}
      \mathcal{Q} = \mathrm{span}_\R \cbracket{ X^{\alpha} ~~|~~ \alpha \in S} ~,
  \end{align*}
and the operator in \eqref{radontransformdef} defined for each $P \in \mathcal{Q}$. Then, for each $1 < p< \infty$, there exists a constant $C$ which may depend on $n , d, p$, such that
\begin{align*}
  \norm{ \sup_{P \in \mathcal{Q}} \abs{R_{P} f}}_{L^p(\R^{n+1})} \le C \norm{f}_{L^p(\R^{n+1})}~,~\mbox{ for all $f \in \mathcal{S}(\R^{n+1})$~.}
\end{align*}

\end{cor}

For example, when $d = 2$, we can take
\begin{align*}
  T = \cbracket{2e_i ~~|~~1 \le i \le n}~,~ S = \cbracket{e_i + e_j ~~|~~ 1 \le i < j \le n}~,
\end{align*}
so that $T$ is a linearly independent set in $\R^n$ and every element in $T$ has at least one coordinate $\ge 2$. Then, the codimension $n$ subspace 
\begin{align*}
    \mathcal{Q} = \mathrm{span}_\R \cbracket{X_iX_j , 1\le i< j \le n}
\end{align*}
is admissible for Corollary \ref{monomialcor}, and hence the corresponding operator \eqref{operator} for this $\mathcal{Q}$ is bounded on $L^p(\R^{n+1})$, for $1 < p< \infty$.

When $d = n$, there is exactly one multi-index of degree $d$ in $n$ variables such that all components are $\le 1$, namely the multi-index $(1,...,1)$, while when $d>n$, no multi-indices of degree $d$ in $n$ variables have all components $\le 1$ by the pigeonhole principle. Moreover, when $n = 2$, any two distinct multi-indices of degree $d$ in two variables are linearly independent. Therefore, any $d-1$ dimensional subspace $\mathcal{Q} \subseteq \mathcal{P}^{2,d}$ which admits a monomial basis is admissible for Corollary \ref{monomialcor}. We therefore have the following special case for $n = 2$.
\begin{cor} \label{R2exmp}
Fix $d \ge 2$. Let $S \subseteq \cbracket{\alpha \in \Z^{2}_{\ge 0}~~|~~ \abs{\alpha} = d}$ be an \textbf{arbitrary} set of multi-indices in two variables of degree $d$ with $|S| = d-1$. Then, the corresponding codimension $2$ subspace
\begin{align*}
    \mathcal{Q} : = \mathrm{span}_\R \cbracket{ X^{\alpha} ~~|~~ \alpha \in S}
\end{align*}
is admissible for Corollary \ref{monomialcor}, and thus generates a bounded Carleson operator along the paraboloid on $L^p(\R^3)$, for all $1 < p < \infty$. 
\end{cor}

\subsection{Relation to Previous Work}

Carleson operators originated from Carleson's celebrated theorem \cite{MR199631}, which established almost
everywhere convergence of Fourier series for functions in
$L^2\bracket{[0,1]}$. Define an integral operator for each polynomial $P \in \R[X_1,...,X_n]$ by 
\begin{align}
  T_Pf(x) = \mathrm{p.v.} \int_{\mathbb{R}^n} e^{iP(y)}f(x-y)K(y)~\dd y~, \label{carlesonoperator}
\end{align}
where $K$ is a Calderón-Zygmund kernel. 
In modern terminology, Carleson's theorem can be summarized as the $L^2$ boundedness of the maximal operator $f \mapsto \sup_{P \in Y} \abs{T_Pf}$
in the case when $n=1$, $K(y)=1/y$, and $Y = \mathrm{span}_{\R} \{X_1\}$.
Hunt \cite{MR238019} subsequently extended Carleson's result to $L^p(\R)$ for $1<p<\infty$, while Sj\"olin \cite{MR336222} proved the corresponding result in $\R^n$ for general $n$. Further approaches to the Carleson operator were developed by Fefferman \cite{MR340926} and Lacey and Thiele \cite{MR1783613}, introducing and developing the use of time-frequency analysis methods in the treatment of maximally modulated operators.

Stein expanded the scope of Carleson's theorem by asking whether the operator $f \mapsto \sup_{P \in Y} \abs{T_Pf}$ could still be bounded when $Y$ is taken to be a larger subset of polynomials of bounded degree. It was first shown by Ricci and Stein \cite{MR890662} that the operator
\begin{align*}
    Tf(x) = \mathrm{p.v.} \int_{\mathbb{R}^n} e^{iP(x,y)}f(x-y)K(y)~\dd y~, ~ P \in \R[X_1,..., X_n , Y_1,..., Y_n]~,
\end{align*}
is bounded on $L^p(\R^n)$ for $1<p<\infty$, with operator norm depending only on the degree of the polynomial $P$, but not on the coefficients. Stein \cite{MR1364908} then showed that the maximal operator $f \mapsto \sup_{P \in Y} \abs{T_Pf}$ where $T_P$ is defined as in \eqref{carlesonoperator} is bounded on $L^p(\R^n)$ for $1<p<\infty$, in the case $n = 1$, $K(y) = 1/y$, and $Y = \mathrm{span}_\R \{ X_1^2\}$. In 2001,
Stein and Wainger \cite{MR1879821} introduced an approach to polynomial Carleson operators based on the $TT^*$ method and Van der Corput estimates for oscillatory integrals, and proved the following result.
\begin{thm} \label{steinwaingerthm}
\cite[Theorem 2]{MR1879821} Assume $n\ge 1$, $d \ge 2$, and \begin{align*}
  Y = \mathrm{span}_\R \cbracket{ X^\alpha ~~|~~ 2 \le \abs{\alpha} \le d}~,
\end{align*}
where $X^\alpha = X_1^{\alpha_1} \cdots X_n^{\alpha_n}$ and $\abs{\alpha} = \alpha_1 + ... + \alpha_n$. Consider \eqref{carlesonoperator} defined for each $P \in Y$. Then, for each $1 < p< \infty$, the following estimate holds for each Schwartz function $f \in \mathcal{S}(\R^n)$:
\begin{align*}
  \norm{\sup_{P \in Y} \abs{T_Pf}}_{L^p(\R^n)} \le C_{n,p,d} \norm{f}_{L^p(\R^n)}~.
\end{align*}
\end{thm}
The exclusion of linear polynomials from $Y$ is a natural limitation for this method: through a stopping time argument and $TT^*$, Stein and Wainger derived their result from the pointwise decay of an oscillatory integral kernel roughly of the form
\begin{align*}
    K_\flat^{\nu , \mu} (u) = \int_{\R^n} e^{i(P_\nu(u+z) -P_\mu(z))} \Psi(u , z) ~\dd z~,~ \mbox{ for } \Psi \in \CC^1_c(\R^n \times \R^n)~ \mbox{ and } P_\nu , P_\mu \in Y~.
\end{align*}
The obstruction is that one cannot show decay of this kernel if $Y$ contains any nonzero linear polynomial;
see \cite[Page 794]{MR1879821} for a more detailed explanation. For this reason, Stein and Wainger's result strictly does not include the classical Carleson theorem, which corresponds to taking $Y = \mathrm{span}_\R \{ X_1,... , X_n\}$ in Theorem \ref{steinwaingerthm}. Via time-frequency analysis methods, the case in which
\begin{align*}
  Y = \mathrm{span}_\R \cbracket{ X^\alpha ~~|~~ 1 \le \abs{\alpha} \le d}~,~ d\ge 1~
\end{align*}
in Theorem \ref{steinwaingerthm} was later treated by Lie \cite{MR2545246,MR4125450} for $n=1$, and independently by Lie and Zorin-Kranich \cite{lie2017notepolynomialcarlesonoperator,MR4270523} for general $n$, leading to a complete resolution of the boundedness problem for polynomial Carleson operators in arbitrary dimensions.

Pierce and Yung \cite{MR3945729} introduced polynomial Carleson operators along the paraboloid defined in \eqref{operator}. Motivated by Stein and Wainger, they used a combination of square functions, Littlewood-Paley theory, and $TT^*$ methods, and proved the following result for $n \ge 2$. 
\begin{thm}\cite[Theorem 1.1]{MR3945729}\label{pierceyungthm}
  Assume $n , d \ge 2$, and \begin{align*}
  Y = \mathrm{span}_{\mathbb R} \{p_2,\ldots,p_d\}~,
\end{align*}
where $p_i\in\mathcal{P}^{n,i}$ for each $2 \le i \le d$, and $p_2 \notin \mathrm{span}_\R \{|X|^2\}$ if $p_2$ is not the zero polynomial.

Consider the operator defined in \eqref{operator}.
Then, for each $1 < p< \infty$, the following \textit{a priori} estimate holds for each Schwartz function $f \in \mathcal{S}(\R^{n+1})$:
\begin{align}
  \norm{\sup_{P \in Y} \abs{R_P f}}_{L^p(\R^{n+1})} \le C_{n , p , d } \norm{f}_{L^p(\R^{n+1})}~.\label{boundednesspierceyung}
\end{align}
The constant $C_{n , p , d }$ may depend on $n , p , d$, and the choice of basis $\{p_2,\ldots,p_d\}$ of $Y$.
\end{thm}
Central to the result of \cite{MR3945729} is the reduction of the estimate in \eqref{boundednesspierceyung} to a pointwise estimate of a particular oscillatory integral kernel defined by a subspace $Y$ which satisfies the hypothesis of Theorem \ref{pierceyungthm}. The further assumption that the basis elements of $Y$ have distinct homogeneous degrees is required for the estimates of the oscillatory integral kernel obtained after the reduction. 
Analogous to the obstruction present in \cite{MR1879821}, the requirement that $Y$ does not contain any linear terms arises again in \cite{MR3945729} due to the same obstruction about decay of oscillatory integral kernels obtained through $TT^*$.
The further assumption that $p_2\notin \mathrm{span}_\R \{|X|^2\} \backslash \{0\}$ is required for a similar reason: since the operator integrates along the paraboloid, one is also obstructed from showing decay for the associated oscillatory integral kernel in this case; see \cite[Pages 3-12]{MR3945729} for a more detailed discussion.

Following the work of Pierce and Yung, several authors have studied polynomial Carleson operators along submanifolds of $\R^{n+1}$, covering cases which are not included by their original setting.
We give a brief overview of some of these results. 
The result of \cite{MR3945729} was extended by Anderson, Maldague, Pierce, and Yung \cite{MR4790541} to the case when integration along the paraboloid $(y, |y|^2) $ in  \eqref{radontransformdef} is replaced by an arbitrary nondegenerate quadratic hypersurface $\gamma(y)$ in $\R^{n+1}$, that is $\gamma(y) = (y, y^TQ y)$ for an $n \times n$ symmetric matrix $Q$ such that $\det (Q) \neq 0$.

In the case $n=1$, the original method employed by \cite{MR3945729} is not directly applicable due to the lack of sufficient dimensions to isolate an oscillatory integral kernel. Combining Carleson's theorem with $TT^*$ estimates in the spirit of Stein and Wainger, Guo, Pierce, Roos, and Yung \cite[Theorem 1.1]{MR3708001} obtained partial results on polynomial Carleson operators along monomial curves in the plane. Using new oscillatory integral estimates developed in \cite{MR4998907}, Beltran, Guo, and Hickman \cite{MR5056223} subsequently proved $L^p(\R^{n+1})$, $1< p< \infty$, boundedness of the operator in \eqref{operator} when $n = 1$, $K(y) = 1/y$ and $Y = \mathrm{span}_\R\{ X_1^3\}$.

The approach of \cite{MR3945729} also cannot extend to the case in which $Y$ contains nonzero linear polynomials or $|X|^2$ in the operator \eqref{operator}. This case was treated via time-frequency analysis approaches. Partial results were obtained by Becker in \cite[Corollary A.1]{MR4687376}, which concerns a thickened version of \eqref{operator} for general $n\ge 1$ and the modulation space
\begin{align*}
    Y = \mathrm{span}_\R\{X^\alpha ~~|~~ 1 \le |\alpha| \le d\}, \quad d\ge1~,
\end{align*}
using improved constants for maximally modulated Calderón-Zygmund operators on the torus.
More recently, using time-frequency analysis methods, Hsu and Lie \cite{hsu2026resonantcarlesonradontransformdimensions} established $L^p(\mathbb{R}^{n+1})$ boundedness, for every $1<p<\infty$, of the operator \eqref{operator} for general $n\ge1$ and
\begin{align*}
  Y = \{a_1X_1 + ... + a_nX_n + a_{n+1}|X|^2~~|~~ a \in \mathbb{V}\}~,
\end{align*}
where $\mathbb{V} \subseteq \R^{n+1}$ is a linear subspace satisfying the condition that there is a nontrivial vector $v_0 \in \R^{n} \times \{0\}$ such that $v_0 \perp \mathbb{V}$. In fact, Hsu and Lie actually proved $L^p$ boundedness for a maximal operator which dominates \eqref{operator}; see \cite[Theorem 2.2]{hsu2026resonantcarlesonradontransformdimensions}.

The present paper continues this line of research.
Notice that the result of \cite{MR3945729} stated in Theorem \ref{pierceyungthm} does not permit the modulation space $Y$ to contain more than one linearly independent homogeneous polynomial of the same degree in $\R[X_1,...,X_n]$. A recent result of Zhang \cite[Theorem 4]{zhang2026oscillatoryintegraloperatorsparaboloid}, obtained independently of our work, established boundedness of \eqref{operator} on $L^p(\R^{n+1})$, $1 < p< \infty$, in the case $n = 2$ and $Y = \mathrm{span}_\R\{ X_1^3 , X_2^3\}$, and proposed an alternative perspective on the argument of \cite{MR3945729}.
In the present paper, we consider polynomial Carleson operators along the paraboloid in \eqref{operator} defined for general $n , d \ge 2$ and arbitrary homogeneous subspaces $Y \subseteq \mathcal{P}^{n,d}$.
We import from \cite{MR3945729} and \cite{MR4790541} their reduction argument to reduce our problem to bounding an oscillatory integral kernel defined by $Y$. Following their setup, we assume that $Y$ does not contain any nonzero linear polynomials as well as $|X|^2$.
The additional assumption of distinct degree among the basis elements of $Y$ underlying the kernel estimates in \cite{MR3945729} is no longer available in our setting. We therefore prove a new kernel estimate which serves as a homogeneous analogue of \cite[Theorem 7.1]{MR4790541}.
As a consequence, we obtain a new class of polynomial subspaces for which the polynomial Carleson operator along the paraboloid is bounded on $L^p(\R^{n+1})$, $1<p<\infty$.

\subsection{Outline of the Argument}\label{summary}

The proof of Theorem \ref{mainresult} is based on the method introduced in \cite{MR3945729}, also summarized in \cite{MR4790541}. This uses a combination of Littlewood-Paley theory, square function estimates, and the method of $TT^*$ to show that for $n \ge 2$, the $L^p(\R^{n+1})$, $1< p< \infty$ boundedness of the operator in \eqref{operator} with modulation space $Y\subseteq \R[X_1,..., X_n]$ is implied by a pointwise estimate for a particular oscillatory integral kernel defined in terms of $Y$ (see Section \ref{reductions} for a precise definition). We state the reduction result of \cite{MR3945729} in Section \ref{reductions}, and refer to \cite[Steps 1-4]{MR4790541} for its full proof. In particular, the reduction requires no other assumption on $Y$, other than the absence of linear and constant polynomials. 

The main technical contribution of our proof is a pointwise estimate of an oscillatory integral kernel defined by a codimension $n$ linear subspace $\mathcal{Q}\subseteq \mathcal{P}^{n,d}$. To illustrate the main idea more clearly, we introduce a simplified version of the kernel of interest. In the following, we denote
\begin{align*}
    N := \dim(\mathcal{Q}) = \dim(\mathcal{P}^{n,d}) - n ~,
\end{align*}
and fix a basis $\mathfrak{B} = \{p_1 , ..., p_N \}$ of $\mathcal{Q}$. Denote
\begin{align*}
    P_\lambda = \sum_{i=1}^{N} \lambda_i p_i~,~ \mbox{ for } \lambda \in \R^N~.
\end{align*}
Let $\nu , \mu \in \R^N$ and $(u,\tau) \in \R^n \times \R$. The required kernel is roughly of the form
\begin{align*}
  K_{\sharp}^{\nu , \mu}(u , \tau) := \int_{\R^{n-1}} e^{i(P_\nu(u+z) -P_\mu(z))} \Psi(u , z) ~\dd \sigma~,
\end{align*}
where $\Psi \in \CC^{\infty}_c(\R^n \times \R^n)$ is a smooth cutoff function, and $\sigma \in \R^{n-1}$ is defined implicitly by a particular change of variables $z \mapsto (\sigma, \tau)$ in $\R^{n-1}\times \R$ determined by the paraboloid. When the parameters $\nu , \mu \in \R^N$ satisfy the size condition $\abs{\nu}, \abs{\mu} \sim r$ for $r \ge 1$, we require $K^{\nu , \mu}_\sharp(u , \tau)$ to satisfy a corresponding pointwise upper bound of the form
\begin{align*}
    |K^{\nu , \mu}_\sharp(u , \tau)| \le C r^{-\delta} \quad \text{for some } \delta > 0~.
\end{align*}
We could allow the above pointwise upper bound to fail on an exceptional set of $(u,\tau)$ of small measure (more precisely, of measure $\le C r^{-\delta}$), but we require the measure of this exceptional set to be small \textbf{uniformly} in $\mu \in \R^N$. We refer the reader to Remark \ref{stoppingtime} below for a brief explanation of why this uniformity is needed in a $TT^*$ argument; see also \cite[Pages 10-11]{MR3945729}.

Implicitly, we consider the phase $P_\nu(u+z) -P_\mu(z)$ of $K_\sharp^{\nu , \mu}$ as a polynomial in $\sigma \in \R^{n-1}$. For $\beta \in \Z_{\ge 0}^{n-1}$, denote the coefficient of $\sigma^\beta$ by $C[\sigma^\beta](u, \tau)$. Then,
\begin{align}
  P_\nu(u + z) - P_\mu(z) = \sum_{0 \le \abs{\beta} \le d} C[\sigma^\beta](u, \tau) \sigma^\beta~. \label{poly}
\end{align}
Observe that the kernel $K_{\sharp}^{\nu , \mu}(u , \tau)$ is an $(n-1)$ dimensional integral, in contrast to the $n$ dimensional integral kernel
\begin{align*}
    K_\flat^{\nu , \mu} (u) = \int_{\R^n} e^{i(P_\nu(u+z) -P_\mu(z))} \Psi(u , z) ~\dd z
\end{align*}
considered in the original $TT^*$ method of Stein and Wainger \cite{MR1879821}. This reflects the phenomenon that phase polynomials from $\mathcal{Q}$ will interact implicitly with the paraboloid through the change of variables $z \mapsto (\tau , \sigma)$, and affect the decay of $K_\sharp^{\nu , \mu}(u, \tau)$.

To prove the required pointwise estimate for $K_{\sharp}^{\nu , \mu}(u , \tau)$, by a Van der Corput estimate, we need to show that for $(u,\tau)$ outside of an exceptional set, there is a multi-index $\beta \neq 0$ which is ``large'' in the sense that $\abs{C[\sigma^\beta](u, \tau)} > r^\epsilon$, where $\epsilon> 0$ is a fixed small constant. To this end, for fixed $(u,\tau)$ and $r\le |\nu| , |\mu| \le 2r$, we divide into the following cases based on the size of the coefficients $C[\sigma^\beta](u, \tau)$ for $1 \le |\beta| \le d$. 
\begin{enumerate}[label=\textbf{Case \arabic*:}, ref=\textbf{Case \arabic*:}]
  \item $\abs{C[\sigma^\beta](u, \tau)} > r^\epsilon$ for some $1 \le |\beta| \le d$ 
  \item $\abs{C[\sigma^\beta](u, \tau)} \le r^\epsilon$ for all $1 \le |\beta| \le d$
  \end{enumerate}

If $(u, \tau)$ satisfies Case 1, then Van der Corput estimates apply immediately to extract decay in $r$ from $K_{\sharp}^{\nu , \mu}(u ,\tau)$. The difficulty is to show that the set of $(u , \tau)$ satisfying Case 2 has small measure uniformly as $\mu \in \R^N$ varies. We achieve the exceptional set estimate in two steps: first, for each $(u , \tau)$ in Case 2, find an explicit $\mu$-independent representation of $C[\sigma^\beta](u, \tau)$ for some multi-indices $\beta$; then, estimate the set of $(u , \tau)$ in Case 2 using this representation.

We look for $\mu$-independent representations of the coefficients  $C[\sigma_s](u, \tau)$, $1 \le s \le n-1$. Expanding the change of variables $z \mapsto (\sigma , \tau)$ explicitly, we obtain that every coefficient $C[\sigma_s](u , \tau)$ depends linearly on both $\nu$ and $\mu$. We denote
\begin{align}
  \sum_{s = 1}^{n-1} C[\sigma_s](u , \tau)\sigma_s = L(u,\tau)P_\nu(\sigma) + \tilde{L}(u,\tau) (P_\nu - P_\mu)(\sigma) ~,\label{findsimplerrep}
\end{align}
where
\begin{align*}
    L(u,\tau), \tilde{L}(u,\tau): \mathcal{P}^{n,d}\to  \mathrm{span}_\R\{\sigma_1, \dots, \sigma_{n-1}\}
\end{align*}
are linear maps between polynomial subspaces which depend only on $u,\tau$ and the basis $\mathfrak{B}$ of $\mathcal{Q}$, but not otherwise on $\nu,\mu$. Notice that the only dependence on $\mu$ in \eqref{findsimplerrep} is in the term $\tilde{L}(u,\tau)(P_\nu - P_\mu)$.
We aim to use the assumption that $(u,\tau)$ is in Case $2$ to eliminate the dependence on $\mu$ in this term. To do so, using the change of variables $z \mapsto (\sigma, \tau)$, we construct linear maps which capture the higher order $\sigma$ terms
\begin{align}
  A(u,\tau) ~,~ \tilde{A}(u,\tau): \mathcal{P}^{n , d} \to \mathrm{span}_\R\{\sigma^\beta~~|~~ 2 \le \abs{\beta} \le d\}~,\label{restrict}
\end{align}
so that
\begin{align*}
  \sum_{2 \le |\beta| \le d} C[\sigma^\beta](u, \tau) \sigma^\beta = A(u,\tau)(P_\nu - P_\mu)(\sigma) + \tilde{A}(u,\tau) P_\nu(\sigma)~.
\end{align*}
Again, $A(u,\tau) , \tilde{A}(u,\tau)$ depend only on $u,\tau$ and the basis $\mathfrak{B}$ of $\mathcal{Q}$, but not otherwise on $\nu,\mu$. In Theorem \ref{mainresult}, the tangential nondegeneracy condition \eqref{cond1} ensures that $A(u,\tau)$ is injective on $\mathcal{Q}$ for almost every $u,\tau$, while the codimension $n$ assumption on $\mathcal{Q}$ ensures a bijection. We can thus invert $A(u,\tau)|_{\mathcal{Q}}$ on its image using Cramer's rule, so that
\begin{align*}
    P_\nu - P_\mu = (A(u,\tau)|_{\mathcal{Q}})^{-1} \bracket{\sum_{2 \le |\beta| \le d} C[\sigma^\beta](u, \tau) \sigma^\beta - \tilde{A}(u,\tau) P_\nu}~.
\end{align*}
If $(u , \tau)$ is in Case $2$, then $\abs{C[\sigma^\beta](u, \tau)} \le r^\epsilon$ for $2 \le |\beta| \le d$, so
\begin{align*}
    P_\nu - P_\mu = \text{error term bounded by $Cr^{\epsilon}$} - (A(u,\tau)|_{\mathcal{Q}})^{-1} \tilde{A}(u,\tau) P_\nu~.
\end{align*}
We then substitute $P_\nu - P_\mu$ back into \eqref{findsimplerrep}, and
represent each $C[\sigma_s](u , \tau)$, $1 \le s \le n-1$ roughly in the following form: 
\begin{align}
  C[\sigma_s](u , \tau) = \text{error term bounded by $Cr^{\epsilon}$} + L_s^{(1)}P_\nu(u) + \tau L_s^{(2)}P_\nu(u)~, \label{expressionofcoeff}
\end{align}
where $L^{(1)}_s, L^{(2)}_s: \mathcal{P}^{n,d}\to \mathcal{P}^{n,d}$ 
are linear maps which depend only on the basis $\mathfrak{B} $ of $\mathcal{Q}$, but not otherwise on $\nu, \mu , u, \tau$. We note that the precise substitution step in our proof is more technical than just inverting the operator $A(u,\tau)|_{\mathcal{Q}}$, since we need to ensure that all the relevant terms remain polynomials in $u , \tau$. 

Finally, we take the linear operators $L^{(1)}_s, L^{(2)}_s$ in \eqref{expressionofcoeff}, and use Van der Corput sublevel set estimates to show that the set of $(u,\tau)$ in Case 2 is small uniformly in $\mu$ if the linear map
\begin{align}
  \mathcal{P}^{n , d} & \longrightarrow \bracket{\mathcal{P}^{n , d}}^{2n - 2} \nonumber\\
  p & \longmapsto (L^{(1)}_1 p ~,~ ...~,~ L^{(1)}_{n-1} p ~,~ L^{(2)}_1 p ~,~ ...~,~ L^{(2)}_{n-1} p )~, \label{jointlinearmap}
\end{align}
where the codomain is taken to be the direct sum of $2n-2$ copies of $\mathcal{P}^{n,d}$,
is injective on $\mathcal{P}^{n , d}$. This is the central new idea of the paper parallel to the argument based on distinct degree assumptions in \cite{MR3945729,MR4790541}. The strategy we use to prove injectivity of \eqref{jointlinearmap} is as follows.
We will show that when $d$ is odd, the kernel of $L^{(1)} = ( L^{(1)}_1,..., L^{(1)}_{n-1})$ is trivial, so injectivity is automatic. When $d$ is even, the kernel of $L^{(1)}$ is nontrivial, so we require the radial nondegeneracy condition \eqref{cond2} to show that $L^{(2)} = (L^{(2)}_1,..., L^{(2)}_{n-1})$ is injective on the kernel of $L^{(1)}$.

In the case when $d$ is even, if the tangential nondegeneracy condition \eqref{cond1} is not satisfied, then $A(u,\tau)|_{\mathcal{Q}}$ in \eqref{restrict} may not be invertible. Then, there will be no available simple representation for the linear coefficients $C[\sigma_s](u , \tau)$ in \eqref{expressionofcoeff}, and we would not be able to proceed with the same argument. 
If the radial nondegeneracy condition \eqref{cond2} is not satisfied, then the linear map \eqref{jointlinearmap} is not guaranteed to be injective. In the worst case scenario, it is possible that there is some $|\nu| \sim r$ such that $P_\nu$ is in the kernel of $(L^{(1)}, L^{(2)})$. In this case, \textbf{all} $(u , \tau)$ will end up in Case 2, and the exceptional set would fail to be small.
A similar obstruction appears in \cite{MR3945729,MR4790541}, where the authors require that their homogeneous quadratic part $p_2$ to not be a multiple of $|X|^2$. Thus, the radial nondegeneracy condition \eqref{cond2} is a higher degree analogue of this restriction.

We remark that the choice to state the nondegeneracy conditions in their current form is simply because they are the easiest to compute with. When treating Case 2, we could have chosen to substitute into a different set of coefficients $C[\sigma^{\beta(i)}](u , \tau)$ for $i = 1,...,m$, instead of the coefficients $C[\sigma_s](u,\tau)$, $s = 1, ... , n-1$. There would naturally be a different tangential nondegeneracy condition on $\mathcal{Q}$ to ensure bijectivity of $A(u,\tau)|_{\mathcal{Q}}$. After substitution, we would have obtained a different representation of the coefficients $C[\sigma^{\beta(i)}](u , \tau)$ to \eqref{expressionofcoeff}. We could then impose different radial nondegeneracy conditions to ensure injectivity of the corresponding linear map. A similar argument would potentially still work, but it appears that the computation of \eqref{expressionofcoeff} and the kernel of $(L^{(i)})$ would be substantially more difficult.

\subsection*{Structure of the paper}
In Section \ref{preliminaries}, we survey notation used throughout the paper, as well as the crucial Van der Corput estimates needed for the oscillatory integral bound later in Section \ref{oscillatoryint}. In Section \ref{reductions}, we state the reduction result of Steps 1-4 in \cite{MR4790541}, which provides a general framework to reduce the proof of Theorem \ref{mainresult} to a particular estimate of an oscillatory integral kernel. In Section \ref{oscillatoryint}, we bound the kernel assuming a few lemmas, which we then prove in Section \ref{lemma}. Sections 3-5 together finish the proof of Theorem \ref{mainresult}. In Section \ref{subspacesec}, we prove Proposition \ref{derivationofmainthm}, which derives Theorem \ref{mainresult2} from the more general Theorem \ref{mainresult}. In Section \ref{monomialsec}, we state the complete proof for Proposition \ref{monomials}, which illustrates a special case of Theorem \ref{mainresult} under the condition that $\mathcal{Q}$ admits a monomial basis.

\subsection*{Acknowledgement}
The author was partially supported during this research by the National Science Foundation under DMS-2200470. The author would like to thank Lillian Pierce for suggesting this problem, for her guidance throughout the course of the project, and for her extremely valuable comments on earlier drafts of this manuscript. The author is also grateful to HaoYun Yao for valuable discussions about the exposition and insights about the monomial examples, and Jiahui Yu for comments on how to effectively present the method. 

\subsection*{AI Disclosure}
The author has used ChatGPT 5.5 to assist some of the work of this paper. Specifically, AI was used to generate MATLAB code to explicitly compute the determinant $D(X)$ in Lemma \ref{minorexp} for all possible partitions $\{\alpha \in \Z_{\ge 0}^n ~~|~~ \abs{\alpha} = d\} = T \sqcup S$
in Proposition \ref{monomials}, for the range $2 \le n,d \le 6$. AI was instrumental to identifying a pattern in the permissible subspaces in the early stage of the project. Moreover, AI was used to survey techniques for determinant computation when the author was working on proving an example of the main result in $\R^2$ (i.e., Corollary \ref{R2exmp}) using a slightly different approach. Relevant techniques were later adapted to prove of some of the lemmas in this paper. More precisely, these techniques refer to the use of Schur's formula to compute the determinant of $\widetilde{M}_s(X)$ in the proof of Lemma \ref{L1}, as well as the isolation of the leading monomial in the determinant $\det(\partial_jf_i(X))$ in the proof of Lemma \ref{lemma:subspace}. The author has confirmed all AI-generated suggestions through manual verification, and no part of this manuscript was written by AI.

\section{Preliminaries}\label{preliminaries}
This section surveys notation used throughout the paper, as well as the Van der Corput estimates needed for later proofs.
For $A , B \in \C$, we use the notation $B = O(A)$ if there exists a positive constant $C$ such that $\abs{B} \le C\abs{A}$; the value of $C$ may change from line to line.

Let $m , k \ge 2$ be fixed integers and $\R[X_1, ... ,X_m]$ denote the class of real-valued polynomials in $m$ variables. Let
\begin{align}
  Q_\lambda(X) = \sum_{0 \le \abs{\alpha} \le k}^{} \lambda_\alpha X^\alpha~\label{Qlambda}
\end{align}
be a polynomial of degree at most $k$ in $\R[X_1,..., X_m]$, where $\lambda = (\lambda_\alpha)_{0 \le \abs{\alpha} \le k}$ is a real-valued vector. Let
\begin{align*}
  \norm{\lambda} : = \sum_{1 \le \abs{\alpha} \le k}^{} \abs{\lambda_\alpha}
\end{align*}
denote the isotropic norm of the coefficients of nonconstant terms, and let
\begin{align}
  [\![\lambda]\!]: = \sum_{0 \le \abs{\alpha} \le k}^{} \abs{\lambda_\alpha} = \norm{\lambda} + \abs{Q_\lambda(0)} \label{isotropicnorm}
\end{align}
denote the isotropic norm, including the constant term. We
will use $[\![Q_\lambda(X)]\!]_X$ and $[\![\lambda]\!]$ interchangeably. The presence of the variable $X$ in the notation $[\![\cdot]\!]_X$ here serves to emphasize that the norm sees the input as a polynomial in $X$. In particular, since every polynomial can be written in the form \eqref{Qlambda}, the function $p \mapsto [\![p(X)]\!]_X$ is a well-defined norm on the vector space $\R[X_1,...,X_m]$.

Let $\mathcal{P}^{m,k} \subseteq \R[X_1, ... ,X_m]$ denote the $\binom{k+m-1}{m-1}$ dimensional vector space of real-valued homogeneous polynomials of degree $k \ge 2$ in $m\ge 2$ variables. We later will need to consider orthogonal complements of various subspaces in $\mathcal{P}^{m,k}$. We take the inner product underlying the orthogonal projections to be the Fischer inner product: for $p(X) = \sum_{\alpha} \lambda_\alpha X^\alpha, q(X) = \sum_{\alpha} \nu_\alpha X^\alpha \in \mathcal{P}^{m , k}$, let 
\begin{align}
  \inner{p }{q} : = p(\partial X) q(0) = \sum_{\alpha}^{} \lambda_\alpha \partial^\alpha q(0)~.  \label{innerprod}
\end{align}
We point out that any inner product on $\mathcal{P}^{m,k}$ works for the proof of Theorem \ref{mainresult}. The specific choice of inner product in \eqref{innerprod} is in fact only needed for the proof of Propositions \ref{derivationofmainthm} and \ref{monomials} (both are applications of Theorem \ref{mainresult}), since we will need to be able to compute the orthogonal projection of some polynomials explicitly.

For ease of notation, in the proof of the oscillatory kernel estimate in Section \ref{oscillatoryint}, we will need to consider a norm on $\mathcal{P}^{m , k}$ different from both the isotropic norm in \eqref{isotropicnorm} and the $\ell^2$ norm associated to the inner product in \eqref{innerprod}. 
For a codimension $m$ linear subspace $\mathcal{Q}\subseteq \mathcal{P}^{m,k}$ which has 
dimension $N : = \binom{k+m-1}{m-1} - m$, equipped with a basis $\{p_1,...,p_N\}$, we take $\{p_{N+1}, ... ,p_{N+m}\}$ to be a basis for the orthogonal complement $\mathcal{Q}^\perp$ in $\mathcal{P}^{m , k}$, so that $\{p_1,...,p_N, p_{N+1} ,..., p_{N+m}\}$ is a basis for $\mathcal{P}^{m , k}$.
We equip $\mathcal{P}^{m , k}$ with the $\ell^1$ norm corresponding to this particular basis, and denote
\begin{align}
  \norm{p}_{\ell^1(\mathcal{Q})} = \sum_{i =1}^{N+m} \abs{\nu_i}~,~\mbox{ if  } ~~~p = \sum_{i = 1}^{N+m} \nu_i p_i \in \mathcal{P}^{m , k}~. \label{l1norm}
\end{align}

Next, we survey two Van der Corput estimates which will be used later in the paper. The first estimate is Proposition
2.1 of \cite{MR1879821}:
\begin{lem} \label{Vandecorput1}
Let $Q_{\lambda}(X) : = \sum_{0 \le \abs{\alpha} \le k}^{} \lambda_\alpha X^\alpha \in \R[X_1,...,X_m]$.
For any $\CC^1$ function $\psi$ defined on the unit ball
$B_1(\mathbb{R}^m)$, with $\|\psi\|_{\CC^1} \leq 1,$
and for any convex subset $\Omega \subseteq B_1(\mathbb{R}^m)$, we have
\begin{equation*}
\left|
\int_{\Omega} e^{iQ_{\lambda}(x)}\psi(x)\,dx
\right|
\leq C\|\lambda\|^{-1/k}~,
\end{equation*}
where the constant $C$ depends on the dimension $m$ and the degree $k$
of $Q_{\lambda}$, but not otherwise on $Q_{\lambda}$, $\psi$, or
$\Omega$. This upper bound is in particular independent of the constant term $Q_{\lambda}(0)$.
\end{lem}
The result of the lemma continues to hold with the unit ball $B_1(\mathbb{R}^m)$ replaced by a Euclidean ball of any other fixed radius such as $B_{C_0}(\R^m)$, with an implicit constant dependent on $C_0$.
We also need an estimate for the measure of the set where $Q_{\lambda}$ takes small values, which includes the case in which the constant term of $Q_\lambda$ is large.
\begin{lem}\label{Vandecorput2}
Let $Q_{\lambda}(X) : = \sum_{0 \le \abs{\alpha} \le k}^{} \lambda_\alpha X^\alpha \in \R[X_1,...,X_m]$. For every $\rho>0$, one has
\begin{equation*}
\left|
\left\{
x\in B_1(\mathbb{R}^m): |Q_{\lambda}(x)|\leq \rho
\right\}
\right|
\leq 
\begin{cases}
  C'\rho^{1/k}[\![\lambda]\!]^{-1/k} & \mbox{ if } k\ge 1\\
  C'\rho [\![\lambda]\!]^{-1} & \mbox{ if } k = 0
\end{cases}~,
\end{equation*}
where the constant $C'$ depends only on the dimension $m$ and the degree
$k$ of $Q_{\lambda}$.
\end{lem}
\begin{proof}
The case in which $k \ge 1$ is Lemma 3.3 of \cite{MR3945729}. We consider the case in which $Q_\lambda(x) = \lambda_0 \neq 0$ is a constant polynomial. For each $\rho > 0$, we have
\begin{align*}
  \abs{\cbracket{x\in B_1(\mathbb{R}^m): |Q_{\lambda}(x)|\leq \rho} } & = \begin{cases}
    \abs{B_1(\mathbb{R}^m)} & \mbox{ if } \inorm{\lambda} \le \rho\\
    0 & \mbox{ if } \inorm{\lambda} > \rho
  \end{cases} = \abs{B_1(\mathbb{R}^m)} \chi_{\rho \ge \inorm{\lambda}}(\rho)~.
\end{align*}
Then the result of the lemma follows from the observation that if $\lambda_0 \neq 0$, then $\chi_{\rho \ge \inorm{\lambda}}(\rho) \le \rho \inorm{\lambda}^{-1}$ as a function of $\rho > 0$.

\end{proof}

Again, the result of the above lemma continues to hold with the unit ball $B_1(\mathbb{R}^m)$ replaced by a Euclidean ball of any other fixed radius $C_0$, with an implicit constant dependent on $C_0$.

\section{Reduction Theorem} \label{reductions}
In this section, we state the key theorem that will be used to reduce the proof of Theorem \ref{mainresult} to a particular oscillatory kernel estimate. This is exactly the result of Steps 1-4 (corresponding to Sections 3-6) summarized in \cite{MR4790541}. In particular, the reduction makes no extra assumption yet on the structure of the subspace $\mathcal{Q} \subseteq \R[X_1,...,X_n]$ in the definition of the operator in \eqref{operator}, other than the requirement that $\mathcal{Q}\cap \mathrm{span}_\R\{X_1,..., X_n\} = \{0\}$.

Let $n\ge 2$. We first introduce notation required to precisely define the $K^{\nu , \mu}_{\sharp}(u,\tau)$ kernel.
For any variable $u \in \R^n$, let $u^{(l)}$ denote the variable in $\R^{n-1}$ that omits the $l$-th coordinate of $u$. Similarly, for any multi-index $\alpha \in \Z^{n}_{\ge 0}$, we let $\alpha^{(l)} \in \Z^{n-1}_{\ge 0}$ denote the multi-index omitting $\alpha_l$.
Given coordinates $(u , z) \in (\R^n \backslash \{0\} ) \times \R^n $, we fix a coordinate $1 \le l \le n$ where $u_l \neq 0$, and make the change of variables $z \mapsto (\tau, \sigma) \in \R \times \R^{n-1}$ by
\begin{align}
  \tau = \frac{u \cdot z}{\abs{u}}~,~ \quad \sigma = \frac{u^{(l)} \tau - \abs{u}z^{(l)}}{u_l}~.\label{changeofvar}
\end{align}
so that $z$ is defined implicitly by the relations
\begin{align*}
  z^{(l)}  &=  \frac{u^{(l)}\tau - u_l \sigma}{\abs{u}}~,~\quad
  z_l = \frac{u_l \tau + u^{(l)}\cdot \sigma}{\abs{u}}~.
\end{align*}
By explicit computation, the Jacobian associated with this change of variables is $(\abs{u_l} / \abs{u})^{n-2}$. 
Notice that if $(u,z) \in B_2(\R^n) \times B_1(\R^n)$, with $\abs{u_l} \ge c$ for some fixed constant $c> 0$, then $\abs{\tau} \le 1$, and
\begin{align*}
    \abs{\sigma} \le \frac{\abs{u^{(l)}} \abs{\tau} + \abs{u} \abs{z^{(l)}}}{\abs{u_l}} \le \frac{2\abs{u}}{c} \le \frac{4}{c}~.
\end{align*}
Thus compact support in the $(u, z)$ variables is preserved under this change of variables, provided that $(u, z)$ ranges over a bounded set depending on $c$ and $l$.

Next, fix a partition of unity for $\mathbb{S}^{n-1} \subseteq \R^n$, say
\begin{align}
  1 = \sum_{l = 1}^{n} W_l(s)~, \label{partition}
\end{align}
such that for each $1 \le l \le n$, $W_l \in \CC^\infty_c(\mathbb{S}^{n-1})$ and $W_l$ is supported where $\abs{s_l} \ge c_0$, for some fixed $c_0 > 0$. 

Let $\Psi (u ,z)$ be a $\CC^1$ function supported on $B_2(\R^n) \times B_1(\R^n)$, satisfying the bound \begin{align*}
  \norm{\Psi}_{\CC^1}: = \sup_{(u,z) \in B_2(\R^n) \times B_1(\R^{n})} \bracket{\abs{\Psi(u, z)} + \abs{\frac{\partial}{\partial \sigma} \Psi(u,z)}}~\le~ 1~.
\end{align*}
Let $\{p_1, ..., p_N\}$ be a set of linearly independent homogeneous polynomials of degree at least $2$. For $\nu , \mu \in \R^N$, define
\begin{align*}
  P_{\nu} : =\sum_{i=1}^N\nu_i p_i ~,~ P_{\mu} : =\sum_{i=1}^N\mu_i p_i \in \R[X_1,...,X_n]~.
\end{align*}
Define the kernel $K_{\sharp, l}^{\nu, \mu}(u,\tau)$ by 
\begin{align}
  K_{\sharp, l}^{\nu, \mu}(u, \tau) := \bracket{\frac{\abs{u_l}}{\abs{u}}}^{n-2} W_l \bracket{ \frac{u}{\abs{u}}} \int_{\R^{n-1}} e^{i\bracket{P_\nu(u + z) - P_\mu(z)}} \Psi(u, z) \dd \sigma~, \label{kernel}
\end{align}
where the $z$ in the integral is defined implicitly via the change of variables $z \mapsto (\tau, \sigma) \in \R \times \R^{n-1}$. Observe that $K_{\sharp, l}^{\nu, \mu}(u, \tau)$ is supported on $B_2(\R^n) \times B_1(\R)$ as a function of $(u , \tau)$ since $\Psi (u ,z)$ is a function supported on $B_2(\R^n) \times B_1(\R^n)$. Moreover, the region of integration in $\sigma$ is compact, because the change of variables preserves compact support if $u/|u|$ lies in the support of $W_l$. This is important because it allows Van der Corput estimates to be applied effectively to the integral in $\sigma$.

With the definition of the $K_{\sharp,l}^{\nu, \mu}$ kernel stated, we are now ready to state the reduction theorem (Theorem \ref{reductionthm} below). This is exactly Theorem 6.1 of \cite{MR4790541} and we refer the reader to \cite{MR4790541} for a complete proof.

\begin{thm}
\label{reductionthm}
For each $1\le l\le n$, for $u \in \R^n$ such that $u/|u|$ lies in the support of the partition $W_l$ defined in \eqref{partition}, consider the kernel $K_{\sharp,l}^{\nu , \mu}(u,\tau)$ as defined in \eqref{kernel}. Suppose there exists a constant $\delta >0$, such that for all $r\ge 1$ and $\mu , \nu \in\R^N$ satisfying
$r \le \norm{\nu} , \norm{\mu} \le 2 r$,
there exists a measurable set $G^\nu \subseteq B_2(\R^n)$ (depending on $\nu$ but independent
of $r, \mu$) and for each $u \in B_2(\R^n)$, a measurable set $F_u^\nu \subseteq B_1(\R)$ (depending on $\nu , u$ but independent
of $ r, \mu$) such that
\begin{align}
  \abs{G^\nu} \le C r^{-\delta}~,~ \abs{F_u^\nu} \le Cr^{-\delta} \mbox{ for all } u \in B_2(\R^n)~,\label{measureR2}
\end{align}
and
\begin{align}
  \abs{K_{\sharp,l}^{\nu, \mu}(u, \tau)} \le C \bracket{r^{-\delta} \chi_{B_2}(u) \chi_{B_1}(\tau) + \chi_{G^\nu}(u) \chi_{B_1}(\tau) + \chi_{B_2}(u) \chi_{F_u^\nu}(\tau)}~.\label{KsharpestimateR2}
\end{align}
Then, the maximally modulated operator 
\begin{align*}
  f \mapsto \sup_{\lambda \in \R^N} \abs{R_{P_\lambda}f}
\end{align*}
satisfies an \textit{a priori bound} on $L^p(\R^{n+1})$, for each $1 < p < \infty$.
\end{thm}

\begin{rem}\label{stoppingtime}
  Notice that the exceptional sets $G^\nu$ and $F_u^\nu$ are independent of the parameter $\mu$. Below, we briefly explain how this independence is crucial for a $TT^*$ argument. This is one of the fundamental steps in the reduction argument which derives Theorem \ref{reductionthm}.

We first introduce a model case to illustrate how a $TT^*$ argument is applied to a maximally modulated operator.
Let $m\ge 1$, $d \ge 2$ be integers and let $\alpha \in \Z_{\ge 0}^m$ denote multi-indices in $m$ variables. Let $M\ge 1$ be the integer such that
\begin{align*}
    \{ (\lambda_\alpha)_{2 \le |\alpha| \le d}~~|~~ \lambda_\alpha \in \R\} = \R^M~.
\end{align*} 
A central step in the original method of Stein and Wainger \cite{MR1879821} considers a maximal operator of the form
\begin{align}
    f \mapsto \sup_{ k \in \Z} \sup_{r \le |\lambda | \le 2r} |I_{\lambda , 2^k} f |~,\label{Ioperator}
\end{align}
where for each $a \in \Z$, $\lambda \in \R^M$, $I_{\lambda, a}$ is a particular linear operator on $\SSS(\R^m)$.
Stein and Wainger showed that $L^p(\R^m)$, $1 < p< \infty$, boundedness of the Carleson operator in Theorem \ref{steinwaingerthm} follows from an $L^2$ bound of \eqref{Ioperator} with operator norm which decays as $r \to + \infty$.

Fix $r \ge 1$.
Take arbitrary measurable functions $a: \R^m \to \Z$, $\lambda: \R^m \to \R^M$ such that $r \le |\lambda | \le 2r$. To obtain the desired bound of the maximal operator \eqref{Ioperator}, it suffices to bound the linearized operator
\begin{align}
  Tf(x) = I_{\lambda(x), a(x)}f(x)~, \label{linearize}
\end{align}
on $L^2$ with operator norm independent of the functions $a , \lambda$ which also decays in $r$. Indeed, if this is true, then we can construct functions $a(x) , \lambda(x)$ such that for each $x$,
\begin{align*}
    |Tf(x)| \ge \frac{1}{2}\sup_{ k \in \Z} \sup_{r \le |\lambda | \le 2r} |I_{\lambda , 2^k} f (x)|~~,
\end{align*}
so that $T$ dominates the maximal operator \eqref{Ioperator}.
We call the functions $a , \lambda$ \textbf{stopping times}. 

To bound $T$, it suffices to bound $TT^*$, where $T^*$ is the adjoint of $T$. The advantage of bounding $TT^*$ is that it can be written as integration agains a kernel which we can now analyze. For the model case operator, computing $TT^*$ of \eqref{linearize} gives that
\begin{align}
    TT^*f (x)= \int_{\R^m} f(x - y) K^{\lambda(x) , \lambda(x- y )}_{a(x) , a(x - y)}(y) ~\dd y~,\label{Tkernel}
\end{align}
where $K^{\nu , \mu}_{a_1 , a_2}(y)$ is a particular oscillatory integral kernel of the form
\begin{align*}
    K^{\nu , \mu}_{a_1 , a_2}(y) = a_1^{-m}a_2^{-m}\int_{\R^m} e^{i(P_\nu((y+z)/a_1) -P_\mu(z/a_2))} \Psi((y+z)/a_1)\Psi(z/a_2) ~\dd z~,~ \mbox{ for } \Psi \in \CC^1_c(\R^m)~,
\end{align*}
and $P_\lambda(y): = \sum_{2 \le |\alpha| \le d} \lambda_\alpha y^\alpha$.
Stein and Wainger proved that once one shows that for every $r\ge1$ and $r \le ||\nu|| , ||\mu|| \le 2r$, there exist small sets $E^\nu , E^\mu \subseteq B_2(\R^m)$ such that $|E^\nu| , |E^\mu| \le Cr^{-\delta}$ and a pointwise upper bound of the form
\begin{align}
  |K^{\nu , \mu}_{a_1 , a_2}(x)| \le Cr^{-\delta}(a_1^{-m}\chi_{B_2}(x/a_1) + a_2^{-m}\chi_{B_2}(x/a_2)) + C( a_1^{-m}\chi_{E^\nu}(x/a_1) + a_2^{-m}\chi_{E^\mu}(x/a_2))~,\label{Kbound}
\end{align}
then an $L^2$ bound for $TT^*$ of the form
\begin{align*}
    ||TT^* f||_{L^2(\R^m)} \le Cr^{-\delta} ||f||_{L^2(\R^m)}
\end{align*}
can be obtained by comparison to certain maximal functions. Here, $\delta > 0$ is an absolute constant which only depends on $m$ and $d$. The reason why Stein and Wainger needed exceptional sets to depend separately on $\nu$ and $\mu$  is because a key step in their proof is a bilinear argument, which views the value $||TT^*f ||_{L^2(\R^m)}$ as the supremum over all $g \in \SSS(\R^m)$ with $||g||_{L^2(\R^m)} = 1$ of
\begin{align*}
  \abs{\int_{\R^m} TT^* f(x) \bar{g} (x) \dd x} = \abs{\int_{\R^m}\int_{\R^m} \bar{g}(x) f(x - y) K^{\lambda(x) , \lambda(x - y)}_{a(x) , a(x - y)}( y) ~\dd x \dd y } ~.
\end{align*}
Applying the pointwise upper bound of $K^{\nu , \mu}_{a , b}$ to the above, we obtain
\begin{align*}
    |\int_{\R^m}& TT^* f(x)  \bar{g} (x) \dd x | \\
    \le & ~ C r^{-\delta}\int_{\R^m}\int_{\R^m} |\bar{g}|(x) |f|(x - y) \bracket{a(x)^{-m}\chi_{B_2}(y/a(x)) + a(x-y)^{-m}\chi_{B_2}(y/a(x-y))} ~ \dd y\dd x\\ 
     & + C \int_{\R^m} \int_{\R^m} |\bar{g}|(x) |f|(x - y) \bracket{ a(x)^{-m} \chi_{E^{\lambda(x)}}(y/a(x)) + a(x-y)^{-m}\chi_{E^{\lambda(x - y)}}(y/ a(x-y)) } ~ \dd y\dd x~.
\end{align*}
We can then use a linear change of variables $x \mapsto (x+ y)$ and exchange the order of integration, so that
\begin{align*}
   |\int_{\R^m} &  TT^* f(x)  \bar{g} (x) \dd x |\\ \le &  C r^{-\delta}\int_{\R^m} |\bar{g}|(x) \bracket{a(x)^{-m}\int_{\R^m}  |f|(x - y) \chi_{B_2}(y/a(x))~\dd y} \dd x \\ & + C \int_{\R^m} |\bar{g}|(x) \bracket{a(x)^{-m}\int_{\R^m} |f|(x - y) \chi_{E^{\lambda(x)}}(y/a(x)) ~\dd y} \dd x\\
   &+ C r^{-\delta}\int_{\R^m} |f|(y) \bracket{a(y)^{-m}\int_{\R^m}  |g|(x + y) \chi_{B_2}(x/a(y))~\dd x} \dd y \\ & + C \int_{\R^m} |f|(y) \bracket{ a(y)^{-m}\int_{\R^m} |g|(x + y)\chi_{E^{\lambda(y)}}(x/a(y))~\dd x} \dd y~.
\end{align*}
By Hölder's inequality, we get
\begin{align*}
    \abs{\int_{\R^m} TT^* f(x) \bar{g} (x) \dd x} \le & ~ C r^{-\delta}\norm{g}_{L^2(\R^m)} \norm{\mathcal{M} f }_{L^2(\R^m)} + \norm{\mathcal{M} g}_{L^2(\R^m)} \norm{f }_{L^2(\R^m)} \\ 
    & + C \norm{g}_{L^2(\R^m)} \norm{\mathcal{M}_r f }_{L^2(\R^m)}\norm{\mathcal{M}_r g}_{L^2(\R^m)} \norm{f }_{L^2(\R^m)}~.
\end{align*}
Here, $\mathcal{M}$ is the Hardy-Littlewood maximal operator, and $\mathcal{M}_r$ is a bounded maximal operator with operator norm decaying in $r$. If the exceptional sets in the estimate of $K_{a_1 , a_2}^{\nu , \mu}$ did not depend separately on $\nu$ and $\mu$, then one could not consider the inner integrals as bounded maximal functions.

Analogously, the kernel $K_{\sharp,l}^{\nu, \mu}(u,\tau)$ defined in \eqref{kernel} is an integration kernel of a $TT^*$ operator similar to \eqref{Tkernel}, but now one also has to take into account Radon behavior introduced by integration over the paraboloid. Hence, exceptional sets cannot depend on both $\nu$ and $\mu$ for the exact same reason as that of the model case. 
Let $n , d \ge 2$,
$N = \dim(\mathcal{Q})$, $r \ge 1$ and $\lambda : \R^{n+1} \to \R^N$ be an arbitrary stopping time which satisfies $r \le |\lambda| \le 2r$. For the restricted setting of taking $a = 1$, the analog of \eqref{Tkernel} in the present paper takes the form
\begin{align*}
    II^*f(x , t) : = \int_{\R^n} \int_{\tau \in \R} f(x - u , t - |u|^2 - 2|u|\tau)  K_{\sharp,l}^{\lambda(x,t), \lambda(x - u , t - |u|^2 - 2|u|\tau)}(u,\tau) ~\dd u \dd \tau~,
\end{align*} 
where $K_{\sharp,l}^{\nu, \mu}(u,\tau)$ is the kernel given by \eqref{kernel}.
Suppose now for each $r\ge 1$ we had a small set $H^{\nu, \mu} \subseteq B_{2}(\R^{n})\times B_{1}(\R)$, such that $\abs{H^{\nu, \mu}} \le Cr^{-\delta}$ and
\begin{align*}
  \abs{K_{\sharp,l}^{\nu, \mu}(u,\tau)} \le C (r^{-\delta } \chi_{B_2}(u) \chi_{B_1}(\tau) + \chi_{H^{\nu ,\mu}}(u,\tau))~.
\end{align*}
This bound is achievable through a direct application of the Van der Corput lemmas introduced in Section \ref{preliminaries} to $K_{\sharp,l}^{\nu, \mu}$. However, the dependence of $H^{\nu ,\mu}$ on both $\nu$ and $\mu$ is problematic if we later want to apply Stein and Wainger's maximal function argument.
Applying this pointwise bound gives
\begin{align*}
    |\int_{\R^{n+1}} & II^*f(x , t) \bar{g}(x,t) ~\dd x \dd t |\\
    \le & ~ C r^{-\delta} \int_{\R^{2n+2}} |f|(x - u , t - |u|^2 - 2|u|\tau)|\bar{g}|(x,t)\chi_{B_2}(u) \chi_{B_1}(\tau)~ \dd u \dd \tau \dd x \dd t \\ 
    &+ C\int_{\R^{2n+2}} |f| (x - u , t - |u|^2 - 2|u|\tau)|\bar{g}|(x,t) \chi_{H^{\lambda(x,t), \lambda(x-u, t-|u|^2-2|u|\tau)}}(u,\tau) ~\dd u \dd \tau \dd x \dd t~.
\end{align*}
Upon setting the linear (with respect to $(x,t)$) change of variables $(z,s) = (x-u, t - |u|^2 - 2|u|\tau)$, the second term becomes
\begin{align*}
    C \int_{\R^{n+1}}|f|(z,s) \bracket{\int_{\R^{n+1}} |\bar{g}|(z+u,s+|u|^2+2|u|\tau) \chi_{H^{\lambda(z+u,s+|u|^2+2|u|\tau), \lambda(z,s)}}(u,\tau) ~\dd u \dd \tau} \dd z \dd s~.
\end{align*}
We cannot expect the inner integral here to be dominated by a bounded maximal operator with operator norm which decays in $r$, since now the small set $H^{\lambda(z+u,s+|u|^2+2|u|\tau), \lambda(z,s)}$ depends also on the integration variables: the size of this set is small for every fixed $z ,s,u,\tau$, but the integral over $u$ and $\tau$ could still be large.

Therefore, we have shown that the dependence of exceptional sets is what obstructs us from being able to directly apply the original method of $TT^*$ of Stein and Wainger to an operator which involves integration over a curved surface. To circumvent this issue, in the reduction argument introduced by \cite{MR3945729}, 
Pierce and Yung used a combination of $TT^*$, square functions, and sophisticated stopping time arguments to prove that
once a $\mu$-\textbf{independent} bound proposed by Theorem \ref{reductionthm} is obtained for $K_{\sharp,l}^{\nu, \mu}(u,\tau)$, then one can obtain an $L^p$, $1 < p< \infty$ bound for the operator \eqref{operator}. This is why we emphasize the importance of independence in $\mu$ in our oscillatory integral bound. See \cite[Pages 10-11]{MR3945729} for a more detailed explanation.
\end{rem}

\section{Oscillatory Integral Bound}\label{oscillatoryint}

To prove Theorem \ref{mainresult}, by Theorem \ref{reductionthm}, it suffices to verify when $\mathcal{Q} \subseteq \mathcal{P}^{n, d}$ is a codimension $n$ subspace, the associated kernel $K_{\sharp, l}^{\nu, \mu}$ satisfies a suitable upper bound. We record the desired result below.
\begin{prop}\label{mainprop}
Fix a dimension $n \ge 2$ and a degree $d \ge 2$. Let $\mathcal{Q}\subseteq \mathcal{P}^{n,d}$ be a
\begin{align*}
 N : = \binom{d+n-1}{n-1} - n
\end{align*}
dimensional linear subspace, equipped with basis $\mathfrak{B} =\{p_1, ... ,p_N\}$.
Suppose there exists $c \in (\R\backslash\{0\})^n$ 
such that $\mathcal{Q}$ satisfies the tangential nondegeneracy condition for $c$, i.e.,
\begin{align*}
  \mathcal{Q} \cap \mathrm{span}_\R \cbracket{ X_j(c \cdot X)^{d-1}~~|~~ 1 \le j \le n} = \{0\}~.
\end{align*}
If $d$ is even, assume further that $\mathcal{Q}$ also satisfies the radial nondegeneracy condition \eqref{cond2} for $c$, i.e.,
\begin{align*}
  \mathcal{Q} \cap \mathrm{span}_\R \cbracket{ (c \cdot X)^d~,~ \abs{X}^{2}(c \cdot X)^{d-2}} = \{0\}~.
\end{align*}

For each $1\le l\le n$, for $u \in \R^n$ such that $u/|u|$ lies in the support of the partition $W_l$ defined in \eqref{partition}, consider the kernel $K_{\sharp,l}^{\nu , \mu}(u,\tau)$ as defined in \eqref{kernel}, where $\{p_1,...,p_N\}$ is taken to be the chosen basis of $\mathcal{Q}$.
Then, there exists a constant $\delta >0$ (depending only on $d$ and $n$), such that the following holds: if $r\ge 1$ and $\mu , \nu \in \R^N$ satisfy $ r \le \norm{\nu} , \norm{\mu} \le 2r$,
then there exists a measurable set $G^\nu \subseteq B_2(\R^n)$ and for each $u \in B_2(\R^n)$, a measurable set $F_u^\nu \subseteq B_1(\R)$ such that
\begin{align}
  \abs{G^\nu} \le C r^{-\delta}~,~ \abs{F_u^\nu} \le Cr^{-\delta} \mbox{ for all } u \in B_2(\R^n)~, \label{setestimate}
\end{align}
and
\begin{align}
  \abs{K_{\sharp,l}^{\nu, \mu}(u, \tau)} \le C \bracket{r^{-\delta} \chi_{B_2}(u) \chi_{B_1}(\tau) + \chi_{G^\nu}(u) \chi_{B_1}(\tau) + \chi_{B_2}(u) \chi_{F_u^\nu}(\tau)}~. \label{fullkernelbound}
\end{align}
The choices of the sets $G^\nu$ and $F_u^\nu$ are independent of both $\mu$ and $\Psi$, and the choice of $C$ depends only on $d, n$ and the chosen basis $\{p_1 , ..., p_N\}$.
\end{prop}

The rest of this section is dedicated to the proof of Proposition \ref{mainprop}.
For ease of notation, we only prove the proposition for the case $l = n$, and remark that the argument is analogous for each $1 \le l\le n$.

\subsection{Notation}
We first state the notation required for the proof of Proposition \ref{mainprop}. First, we substitute the change of variables $z \mapsto (\sigma , \tau)$ defined in \eqref{changeofvar} into the phase polynomial $P_\nu(u + z) - P_\mu(z)$ in $K_{\sharp , n}^{\nu , \mu}$, and consider it as a polynomial of $\sigma = (\sigma_1 , ... ,\sigma_{n-1})$. We adopt the notation
\begin{align}
P_\nu(u + z) - P_\mu(z) = \sum_{0 \le \abs{\beta} \le d} C[\sigma^\beta](u, \tau) \sigma^\beta~, \nonumber
\end{align}
so that $C[\sigma^\beta](u, \tau)$ denotes the coefficient of the monomial $\sigma^\beta$, and is itself a function of $u,\tau$. In the rest of this paper, we will also frequently use the notation $[\sigma^\beta](f(\sigma))$ to denote the coefficient of $\sigma^\beta$ in a polynomial $f(\sigma)$.

We now compute the coefficients $C[\sigma^\beta](u, \tau)$ explicitly. If $X = (X_1 , ..., X_n)$, we use $X^{(n)} = (X_1 , ..., X_{n-1})$ to denote the first $n-1$ components of $X$.
For $\sigma \in \R^{n-1}$, we use $X \mapsto \omega(\sigma , X)$ to denote the following change of variables:
\begin{align}
  \omega(\sigma , X)^{(n)} = X^{(n)} - X_n \sigma~,~ \omega(\sigma , X )_n = X_n + X^{(n)} \cdot \sigma~. \label{changeofvar2}
\end{align}
Using the binomial theorem to expand all the powers, we obtain that if $T \in \R$, $x \in \mathbb{S}^{n-1}$, $\sigma \in \R^{n-1}$, $\alpha \in \Z_{\ge 0}^n$ and $\beta \in \Z_{\ge 0 }^{n-1}$, then
\begin{align*}
  (Tx^{(n)} - x_n \sigma)^{\alpha^{(n)}}(Tx_n + x^{(n)}\cdot \sigma)^{\alpha_n} = \sum_{0 \le \abs{\beta} \le \abs{\alpha}}^{} T^{\abs{\alpha} - \abs{\beta}}[\sigma^\beta]\omega(\sigma , x )^{\alpha} \sigma^\beta~,
\end{align*}
where the sum ranges over all multi-indices $\beta \in \Z_{\ge 0}^{n-1}$ satisfying $0 \le \abs{\beta} \le \abs{\alpha}$.
Then, if $p = \sum_{\abs{\alpha} = d} \lambda_\alpha X^\alpha$, then since all monomials of $p$ are of degree $d$, we have
\begin{align*}
  p(Tx^{(n)} - x_n \sigma&, Tx_n + x^{(n)} \cdot \sigma) \\ &= \sum_{\abs{\alpha} = d} \lambda_\alpha (Tx^{(n)} - x_n \sigma)^{\alpha^{(n)}}(Tx_n + x^{(n)}\cdot \sigma)^{\alpha_n} = \sum_{0 \le \abs{\beta} \le d}^{} T^{d - \abs{\beta}}[\sigma^\beta]p(\omega(\sigma , x)) \sigma^\beta.
\end{align*}
Upon setting $T_1 = \abs{u} + \tau$, $T_2 = \tau$, and $x = u/\abs{u} \in \mathbb{S}^{n-1}$, according to the definition of the change of variables \eqref{changeofvar}, we can compute \begin{align*}
  P_{\nu}(u+z) - P_\mu(z) &= P_{\nu}\bracket{ T_1 x^{(n)} - x_n \sigma , T_1 x_n + x^{(n)} \cdot \sigma}  - P_\mu\bracket{ T_2 x^{(n)} - x_n \sigma , T_2 x_n + x^{(n)} \cdot \sigma} \\
  & = \sum_{0 \le \abs{\beta} \le d}^{} T_1^{d - \abs{\beta}}[\sigma^\beta]P_\nu(\omega(\sigma , x))\sigma^\beta -  T_2^{d - \abs{\beta}}[\sigma^\beta]P_\mu(\omega(\sigma , x))\sigma^\beta\\
   &= \sum_{0 \le \abs{\beta} \le d} \bracket{\sum_{i = 1}^{N} \bracket{T_1^{d - \abs{\beta}} \nu_i - T_2^{d - \abs{\beta}} \mu_i}[\sigma^\beta]p_i(\omega(\sigma , x))}\sigma^\beta~.
\end{align*}
Therefore, for each $0 \le \abs{\beta} \le d$, we obtain the explicit formula for the coefficients:
\begin{align}
    C[\sigma^\beta](u, \tau) = \tau^{d-\abs{\beta}}\sum_{i = 1}^N \bracket{\nu_i - \mu_i} B_{i , \beta}\bracket{\frac{u}{\abs{u}}} + \bracket{(\abs{u} + \tau)^{d-\abs{\beta}} - \tau^{d-\abs{\beta}}} \sum_{i = 1}^N \nu_i B_{i , \beta}\bracket{\frac{u}{\abs{u}}}~, \label{coeffsigma}
\end{align} 
where
\begin{align}
    B_{i , \beta}(X) : = [\sigma^\beta]p_i(\omega(\sigma , X))~. \label{Bdef}
\end{align}
In particular, $B_{i , \beta}(X) \in \mathcal{P}^{n , d}$, since each coordinate of $\omega$ lies in $\mathcal{P}^{n,1}$ as a polynomial of $X$. Moreover, we have $B_{i , \beta}(u / \abs{u})= O(1)$ for all $u \neq 0$, since the polynomials $B_{i , \beta}$, are entirely determined by the fixed basis $\{p_1,... , p_N\}$ of $\mathcal{Q}$.


The explicit expression in \eqref{coeffsigma} reveals that the coefficient $C[\sigma^\beta](u, \tau)$ depends linearly on the vectors $\nu , \mu$. Inspired by this observation and the matrix representation in \cite[Equation (7.9)]{MR4790541}, we now package the coefficients corresponding to a selected set of multi-indices into a matrix form.
Given an arbitrary set of $M$ multi-indices
\begin{align}
    \cbracket{\beta(i) ~|~ 1 \le i \le M} \subseteq \cbracket{ \beta \in \Z^{n-1}_{\ge 0}~|~ 0 \le \abs{\beta} \le d}~, \label{betaindex0}
\end{align}
define $M \times N$ matrices 
\begin{align}
  A(u , \tau) :& = \bracket{\tau^{d- \abs{\beta(i)}} B_{j, \beta(i)}\bracket{\frac{u}{\abs{u}}}}_{1 \le i \le M, 1 \le j \le N}~, \nonumber\\ \tilde{A}(u , \tau) :& = \bracket{\bracket{(\abs{u} + \tau)^{d-\abs{\beta(i)}} - \tau^{d-\abs{\beta(i)}}} B_{j , \beta(i)}\bracket{\frac{u}{\abs{u}}}}_{1 \le i \le M, 1 \le j \le N}~. \label{matrixdef}
\end{align}
Then, the coefficients in \eqref{coeffsigma} can be expressed in terms of matrix multiplication,
\begin{align}
  A(u,\tau) \begin{pmatrix}
    \nu_{1} - \mu_{1}\\
    \vdots \\
    \nu_{N} - \mu_{N}
    \end{pmatrix}
    =\begin{pmatrix}
     C[\sigma^{\beta(1)}](u, \tau)\\
    \vdots \\
    C[\sigma^{\beta(M)}](u, \tau) 
  \end{pmatrix}- \tilde{A} (u , \tau) \begin{pmatrix}
    \nu_{1} \\
    \vdots \\
    \nu_{N} 
    \end{pmatrix}~. \label{matrixexpression}
\end{align}

\subsection{Construction of exceptional sets}
Fix $r\ge 1$, $0< \epsilon_1 < 1$ and $\nu ,\mu \in \R^N$ such that $r \le \norm{\nu}, \norm{\mu} \le 2r$.
Fix $(u, \tau) \in B_2(\R^n) \times B_1(\R)$ where $u/|u|$ lies in the support of the function $W_n$ in \eqref{partition}.
If there is some $\beta \in \Z^{n-1}_{\ge 0}$ such that $1 \le \abs{\beta} \le d$ and $\abs{C[\sigma^{\beta}](u,\tau)} \ge r^{\epsilon_1 }$, then the Van der Corput estimate given by Lemma \ref{Vandecorput1} directly yields the bound
\begin{align*}
  \abs{K_{\sharp , n}^{\nu, \mu}(u, \tau)} \le C \bracket{\sum_{1 \le \abs{\beta} \le d}^{} \abs{C[\sigma^\beta](u,\tau)}}^{-1/d} \le Cr^{-\epsilon_1/ d} \chi_{B_2}(u) \chi_{B_1}(\tau)~,
\end{align*}
for some absolute constant $C$. The characteristic functions in this estimate are justified by the fact that $K_{\sharp, n}^{\nu, \mu}(u, \tau)$ is supported on $B_2(\R^n) \times B_1(\R)$ as a function of $(u , \tau)$. 
Hence, it suffices to construct a small exceptional set $G^\nu$, and for each $u$, a small exceptional set $F_u^\nu$ such that all $(u , \tau)$ satisfying $\abs{C[\sigma^{\beta}](u,\tau)} \le r^{\epsilon_1}$ for all $1 \le \abs{\beta} \le d$ are contained in the set
\begin{align*}
    \bracket{G^\nu \times B_1(\R) } \cup \bracket{\bigcup_{u \in B_2(\R^n)} \cbracket{u} \times F_u^\nu}~.
\end{align*}
If this is true, then since $K_{\sharp,n}^{\nu, \mu}(u, \tau)$ is bounded by some universal constant $C$ (since the compactly supported function $\Psi$ in \eqref{kernel} is assumed to have uniformly bounded $\CC^1$ norm), we obtain the full estimate \eqref{fullkernelbound} with $\delta = \epsilon_1 / d$:
\begin{align*}
  \abs{K_{\sharp , n}^{\nu, \mu}(u, \tau)} \le C \bracket{r^{-\delta} \chi_{B_2}(u) \chi_{B_1}(\tau) + \chi_{G^\nu}(u) \chi_{B_1}(\tau) + \chi_{B_2}(u) \chi_{F_u^\nu}(\tau)}~.
\end{align*}

The definition of the exceptional sets $G^\nu$ and $F_u^\nu$ will rely on an explicit $\mu$-independent representation of the coefficients $C[\sigma_s](u,\tau)$, $1 \le s \le n-1$.
From now on, we assume that the fixed point $(u , \tau)$ also satisfies 
\begin{align}
\abs{C[\sigma^\beta](u,\tau)} \le r^{\epsilon_1} \quad \text{for all } 2 \le \abs{\beta} \le d~, \label{assumptiononutau}
\end{align}
and proceed with finding the explicit representation of the coefficients $C[\sigma_s](u,\tau)$, $1 \le s \le n-1$.

Consider the set of multi-indices
\begin{align*}
  \cbracket{ \beta \in \Z^{n-1}_{\ge 0}~|~ 0 \le \abs{\beta} \le d}~.
\end{align*}
Homogenizing each multi-index $\beta = (\beta_1,...., \beta_{n-1})$ by adding a new component $\beta_n = d - \abs{\beta}$, we obtain that the cardinality of the above set is the same as the cardinality of the set $\cbracket{ \alpha \in \Z^{n}_{\ge 0}~|~ \abs{\alpha} = d}$. The latter correspond to the set of all monomials of degree $d$ in $n$ variables (a basis for $\mathcal{P}^{n, d}$), which has cardinality $\binom{d+n-1}{n-1}$. Therefore, the set of $\beta$ with $2 \le \abs{\beta} \le d$ has cardinality exactly $\binom{d+n-1}{n-1} - n = N = \dim(\mathcal{Q})$. Upon setting $M = N$ in \eqref{betaindex0} so that
\begin{align}
  \cbracket{\beta(i) ~|~ 1 \le i \le N} = \cbracket{ \beta \in \Z^{n-1}_{\ge 0}~|~ 2 \le \abs{\beta} \le d}~, \label{betaindex} 
\end{align}
consider the $N\times N$ matrix representation \eqref{matrixexpression} associated to it.
Since we have assumed \eqref{assumptiononutau}, we have $\abs{C[\sigma^{\beta(i)}] (u , \tau) } \le r^{\epsilon_1}$ for each $1 \le i \le N$, so \eqref{matrixexpression} simplifies to
\begin{align*}
  A(u,\tau) \begin{pmatrix}
    \nu_{1} - \mu_{1}\\
    \vdots \\
    \nu_{N} - \mu_{N}
    \end{pmatrix}
    =\begin{pmatrix}
     O(r^{\epsilon_1})\\
    \vdots \\
    O(r^{\epsilon_1}) 
  \end{pmatrix}- \tilde{A} (u , \tau) \begin{pmatrix}
    \nu_{1} \\
    \vdots \\
    \nu_{N} 
    \end{pmatrix} = : B_\nu(u , \tau)~.
\end{align*}
Now the right hand side only depends on $\nu$ and the fixed point $(u,\tau)$, and no longer involves $\mu$.
This gives rise to a linear system of $N$
equations in which $(\nu_i - \mu_i)$, $1\le i \le N$, are the unknown variables. We may solve this system to represent each $\nu_i - \mu_i$ by an expression which only depends on $\nu$. 
Using Cramer's rule, for each $1 \le i \le N$, we have
\begin{align}
  \det(A(u,\tau))(\nu_i - \mu_i) = \det(A^i_\nu(u,\tau))~, \label{cramerR2}
\end{align}
where $A^i_\nu(u,\tau)$ is the matrix obtained by replacing the $i$-th column of $A(u,\tau)$ by $B_\nu(u, \tau)$. This step serves the same purpose as the downward induction step in Section 9 of \cite{MR3945729}.

For each $1 \le s \le n-1$, expression \eqref{coeffsigma} simplifies to
\begin{align*}
  C[\sigma_s](u, \tau) = \tau^{d-1}\sum_{i = 1}^N \bracket{\nu_i - \mu_i} B_{i , e_s}\bracket{\frac{u}{\abs{u}}} + \bracket{(\abs{u} + \tau)^{d-1} - \tau^{d-1}} \sum_{i = 1}^N \nu_i B_{i , e_s}\bracket{\frac{u}{\abs{u}}}~.
\end{align*}
Notice that the only $\mu$-dependence in this expression comes from the terms $\nu_i - \mu_i$, which can now be replaced using \eqref{cramerR2}:
\begin{align}
    \det(A(u,\tau))C[\sigma_s](u, \tau) = \tau^{d-1} &\sum_{i = 1}^N \det(A^i_\nu(u,\tau)) B_{i , e_s}\bracket{\frac{u}{\abs{u}}}\nonumber \\ & + \bracket{(\abs{u} + \tau)^{d-1} - \tau^{d-1}} \sum_{i = 1}^N \nu_i \det(A(u,\tau)) B_{i , e_s}\bracket{\frac{u}{\abs{u}}}~.\label{CsigmaSexpression}
\end{align}
Multiplication by $\det(A(u,\tau))$ on the left hand side in expression \eqref{CsigmaSexpression} is necessary to ensure that the expression on the right hand side remains a polynomial in $\tau$, which is crucial for the later application of Van der Corput sublevel set estimates.
The determinant $\det(A^i_\nu(u,\tau))$ also depends linearly on the vector $\nu \in \R^N$. We verify this in the following lemma.

\begin{lem} \label{minorexp} Let $D\in \mathcal{P}^{n , Nd}$, $D_{j,i} \in \mathcal{P}^{n , (N-1)d}$, $1 \le i , j \le N$, denote the determinants
\begin{align*}
  D(X) : =  \det \bracket{ \bracket{B_{j, \beta(i)}\bracket{X}}_{1 \le i,j\le N}}~,~ D_{j , i}(X) : =\det\bracket{\bracket{B_{s,\beta(k)}\bracket{X}}_{\substack{1\le k ,s \le N \\ k\neq j, s \neq i}}}~,
\end{align*}
where the polynomials $B_{j, \beta}$ are defined as in \eqref{Bdef}, and the multi-indices $\beta(1) , ..., \beta(N)$ are as in \eqref{betaindex}. Then, one has
\begin{align}
  \det \bracket{A(u , \tau)} = \tau^{Nd - \sum_{i = 1}^{N}\abs{\beta(i)}} \abs{u}^{-Nd }D(u)~.\label{detAexpression}
\end{align}
For a multi-index $\beta \in \Z^{n-1}_{\ge 0}$, let $K^\beta: \mathcal{P}^{n , d} \to \mathcal{P}^{n , d}$ denote the linear map
\begin{align*}
  K^\beta p (X) :=  [\sigma^\beta]p(\omega(\sigma, X))~ ,
\end{align*}
where $\omega$ is the change of variables defined in \eqref{changeofvar2}.
For each $\nu \in \R^N$, set $p_\nu: = \sum_{k = 1}^N \nu_k p_k\in \mathcal{Q}$. Then, for each $1 \le i \le N $, 
\begin{align}
  \det & (A^i_\nu(u,\tau)) = O(r^{\epsilon_1}) \nonumber\\
  & - \tau^{(N - 1)d - \sum_{k = 1}^{N}\abs{\beta(k)}} \abs{u}^{-Nd} \sum_{ j =1}^N \tau^{\abs{\beta(j)} } (-1)^{i + j}D_{j , i}(u)\bracket{(\abs{u} + \tau)^{d - \abs{\beta(j)}} - \tau^{d- \abs{\beta(j)}}} K^{\beta(j)}p_\nu(u)~.\label{minorR2}
\end{align}
\end{lem}
\begin{proof}
Equation \eqref{detAexpression} follows from the multilinear property of the determinant function and homogeneity of the polynomials $B_{j, \beta}$.
We now show \eqref{minorR2}.
Minor expansion along the $i$-th column of $A^i_\nu(u,\tau)$ yields
\begin{align*}
  \det(A^i_\nu(u,\tau)) &= \sum_{j = 1}^{N} (-1)^{i+j} \bracket{O(r^{\epsilon_1}) - \bracket{\tilde{A}(u,\tau) \nu }_{j}}A_{ji}~,
\end{align*}
where $\tilde{A}(u,\tau)$ is the matrix defined in \eqref{matrixdef} and $A_{j,i}$ is the $(j,i)$-th minor of $A(u , \tau)$. Since $B_{i,\beta}$ is homogeneous of degree $d$, we have
\begin{align}
  A_{j,i} &=  \tau^{(N-1)d - \sum_{k \neq j}^{} \abs{\beta(k)}} \abs{u}^{-(N-1)d} \det\bracket{\bracket{B_{s,\beta(k)}\bracket{u}}_{\substack{1\le k,s \le N \\ k\neq j, s \neq i}}} \nonumber\\
  &= \tau^{(N-1)d - \sum_{k=1}^{N} \abs{\beta(k)} + \abs{\beta(j)}} \abs{u}^{-(N-1)d} D_{j,i}(u)~.\label{substitute1}
\end{align}
We have $A_{j,i} = O(1)$, since $(u,\tau) \in B_2(\R^n) \times B_1(\R)$, and all entries of $\bracket{B_{j, \beta(i)}\bracket{u/\abs{u}}}_{1 \le i,j\le N}$ are $O(1)$. It follows that moving the $O(r^{\epsilon_1})$ term outside of the sum in the minor expansion is justified, and we obtain
\begin{align}
  \det(A^i_\nu(u,\tau)) &= O(r^{\epsilon_1}) - \sum_{j = 1}^{N} (-1)^{i+j} \bracket{\tilde{A}(u,\tau) \nu }_{j}A_{ji}~,\label{minorexp1}
\end{align}
Finally, recall from \eqref{matrixdef} that
\begin{align}
  \bracket{\tilde{A}(u,\tau) \nu }_{j} & = \abs{u}^{-d}\bracket{(\abs{u} + \tau)^{d-\abs{\beta(j)}} - \tau^{d-\abs{\beta(j)}}} K^{\beta(j)}p_{\nu}(u)~, \label{substitute2}
\end{align}
The claim of \eqref{minorR2} then follows from substituting the expressions \eqref{substitute1} and \eqref{substitute2} into \eqref{minorexp1}. 
\end{proof}

For fixed $\nu \in \R^N$, let $p_\nu : = \sum_{i = 1}^N \nu_i p_i \in \mathcal{Q}$. By Lemma \ref{minorexp}, substituting the explicit expression \eqref{minorR2} into \eqref{CsigmaSexpression}, we obtain
\begin{align}
  \det(A(u,\tau))& C[\sigma_s](u , \tau)\nonumber \\
  =&~  O(r^{\epsilon_1}) - \tau^{Nd - \sum_{k = 1}^{N} \abs{\beta(k)}} \abs{u}^{-Nd-d}\sum_{ i =1}^N \sum_{ j =1}^N \tau^{\abs{\beta(j)} - 1 }  B_{i , e_s}(u)(-1)^{i + j}D_{j , i}(u) \nonumber\\
  & ~~~\qquad\qquad\qquad\qquad\qquad\qquad\qquad\qquad\qquad\qquad \cdot \bracket{(\abs{u} + \tau)^{d - \abs{\beta(j)}} - \tau^{d- \abs{\beta(j)}}} K^{\beta(j)}p_\nu(u)\nonumber\\
   &\qquad\qquad + \tau^{Nd - \sum_{k = 1}^{N} \abs{\beta(k)}} \abs{u}^{-Nd-d}\bracket{(\abs{u} + \tau)^{d-1} - \tau^{d-1}} D(u)K^{e_s}p_\nu(u)~.\label{CsigmaSexpression2}
\end{align}
The operation of moving the $O(r^{\epsilon_1})$ term outside of the sum is again justified by the observation that $B_{i , \beta}(u / \abs{u})= O(1)$ for all $u \neq 0$. The $O(r^{\epsilon_1})$ term in \eqref{CsigmaSexpression2} is then bounded above by $Cr^{\epsilon_1}$, uniformly over $(u , \tau) \in B_2(\R^n) \times B_1(\R)$, and $\nu , \mu \in \R^N$, with an absolute constant $C$ which only depends on $n , d$ and the fixed basis $\{p_1,...,p_N \}$ of $\mathcal{Q}$.

We use the notation $R_{\nu , u}^{s}(\tau)$ to denote the polynomial in $\tau$ excluding the $O(r^{\epsilon_1})$ term in \eqref{CsigmaSexpression2}, so that
\begin{align*}
\det(A(u,\tau)) C[\sigma_s](u , \tau) = O(r^{\epsilon_1}) + R_{\nu , u}^{s}(\tau)~.
\end{align*}
More explicitly, we have
\begin{align}
  \tau^{- Nd + \sum_{i = 1}^{N} \abs{\beta(i)}}& |u|^{Nd + d} R_{\nu , u}^{s}(\tau) \nonumber\\ = & - \sum_{ i =1}^N \sum_{ j =1}^N \tau^{\abs{\beta(j)} - 1 }  B_{i , e_s}(u)(-1)^{i + j}D_{j , i}(u)\bracket{(\abs{u} + \tau)^{d - \abs{\beta(j)}} - \tau^{d- \abs{\beta(j)}}} K^{\beta(j)}p_\nu(u)\nonumber\\ 
  & +   \bracket{(\abs{u} + \tau)^{d-1} - \tau^{d-1}} D(u)K^{e_s}p_\nu(u)~. \label{linearsigmaterm}
\end{align}
We selected the set $\{\beta(j)\}$ to satisfy $2 \le \abs{\beta(j)} \le d$ for all $1 \le j \le N$. Thus, the terms which are constant in $\tau$ on the right hand side of \eqref{linearsigmaterm} come solely from the last term, while the terms which are linear in $\tau$ come from both the last term and the terms with $\abs{\beta(j)} = 2$ in the double sum.
This allows us to compute the coefficients of the two lowest order monomials in $R_{\nu , u}^{s}(\tau)$ as a polynomial in $\tau$:
\begin{align*}
  R_{\nu , u}^{s}(\tau)  =& ~  \tau^{Nd - \sum_{i = 1}^{N} \abs{\beta(i)}} \abs{u}^{-Nd - 1} {D(u)K^{e_s}p_\nu(u)}\\ &+ \tau^{Nd - \sum_{i = 1}^{N} \abs{\beta(i)}+1}  \abs{u}^{-Nd - 2} \bracket{ (d-1)D(u)K^{e_s}p_\nu(u) -  \sum_{j : \abs{\beta(j)} = 2} D_j^{e_s}(u) K^{\beta(j)}p_\nu(u) } \\
  & + O\bracket{\tau^{Nd - \sum_{i = 1}^{N} \abs{\beta(i)}+2} }~,
\end{align*}
where $D_j^{e_s} \in \mathcal{P}^{n , Nd}$ is the polynomial
\begin{align*}
  D_j^{e_s}(X) : = \sum_{i = 1}^{N} (-1)^{i+j} B_{i, e_s}(X) D_{j,i}(X)~.
\end{align*}

For each $1 \le s \le n-1$, define linear maps $L^{(1)}_s , L^{(2)}_s : \mathcal{P}^{n , d} \to \mathcal{P}^{n , (N+1)d}$, acting on polynomials $q \in \mathcal{P}^{n , d}$, by
\begin{align}
  L^{(1)}_s q  &:=  D K^{e_s} q ~~,~~L^{(2)}_s q := \sum_{j : \abs{\beta(j)} = 2} D_j^{e_s} K^{\beta(j)} q~.\label{L1andL2}
\end{align}
Then, $R_{\nu , u}^{s}(\tau)$ can be expressed in terms of the image of $p_\nu$ under the maps $L^{(1)}_s, L^{(2)}_s$ as follows:
\begin{align}
  \tau^{-Nd + \sum_{i = 1}^{N} \abs{\beta(i)}}&R_{\nu , u}^{s}(\tau) \nonumber \\  = & ~   \abs{u}^{-Nd - 1} L^{(1)}_s p_\nu(u) + \tau  \abs{u}^{-Nd - 2} \bracket{(d-1) L^{(1)}_s p_\nu(u) -  L^{(2)}_s p_\nu(u) } + O\bracket{\tau^2}~. \label{expressionofRtau}
\end{align}
Here, the notation $O(\tau^2)$ refers to a polynomial in $\tau$ such that every nonzero monomial in it is of order at least $2$.

Fix $\epsilon_2$ with $0 < \epsilon_1 < \epsilon_2 < 1$. 
Define the exceptional set
\begin{align}\label{exceptionalsets}
  \begin{split}
    G^\nu := \bigcap_{s=1}^{n-1} &\bigg(\left\{ u \in B_2(\R^n)~~|~~ \abs{L^{(1)}_sp_\nu(u)} \le r^{\epsilon_2}\right\} \\ 
  &\hspace{2cm}\cap \left\{ u \in B_2(\R^n)~~|~~ \abs{(d-1) L^{(1)}_s p_\nu(u) -  L^{(2)}_s p_\nu(u) } \le r^{\epsilon_2}\right\}\bigg)~. 
  \end{split}
\end{align}
For each $u \in B_2(\R^n)$, define 
\begin{align*} 
  F_u^\nu = \left\{\begin{array}{cl}
  \bigcap_{s=1}^{n-1}\cbracket{ \tau \in B_1(\R)~~|~~ \abs{R_{\nu, u }^{s}(\tau)} \le C_0r^{\epsilon_1}} & \mbox{ if } u \in B_2(\R^n) \backslash G^\nu \\ [5pt]
    \emptyset & \mbox{ if } u \in G^\nu
  \end{array}\right. \, ~,
\end{align*} 
where $C_0$ is an absolute constant, independent of $\nu$, to be determined later.

\begin{lem}
For an appropriate choice of $C_0$ which depends only on $n , d$, and the basis $\{p_1,..., p_N\}$ of $\mathcal{Q}$, we have
\begin{align*}
  \{ (u , \tau) \in (B_2(\R^n) \backslash \{0\}) \times B_1(\R)~~|~~ \abs{C[\sigma^\beta](u,\tau)} &\le r^{\epsilon_1}~,~\forall1 \le \abs{\beta} \le d \} \\ & \subseteq \bracket{G^\nu \times B_1(\R)} \cup \bracket{\bigcup_{u \in B_2(\R^n)} \cbracket{u} \times F_u^\nu}~.
\end{align*}

\end{lem}
\begin{proof}
It suffices to show that if $(u,\tau)$ is not in the set
\begin{align*}
  \bracket{G^\nu \times B_1(\R)} \cup \bracket{\bigcup_{u \in B_2(\R^n)} \cbracket{u} \times F_u^\nu}~,
\end{align*}
then there is some multi-index $\beta$ with $1 \le \abs{\beta} \le d$ such that $\abs{C[\sigma^\beta](u,\tau)} > r^{\epsilon_1}$.

Suppose $u \notin G^\nu$ and $\tau \notin F_u^\nu$ and $\abs{C[\sigma^\beta](u,\tau)} \le r^{\epsilon_1}$ for all $2 \le \abs{\beta} \le d$. We will show that there is $1 \le s \le n-1$ such that $\abs{C[\sigma_s](u , \tau)} > r^{\epsilon_1}$. Recall from \eqref{CsigmaSexpression2} that for each $1 \le s \le n-1$, we have the representation
\begin{align*}
    \det(A(u,\tau))C[\sigma_s](u , \tau) = O(r^{\epsilon_1}) + R_{\nu , u}^{s}(\tau)~,~
\end{align*}
for $(u,\tau) \in B_2(\R^n) \times B_1(\R)$ with $\abs{C[\sigma^\beta](u,\tau)} \le r^{\epsilon_1}$ for all $2 \le \abs{\beta} \le d$.
By construction of the exceptional sets, there is some $1 \le s_{u,\tau}\le n-1$, potentially depending on both $u,\tau$, such that $\abs{R_{\nu, u }^{s_{u , \tau}}(\tau)} > C_0r^{\epsilon_1}$. Hence, we have
\begin{align}
  \abs{C[\sigma_{s_{u,\tau}}](u , \tau)} \ge C \abs{\det(A(u , \tau))C[\sigma_{s_{u,\tau}}](u , \tau) }\label{C0}
  &\ge C\abs{R_{\nu , u}^{{s_{u,\tau}}}(\tau)} - O(r^{\epsilon_1}) \ge C_0 r^{\epsilon_1} - O(r^{\epsilon_1})~.
\end{align}
The first inequality in \eqref{C0} is justified by the observation that if $u \neq 0$ and $\tau \in B_1(\R)$, then
\begin{align*}
  \det \bracket{A(u , \tau)} = \tau^{Nd - \sum_{i = 1}^{N}\abs{\beta(i)}} \det \bracket{ \bracket{ B_{i , \beta(j)}\bracket{\frac{u}{\abs{u}}} }_{1 \le i, j \le N} } = O(1)~,
\end{align*}
with an upper bound that depends only on combinatorial constants and the basis $\{p_1,..., p_N\}$ of $\mathcal{Q}$. If $C_0$ is sufficiently large, then
\begin{align}
    \abs{C[\sigma_{s_{u,\tau}}](u , \tau)}  & \ge C_0 r^{\epsilon_1} - O(r^{\epsilon_1})\nonumber \ge 2 r^{\epsilon_1}~.\nonumber
\end{align}
Note that the choice of $C_0$ here depends only on the $O(r^{\epsilon_1})$ term in the expansion for each $1 \le s \le n-1$ from \eqref{CsigmaSexpression2}, which is indeed uniform in $u , \tau, \nu , \mu$. Therefore, we have successfully produced a coefficient $C[\sigma_{s_{u,\tau}}](u , \tau)$ that is large.

\end{proof}

\subsection{Estimate of measure of exceptional sets}

It remains to verify that $G^\nu$, $F_u^\nu$ are sets of small measure. We first look at $F_u^\nu$. 
\begin{lem}\label{mainlemFnu}
  For each $r\ge 1$, $ \nu \in \R^N$ such that $r \le ||\nu || \le 2r$, and $u \in B_2(\R^n)\backslash G^\nu$, we have
  \begin{align*}
      \abs{F_u^\nu} \le C r^{(\epsilon_1 - \epsilon_2)/\delta'}~.
  \end{align*}
  Here, $C > 0$ is an absolute constant depending only on $n , d$, and the basis  $\{p_1,..., p _N \}$ of $\mathcal{Q}$, and
  \begin{align*}
    \delta' := \max\cbracket{(N+1)d - \sum_{i = 1}^{N} \abs{\beta(i)} -2, 1}~.
  \end{align*}
  
\end{lem}
\begin{proof}
  By construction, $F_u^\nu$ is the empty set when $u \in G^\nu$, so it suffices to only consider the case $u \notin G^\nu$.
  For each fixed $u \in B_2(\R^n)\backslash G^\nu$, there is $1 \le s_u \le n-1$ such that either
\begin{align*}
  \abs{L^{(1)}_{s_u}p_\nu(u)} > r^{\epsilon_2}~,~ \mbox{ or }\abs{(d-1) L^{(1)}_{s_u} p_\nu(u) -  L^{(2)}_{s_u} p_\nu(u) } > r^{\epsilon_2}~.
\end{align*} 
Recall the notation that $\inorm{q(\tau)}_{\tau}$ denotes the sum of the absolute values of the coefficients of the polynomial $q(\tau)$ from \eqref{isotropicnorm}.
Referring back to the expression of $R_{\nu, u}^{s}(\tau)$ in \eqref{expressionofRtau}, since $\abs{u} \le 2$, we have $\abs{u}^{-1} \ge 1/2$. 
Hence, we obtain the lower bound
\begin{align*}
  \inorm{ R_{\nu , u }^{{s_u}}(\tau)}_{\tau} \ge \abs{u}^{-Nd - 1} \abs{L^{(1)}_{s_u}p_\nu(u)} + \abs{u}^{-Nd - 2}\abs{(d-1) L^{(1)}_{s_u} p_\nu(u) -  L^{(2)}_{s_u} p_\nu(u) } > Cr^{\epsilon_2}~,
\end{align*}
for a constant $C > 0$ depending only on $d$ and $N$, so that $R_{\nu , u }^{s_u}(\tau)$ is a nonzero polynomial in $\tau$.
The Van der Corput sublevel set estimate of Lemma \ref{Vandecorput2} then yields
\begin{align*}
  \abs{F_u^\nu} & \le \abs{\cbracket{ \tau \in B_1(\R)~~|~~ \abs{R_{\nu, u }^{{s_u}}(\tau)} \le C_0r^{\epsilon_1}}} \le \begin{cases}
    C r^{\epsilon_1/\delta_u}\inorm{R_{\nu , u}^{{s_u}}(\tau)}_{\tau}^{- 1/\delta_u} & \mbox{ if } \delta_u \ge 1\\
    C r^{\epsilon_1}\inorm{R_{\nu , u}^{{s_u}}(\tau)}_{\tau}^{- 1} & \mbox{ if } \delta_u =0
  \end{cases}~, 
\end{align*}
where $\delta_u$ denotes the degree of $R_{\nu , u}^{{s_u}}(\tau)$ as a polynomial in $\tau$. The range of $\delta_u$ is given by
\begin{align*}
  0 \le Nd - \sum_{i = 1}^{N} \abs{\beta(i)} \le \delta_u & \le \max\cbracket{(N+1)d - \sum_{i = 1}^{N} \abs{\beta(i)} -2, 1}~.
\end{align*}
We then define $\delta': = \max\cbracket{(N+1)d - \sum_{i = 1}^{N} \abs{\beta(i)} -2, 1}$
to remove dependence on $u$, and obtain
\begin{align*}
  \abs{F_u^\nu} \le C r^{(\epsilon_1 - \epsilon_2)/\delta'}~.
\end{align*}

\end{proof}

Finally, we bound the measure of $G^\nu$. This is the central new idea of the paper, relative to the previous works by \cite{MR3945729} and \cite{MR4790541}.
\begin{lem}\label{mainlemGnu}
Fix $\epsilon_3$ such that $0 < \epsilon_1 < \epsilon_2 < \epsilon_3<1$, where $\epsilon_2$ is the previously chosen parameter in the definition of the exceptional sets \eqref{exceptionalsets}.
Then, for each $r\ge 1$ and $ \nu \in \R^N$ such that $r \le ||\nu || \le 2r$, we have
  \begin{align*}
    \abs{G^\nu} \le C r^{(\epsilon_2 - \epsilon_3)/(N+1)d}~.
  \end{align*}
  Here, $C > 0$ is an absolute constant depending only on $n , d$, and the basis  $\{p_1,..., p _N \}$ of $\mathcal{Q}$.
\end{lem}

To prove Lemma \ref{mainlemGnu}, we require the following lemmas about the maps $L^{(1)}_s $ and $ L^{(2)}_s$ in \eqref{L1andL2} to effectively apply the Van der Corput estimate in Lemma \ref{Vandecorput2} to $G^\nu$. The proof of these lemmas is deferred to Section \ref{lemma}. 

\begin{lem}\label{L1} 
Let $\mathcal{Q}\subseteq \mathcal{P}^{n , d}$ be a $N: = \binom{d+n-1}{n-1}-n$ dimensional subspace with basis $\{p_1,..., p _N \}$. Assume that $\mathcal{Q}$ satisfies the nondegeneracy condition \eqref{cond1} for $c \in (\R \backslash \{0\})^n$.
Consider the linear operator $L^{(1)}: \mathcal{P}^{n , d} \to \bracket{\mathcal{P}^{n , (N+1)d}}^{n-1}$ defined by
\begin{align*}
  L^{(1)}q:= (L^{(1)}_1 q , ... ,L^{(1)}_{n-1} q)~,
\end{align*}
where $\bracket{\mathcal{P}^{n , (N+1)d}}^{n-1}$ is taken to be the vector space formed by a direct sum of $n-1$ copies of $\mathcal{P}^{n , (N+1)d}$, and $L_s^{(1)}$, $1\le s \le n-1$ is defined as in \eqref{L1andL2}.
Then, \begin{align*} 
  \ker \bracket{L^{(1)}}   = \left\{\begin{array}{cl}
    \mathrm{span}_\R \cbracket{ \abs{X}^d }  & \mbox{ if $d$ is even }  \\ [5pt]
    \{0\} & \mbox{ if $d$ is odd }
  \end{array}\right. \, .
\end{align*} 

\end{lem}

\begin{lem}\label{L2} Let $\mathcal{Q}\subseteq \mathcal{P}^{n , d}$ be a codimension $n$ subspace equipped with basis $\{p_1,..., p _N \}$. Assume that $\mathcal{Q}$ satisfies the hypothesis of Proposition \ref{mainprop}.
Consider the operators $L^{(2)}_s : \mathcal{P}^{n , d} \to \mathcal{P}^{n , (N+1)d}$, $1\le s \le n-1$, defined in \eqref{L1andL2}.
Then, there is $1 \le s_* \le n-1$ depending only on $n ,d$ and the basis $\{p_1,..., p_N\}$ of $\mathcal{Q}$, and some constant $C_{s_*}> 0$ depending only on $s_*$, such that for all $p \in \ker \bracket{L^{(1)}}$, one has
\begin{align*}
  \inorm{L^{(2)}_{s_*}p(X)}_{X}  = C_{s_*} \norm{p}_{\ell^1(\mathcal{Q})}~.
\end{align*}
Recall that the norms $\inorm{\cdot}_{X}$, $\norm{\cdot}_{\ell^1(\mathcal{Q})}$ are defined respectively in \eqref{isotropicnorm} and \eqref{l1norm} in Section \ref{preliminaries}.

\end{lem}
Note that Lemma \ref{L1} claims that the kernel of $L^{(1)}$ is either trivial or one-dimensional, depending on whether $d$ is odd or even. Hence, the result of Lemma \ref{L2} is in fact only relevant when $d$ is even.

The intuition for applying Lemmas \ref{L1} and \ref{L2} to prove Lemma \ref{mainlemGnu} is as follows. First, for $\rho>0$ and an arbitrary homogeneous polynomial $q \in \mathcal{P}^{n , (N+1)d} \backslash\{0\}$, the Van der Corput estimate in Lemma \ref{Vandecorput2} gives
\begin{align*}
  \abs{\cbracket{ u \in B_2(\R^n)~~|~~ \abs{q(u)} \le \rho}} & \le C\rho^{1/(N+1)d}[\![q(X)]\!]_X^{-1/(N+1)d}~.
\end{align*}
As mentioned earlier in Section \ref{preliminaries}, $[\![q(X)]\!]_X$ on the right hand side is a well-defined norm on $\mathcal{P}^{n , (N+1)d}$ while $\norm{\cdot}_{\ell^1(\mathcal{Q})}$ is a well-defined norm on $\mathcal{P}^{n , d}$. In view of the definition of $G^\nu$, we want to take $\rho = r^{\epsilon_2}$ and $q = Lp_\nu$, while either $L = L^{(1)}_s$ or $L = (d-1)L^{(1)}_s - L^{(2)}_s$, for some $1\le s \le n-1$. If the map
\begin{align*}
  L: \mathcal{P}^{n , d} \to \mathcal{P}^{n , (N+1)d}~,~ 1 \le s \le n-1
\end{align*}
is injective on $\mathcal{Q}$, then the inverse map
\begin{align*}
  (L|_{\mathcal{Q}})^{-1} : \bracket{L|_{\mathcal{Q}}(\mathcal{Q}), \inorm{\cdot}_{X}} \to \bracket{\mathcal{Q},\norm{\cdot}_{\ell^1(\mathcal{Q})}}~
\end{align*}
is well-defined and necessarily bounded with the indicated norms since all linear maps between finite-dimensional normed vector spaces are continuous and all norms on a finite-dimensional vector space are equivalent. If injectivity can be achieved, then we can make use of the assumption that for the given $r\ge 1$,
\begin{align*}
  r \le \norm{\nu}= \sum_{i=1}^{N} \abs{\nu_i}= \norm{p_\nu}_{\ell^1(\mathcal{Q})}~~
\end{align*}
and the Van der Corput sublevel set estimate from Lemma \ref{Vandecorput2}
to obtain the desired bound
\begin{align*}
    \abs{G^\nu} \le C r^{\epsilon_2/(N+1)d} \inorm{L p_\nu(X)}_X^{-1/(N+1)d}~\le  C r^{(\epsilon_2-1)/(N+1)d}~.
\end{align*}

We are unable to guarantee existence of a particular $1\le s \le n-1$ such that $L^{(1)}_s$ or $(d-1)L^{(1)}_s - L^{(2)}_s$ is injective on $\mathcal{Q}$. However, Lemmas \ref{L1} and \ref{L2} together imply the weaker claim that there is at least one $1 \le s \le n-1$ such that $L^{(2)}_s$ is injective on the kernel of
\begin{align*}
  L^{(1)} = (L^{(1)}_1, \dots, L^{(1)}_{n-1}) : \mathcal{P}^{n , d} \to \bracket{\mathcal{P}^{n , (N+1)d}}^{n-1}~,
\end{align*}
where the codomain is the vector space formed by a direct sum of $n-1$ copies of $\mathcal{P}^{n , (N+1)d}$. In particular,
\begin{align*}
  \ker{L^{(1)}} = \bigcap_{s= 1}^{n-1} \ker{L^{(1)}_s}~,
\end{align*}
and this kernel is trivial in the case when $d$ is odd and $\mathrm{span}_\R\{\abs{X}^d\}$ in the case when $d$ is even. Given the polynomial $p_\nu \in \mathcal{Q}$, if it's ``far'' from the kernel of $L^{(1)}$, we invert $L^{(1)}$ modulo its kernel
and apply Lemma \ref{Vandecorput2} (for a suitable $s$) to the set
\begin{align*}
  \cbracket{ u \in B_2(\R^n)~~|~~ \abs{L^{(1)}_sp_\nu(u)} \le r^{\epsilon_2}}
\end{align*}
to obtain the desired upper bound for $\abs{G^\nu}$. If $p_\nu$ is ``close'' to the kernel of $L^{(1)}$, we instead invert $L^{(2)}_s$ (for a suitable $s$) on the kernel of $L^{(1)}$ and apply Lemma \ref{Vandecorput2} to the set
\begin{align*}
  \cbracket{ u \in B_2(\R^n)~~|~~ \abs{(d-1) L^{(1)}_s p_\nu(u) -  L^{(2)}_s p_\nu(u) } \le r^{\epsilon_2}}
\end{align*}
to obtain the desired upper bound for $\abs{G^\nu}$.

We now prove Lemma \ref{mainlemGnu}, which provides the desired upper bound for $\abs{G^\nu}$. This finishes the proof of Proposition \ref{mainprop}, aside from the verification of Lemma \ref{L1} and Lemma \ref{L2}.
\begin{proof}

To make precise the notion of being ``close'' and ``far'' from the kernel
\begin{align*}
    K : = \ker(L^{(1)})~,
\end{align*}
we fix an arbitrary inner product on $\mathcal{P}^{n , d}$, and consider the orthogonal decomposition of $p_\nu \in \mathcal{Q}$ with respect to $K$ and $K^\perp$ in $\mathcal{P}^{n,d}$. We adopt the notation
\begin{align*}
  p_\nu (X) = f_\nu(X) + g_\nu(X)~, \mbox{ where } f_\nu(X) \in K ~~,~~ g_\nu(X) \in K^\perp~.
\end{align*}
Fix $0 < \epsilon<1$ to be a small absolute constant to be chosen later.
We consider separately the cases $ \norm{g_\nu}_{\ell^1(\mathcal{Q})} \ge \epsilon r^{\epsilon_3}$ (``far'' from $K$) and $ \norm{g_\nu}_{\ell^1(\mathcal{Q})} <\epsilon r^{\epsilon_3}$ (``close'' to $K$). 

\subsubsection{Case A: $\norm{g_\nu}_{\ell^1(\mathcal{Q})} \ge \epsilon r^{\epsilon_3}$} For each $1 \le s \le n-1$, if $L^{(1)}_s p_\nu$ is not the zero polynomial,
then by Lemma \ref{Vandecorput2}, we have
\begin{align}
  \abs{G^\nu} \le \abs{\bigcap_{s= 1}^{n-1} \cbracket{ u \in B_2(\R^n)~~|~~ \abs{L^{(1)}_{s}p_\nu(u)} \le r^{\epsilon_2}}} & \le \abs{\cbracket{ u \in B_2(\R^n)~~|~~ \abs{L^{(1)}_{s}p_\nu(u)}  \le r^{\epsilon_2}}} \nonumber\\ &\le C r^{\epsilon_2 / (N+1)d} [\![L^{(1)}_{s}p_\nu(X)]\!]_X^{-1 / (N+1)d}~,\label{halfstep}
\end{align}
where $[\![\cdot]\!]_X$ denotes the isotropic norm of a polynomial in $\mathcal{P}^{n,(N+1)d}$ defined in \eqref{isotropicnorm}.
Recall our decomposition $p_\nu = f_\nu + g_\nu$, where $f_\nu \in K$ and $g_\nu \in K^\perp$. Since $L^{(1)}$ vanishes on $K$, we have $L^{(1)}p_\nu = L^{(1)}g_\nu$.
This motivates us to find $s = s_\nu$ such that
\begin{align}
  [\![L^{(1)}_{s_\nu}g_\nu(X)]\!]_X = [\![L^{(1)}_{s_\nu}p_\nu(X)]\!]_X \ge C \norm{g_\nu}_{\ell^1(\mathcal{Q})}~, \label{lowerboundL1}
\end{align}
for an absolute constant $C$ that only depends on $n , d$, and the basis $\{p_1,..., p_N\}$. If a bound of the form \eqref{lowerboundL1} can be achieved, then the assumption $\norm{g_\nu}_{\ell^1(\mathcal{Q})} \ge \epsilon r^{\epsilon_3}$ directly implies that $L^{(1)}_{s_\nu}p_\nu$ is not the zero polynomial, and we can combine with \eqref{halfstep} to obtain
\begin{align*}
    \abs{G^\nu} & \le C r^{\epsilon_2 / (N+1)d} [\![L^{(1)}_{s_\nu}g_\nu(X)]\!]_X^{-1 / (N+1)d} \le C \epsilon^{-1 / (N+1)d} r^{(\epsilon_2 - \epsilon_3) / (N+1)d}~.
\end{align*}

With the goal of finding such an $s_\nu$ for which \eqref{lowerboundL1} holds, we analyze the restricted map \begin{align*}
  L^{(1)}|_{K^\perp} : K^\perp  \longrightarrow  L^{(1)}(K^\perp) \subseteq (\mathcal{P}^{n,(N+1)d})^{n-1}~.
\end{align*}
This is a well-defined bijective linear map between finite-dimensional vector spaces, by definition of the kernel $K$ of $L^{(1)}$.
We equip the external direct sum $(\mathcal{P}^{n, (N+1)d})^{n-1}$ with the $\ell^\infty$ norm composed with the isotropic norm on $\mathcal{P}^{n , (N+1)d}$, so that
\begin{align}
  \norm{ \bracket{q_1, ...  , q_{n-1}}}_{\ell^\infty, X} = \max_{1 \le s \le n-1} \inorm{q_s(X)}_X~, \mbox{ if } (q_1,...,q_{n-1}) \in \bracket{\mathcal{P}^{n, (N+1)d}}^{n-1}~.\label{linfnorm}
\end{align}
Hence, the inverse map
\begin{align*}
  \bracket{L^{(1)}|_{K^\perp}}^{-1} : \bracket{ L^{(1)}(K^\perp) ,  \norm{\cdot}_{\ell^\infty, X}}\longrightarrow \bracket{K^\perp, \norm{\cdot}_{\ell^1(\mathcal{Q})}}
\end{align*}
is well-defined and bounded with respect to the indicated norms. In particular, the operator norm
\begin{align*}
  C' := \norm{\bracket{L^{(1)}|_{K^\perp}}^{-1}}_{(L^{(1)}(K^\perp), \norm{\cdot}_{\ell^{\infty}, X}) \longrightarrow (K^\perp, \norm{\cdot}_{\ell^1(\mathcal{Q})})}
\end{align*}
is dependent solely on the construction of the linear map $L^{(1)}$ in \eqref{L1andL2}, the choice of the subspace $K^\perp$, and the norms $\norm{\cdot}_{\ell^\infty, X}$, $\norm{\cdot}_{\ell^1(\mathcal{Q})}$. All of these objects are dependent only on $n , d$, and the basis $\{p_1,..., p_N\}$ of $\mathcal{Q}$, and independent of $\nu$. Also, Lemma \ref{L1} ensures that $K^\perp \neq \{0\}$, so $C' > 0$.

Therefore, by definition of the norm \eqref{linfnorm}, there is $1 \le s_\nu \le n-1$, depending on $\nu$, such that
\begin{align*}
    \inorm{L^{(1)}_{s_\nu}p_\nu(X)}_X = \max_{1 \le s \le n-1} [\![L^{(1)}_sg_\nu(X)]\!]_X &= \norm{ L^{(1)}g_\nu}_{\ell^\infty,X}~.
\end{align*}
It follows from boundedness of the inverse map that
\begin{align*}
  \inorm{L^{(1)}_{s_\nu}p_\nu(X)}_X &= \norm{ L^{(1)}g_\nu}_{\ell^\infty,X} \ge (C')^{-1} \norm{ \bracket{L^{(1)}|_{K^\perp}}^{-1} \circ L^{(1)}g_\nu}_{\ell^1(\mathcal{Q})}= (C')^{-1} \norm{g_{\nu}}_{\ell^1(\mathcal{Q})} ~,
\end{align*}
which finishes the proof of Case A.

\subsubsection{Case B: $\norm{g_\nu}_{\ell^1(\mathcal{Q})} < \epsilon r^{\epsilon_3}$} 
Using the estimate of Lemma \ref{Vandecorput2}, for any $1 \le s \le n-1$, if $(d-1) L^{(1)}_{s} p_\nu -  L^{(2)}_s p_\nu $ is not the zero polynomial,
then Lemma \ref{Vandecorput2} gives
\begin{align*}
  \abs{G^\nu} &\le \abs{\bigcap_{s=1}^{n-1}\cbracket{ u \in B_2(\R^n)~~|~~ \abs{(d-1) L^{(1)}_{s} p_\nu(u) -  L^{(2)}_s p_\nu(u)} \le r^{\epsilon_2}}} \\ 
  & \le \abs{\cbracket{ u \in B_2(\R^n)~~|~~ \abs{(d-1) L^{(1)}_{s} p_\nu(u) -  L^{(2)}_{s} p_\nu(u)} \le r^{\epsilon_2}}} \\ &\le C r^{\epsilon_2 / (N+1)d}\inorm{(d-1) L^{(1)}_{s} p_\nu(X) -  L^{(2)}_{s} p_\nu(X)}_{X}^{-1 / (N+1)d}~.
\end{align*}
Similar to Case A, we want to find a lower bound for $\inorm{(d-1) L^{(1)}_{s} p_\nu(X) -  L^{(2)}_{s} p_\nu(X)}_{X}$ in terms of $r^{\epsilon_3}$, for some $1 \le s \le n-1$.
Since all linear operators between finite dimensional vector spaces are bounded, we have
\begin{align*}
    C' := \norm{(d-1) L^{(1)}_{s} - L^{(2)}_{s}}_{(\mathcal{P}^{n , d} , \norm{\cdot}_{\ell^1(\mathcal{Q})}) \to (\mathcal{P}^{n , (N+1)d},\inorm{\cdot}_{X})} < \infty~.
\end{align*}
The operator norm $C'$ is dependent only on $n,d$ and $\{p_1,..., p_N\}$ by construction of $L^{(1)}_s$, $L^{(2)}_s$, and in particular independent of $\nu$. 
By the triangle inequality, since $p_\nu = f_\nu + g_\nu$, for each $1 \le s \le n-1$ we have
\begin{align*}
  \inorm{(d-1) L^{(1)}_{s} p_\nu(X) -  L^{(2)}_{s} p_\nu(X)}_{X}  &= \inorm{(d-1) L^{(1)}_{s} g_\nu(X) -  L^{(2)}_{s} f_\nu(X) - L^{(2)}_{s} g_\nu(X)}_{X}\\
  & \ge \inorm{L^{(2)}_{s} f_\nu(X)}_{X} - \inorm{((d-1) L^{(1)}_{s} - L^{(2)}_{s})g_\nu(X)}_{X}\\
  & \ge  \inorm{L^{(2)}_{s} f_\nu(X)}_{X} - C'\norm{g_\nu }_{\ell^1(\mathcal{Q})}~.
\end{align*}
Fix $1 \le s_* \le n-1$ as in Lemma \ref{L2} to be the index such that there exists $C_{s_*}>0$ with
\begin{align*}
  \inorm{L^{(2)}_{s_*}p(X)}_{X} =C_{s_*}\norm{p}_{\ell^1(\mathcal{Q})}~, \mbox{ for all } p \in K~.
\end{align*} 
Since $f_\nu \in K$, we have $\inorm{L^{(2)}_{s_*} f_\nu(X)}_{X} = C_{s_*}\norm{f_\nu}_{\ell^1(\mathcal{Q})}$.
Moreover, since we are working in Case B, we may conclude via \eqref{l1norm} that
\begin{align*}
  r^{\epsilon_3}\le r \le \norm{\nu}  =\norm{p_\nu}_{\ell^1(\mathcal{Q})} & \le \norm{f_\nu}_{\ell^1(\mathcal{Q})} + \norm{g_\nu}_{\ell^1(\mathcal{Q})} \le \norm{f_\nu}_{\ell^1(\mathcal{Q})} + \epsilon r^{\epsilon_3}~,
\end{align*}
that is, $\norm{f_\nu}_{\ell^1(\mathcal{Q})} \ge (1-\epsilon) r^{\epsilon_3}$.  
Therefore, we have obtained
\begin{align*}
  \inorm{(d-1) L^{(1)}_{s_*} p_\nu(X) -  L^{(2)}_{s_*} p_\nu(X)}_{X}  & \ge C_{s_*}\norm{f_\nu}_{\ell^1(\mathcal{Q})} - C' \norm{g_\nu}_{\ell^1(\mathcal{Q})} \ge (C_{s_*} - \epsilon(C_{s_*} + C')) r^{\epsilon_3} ~.
\end{align*}
We now take $\epsilon > 0$ to be small enough so that the right hand side is positive. Again, since neither $C_{s_*}$ nor $C'$ depends on $\nu$, the choice of $\epsilon$ here is independent of $\nu$.

In summary, upon defining $\tilde{C} : = (C_{s_*} - \epsilon(C_{s_*} + C'))^{-1 / (N+1)d}$, we have
\begin{align*}
  \abs{G^\nu}
  & \le C r^{\epsilon_2 / (N+1)d}\inorm{(d-1) L^{(1)}_{s_*} p_\nu(X) -  L^{(2)}_{s_*} p_\nu(X)}_{X}^{-1 / (N+1)d} \le C \tilde{C} r^{(\epsilon_2 - \epsilon_3) / (N+1)d}~.
\end{align*}

\end{proof}

\begin{rem}
  Lemma \ref{L1} claims that the kernel $K = \ker (L^{(1)})$ is trivial when $d$ is odd. Hence, if $d$ is odd, the orthogonal decomposition $p_\nu = f_\nu + g_\nu$ simplifies to $p_\nu = g_\nu$, $f_\nu = 0$. Therefore, the analysis in Case B is in fact only relevant when $d$ is even.
\end{rem}

This concludes the construction of exceptional sets $G^\nu$, $F_u^\nu$ and the verification of the proposed estimates \eqref{setestimate} and \eqref{fullkernelbound} of Proposition \ref{mainprop}. All that remains to prove Proposition \ref{mainprop} is to verify Lemmas \ref{L1} and \ref{L2}.

\section{Proof of Lemma \ref{L1} and \ref{L2}} \label{lemma}
To prove Lemma \ref{L1} and \ref{L2}, we need to analyze the kernel of the operators $L^{(1)}_s, L^{(2)}_s: \mathcal{P}^{n , d} \to \mathcal{P}^{n , (N+1)d}$ in \eqref{L1andL2}.
We first state and prove a lemma about the determinant $D(X)$ from Lemma \ref{minorexp}, which will be used later in computing the kernels of $L^{(1)}_s$ and $L^{(2)}_s$.
\begin{lem}\label{det}
Let $\mathcal{Q}\subseteq \mathcal{P}^{n , d}$ be a $N: = \binom{d+n-1}{n-1}-n$ dimensional subspace with basis $\{p_1,..., p _N \}$. Assume that $\mathcal{Q}$ satisfies the nondegeneracy condition \eqref{cond1} for $c \in (\R \backslash \{0\})^n$.
Consider the determinant $D(X) \in \mathcal{P}^{n , Nd}$ from Lemma \ref{minorexp}. Then, $D(c) \neq 0$. Consequently, $D(X)$ is a nonzero polynomial.
\end{lem}
\begin{proof}
  Fix any constant $c \in \R^n$. 
By definition of $D(X)$ in Lemma \ref{minorexp}, $D(c)$ is the determinant of the matrix
\begin{align*}
  M(c) := \bracket{ \bracket{B_{j , \beta(i)}(c)}_{1 \le i , j \le N}}~.
\end{align*}
Recall the definition of $B_{i , \beta}(c)$ from \eqref{Bdef}.
Hence, $M(c)$
is the matrix representation of the coordinate projection map,
\begin{align}
  L(c): \mathrm{span}_\R\{p_i ~|~ 1 \le i \leq N \} & \longrightarrow \mathrm{span}_\R \{ \sigma^{\beta(i)}~~|~~ 1 \le i \le N\} \label{Lcdef}\\
  p_\lambda = \sum_{j = 1}^{N} \lambda_j p_j & \longmapsto \sum_{i = 1}^{N} \bracket{  \sum_{j = 1}^{N} \lambda_j B_{j,\beta(i)}(c) } \sigma^{\beta(i)}  =  \sum_{i = 1}^N\sbracket{\sigma^{\beta(i)}}\bracket{p_\lambda(\omega(\sigma , c))}\sigma^{\beta(i)} \nonumber~,
\end{align}
with respect to the bases explicitly given by the linear spans. Here, we have used the notation $[\sigma^\beta](f(\sigma))$ to denote the coefficient of $\sigma^\beta$ in the polynomial $f(\sigma)$.
Since the domain and codomain of $L(c)$ have the same dimension, to show that $D(c) = \det(M(c)) \neq 0$, it suffices to show that $L(c)$ is an injective map on $\mathcal{Q} = \mathrm{span}_\R\{p_i ~|~ 1 \le i \leq N \}$. We will show this for the constant $c \in (\R \backslash \{0\})^n$ in the tangential nondegeneracy condition, which satisfies
\begin{align}
  \mathcal{Q} \cap \mathrm{span}_\R \cbracket{ X_j(c \cdot X)^{d-1}~~|~~ 1 \le j \le n} = \{0\}~.\label{zerointersection}
\end{align}
By definition of $D(X)$, we have $\{\beta(i)\} = \{\beta \in \Z^{n-1}_{\ge 0} ~|~ 2 \le \abs{\beta} \le d\}$. So $L(c)$ is the restriction of the linear map 
\begin{align*}
  \tilde{L}(c) : \mathcal{P}^{n , d} &\longrightarrow \mathrm{span} \cbracket{ \sigma^\beta ~|~ 2 \le \abs{\beta} \le d} \\
  p &\longmapsto \sum_{2 \le \abs{\beta} \le d}\sbracket{\sigma^\beta}(p ( \omega(\sigma, c))) \sigma^\beta~
\end{align*}
to the linear subspace $\mathcal{Q} \subseteq\mathcal{P}^{n , d}$. Hence, $L(c)$ is injective if and only if $\ker (\tilde{L} (c)) \cap \mathcal{Q} = \{0\}$.

In the following, we will show  \begin{align*}
  \ker \bracket{\tilde{L} (c)} = \mathrm{span} \cbracket{ X_j\bracket{ c \cdot X}^{d-1}  ~~|~~ 1 \le j \le n }~,
\end{align*}
which immediately implies that $L(c)$ is injective by \eqref{zerointersection}.
We will prove double containment.
Suppose $p\in \ker(\tilde{L} (c))$. We will express $p$ in terms of the basis
\begin{align*}
    \cbracket{ X_j\bracket{ c \cdot X}^{d-1}  ~~|~~ 1 \le j \le n }~.
\end{align*}
First, $\tilde{L}$ is the linear map which takes a polynomial $p \in \mathcal{P}^{ n ,d}$, maps it to a polynomial in $\sigma$ via the composition $p \circ \omega(\sigma,c)$, and then eliminates all monomials $\sigma^\beta$ with $0 \le \abs{\beta} \le 1$.
Hence, by definition of $\tilde{L}(c)$, 
\begin{align*}
    p(\omega(\sigma , c)) = l_0(\sigma) + c_0~,
\end{align*}
where $c_0 \in \R$ and $l_0 \in \mathcal{P}^{n-1,1}$. 
Consider the hyperplane in $\R^n$ defined by
\begin{align*}
    H : = \cbracket{ x \in \R^n ~~|~~ c \cdot x= \abs{c}^2}~.
\end{align*}
This hyperplane has parametrization $H = \cbracket{ \omega(\sigma , c)~~|~~ \sigma \in \R^{n-1}}$.
To prove this, observe that for each $\sigma \in \R^{n-1}$, we have
\begin{align*}
  c\cdot \omega(\sigma,c) = \bracket{ \bracket{ \sum_{i = 1}^{n-1} c_i(c_i - c_n \sigma_i) }+ c_n^2 + c_n \sum_{i = 1}^{n-1} c_i\sigma_i} = \abs{c}^2~.
\end{align*}
Moreover, if $x \in H$, then taking $\sigma = (c^{(n)} - x^{(n)})/c_n $ gives $x = (x_1 , ... ,x_n) = \omega\bracket{\sigma, c}$. The assumption that $c \in (\R \backslash \{0\})^n$ in particular ensures that this operation is well-defined. Hence, for each $x \in H$, there is $\sigma \in \R^{n-1}$ such that
\begin{align}
  p(x) = p(\omega(\sigma , c)) = c_0 + l_0 \bracket{\sigma} = c_1 + l_1 \bracket{x^{(n)}} ~,\label{pxonhu0}
\end{align}
where $c_1 \in \R$ is a constant and $l_1 \in \mathcal{P}^{n-1,1}$. 
Using the fact that $p(x)$ is homogeneous of degree $d$, for any $x \in \R^n$ such that $c\cdot x\neq 0$, we can write
\begin{align*}
  x = \bracket{\frac{c\cdot x}{\abs{c}^2}} y~, ~ \mbox{ where } y := \frac{\abs{c}^2}{c\cdot x} x \in H~.
\end{align*}
Using the expression of $p$ in \eqref{pxonhu0} and the fact that $p$ is homogeneous of degree $d$, for all $x \in \R^n$ such that $c\cdot x\neq 0$, we can write
\begin{align*}
  p(x) = \bracket{\frac{c\cdot x}{\abs{c}^2}}^d p(y)  &=  \bracket{\frac{c\cdot x}{\abs{c}^2}}^d  \bracket{ c_1 + l_1\bracket{y^{(n)}}} =  c_1 \bracket{\frac{c\cdot x}{\abs{c}^2}}^d + l_1 \bracket{x^{(n)}}\bracket{\frac{c\cdot x}{\abs{c}^2}}^{d-1}~.
\end{align*}
It follows that the polynomial
\begin{align*}
  p(X) - c_1 \bracket{\frac{c\cdot X}{\abs{c}^2}}^d - l_1 \bracket{X^{(n)}}\bracket{\frac{c\cdot X}{\abs{c}^2}}^{d-1} \in \mathcal{P}^{n , d}
\end{align*}
is a polynomial that vanishes on the nonempty open set $\cbracket{x \in \R^n~~|~~ c \cdot x \neq 0}$. This polynomial is the zero polynomial since a polynomial vanishing on an open set vanishes identically. Therefore, we have
\begin{align*}
  p(X)= c_1 \bracket{\frac{c\cdot X}{\abs{c}^2}}^d + l_1 \bracket{X^{(n)}}\bracket{\frac{c\cdot X}{\abs{c}^2}}^{d-1} \in  \mathrm{span}_\R \cbracket{ X_j\bracket{c\cdot X}^{d-1}  ~~|~~ 1 \le j \le n }~.
\end{align*}

For the reverse inclusion, assume $p(X) = l(X)\bracket{ \frac{c\cdot X}{\abs{c}^2}}^{d-1}$ for some $l \in \mathcal{P}^{n , 1}$. We will show explicitly that $\tilde{L}(c) p $ is the zero polynomial.
Since $c \cdot \omega(\sigma,c) = |c|^2$, we have
\begin{align*}
  p(\omega(\sigma , c)) = l(\omega(\sigma , c))\bracket{ \frac{c\cdot \omega(\sigma , c)}{\abs{c}^2}}^{d-1} = l(\omega(\sigma , c))~.
\end{align*}
Since $l$ is homogeneous of degree $1$, $[\sigma^\beta](l(\omega(\sigma , c))) = 0$ for all $2 \le  \abs{\beta} \le d$. It follows that $p$ is in the kernel of $ \tilde{L} (c)$.

\end{proof}

\subsection{Proof of Lemma \ref{L1}}\label{lemmaL1proof}
Recall the definition of $L^{(1)}$,
\begin{align*}
    L^{(1)}: \mathcal{P}^{n,d} \longrightarrow \bracket{\mathcal{P}^{n,(N+1)d}}^{n-1}, \quad p \longmapsto (D K^{e_1} p , ... , D K^{e_{n-1}} p)~,
\end{align*}
where $D \in \mathcal{P}^{n,Nd}$ and the maps $K^{e_s}: \mathcal{P}^{n,d} \to \mathcal{P}^{n,d}$ are defined as in Lemma \ref{minorexp}.

We will prove double containment. First, assume that $p \in \ker (L^{(1)})$. The strategy is to show that the function $p : \R^n \to \R$ given by evaluating the polynomial $p$ at points in $\R^n$ is radial, i.e., there is a function $f:[0,\infty) \to \R$ such that $p(u) = f(|u|)$. Indeed, if this is true, then since all radial polynomials in $\R[X_1,...,X_n]$ are of the form $q(|X|^2)$ for some polynomial $q \in \R[t]$, we have
\begin{align*}
    \ker (L^{(1)}) \subseteq \R[\abs{X}^2]\cap \mathcal{P}^{n,d} = \left\{\begin{array}{cl}
    \mathrm{span}_\R \cbracket{ \abs{X}^d }  & \mbox{ if $d$ is even }  \\ [5pt]
    \{0\} & \mbox{ if $d$ is odd }
  \end{array}\right. \, ~.
\end{align*}
If $L^{(1)}p = 0$, then $D K^{e_s}p = 0$ as a polynomial for each $1 \le s \le n-1$. Since the polynomial ring $\R[X]$ is an integral domain and Lemma \ref{det} has shown that $D (X)$ is not the zero polynomial, $K^{e_s}p = 0$ as a polynomial for each $1 \le s \le n-1$. It follows that for each $u \in \R^n$, $1\le s \le n-1$, one has
\begin{align*}
  0 = K^{e_s} p (u) =  [\sigma_s]p(\omega(\sigma, u)) & = \frac{\partial}{\partial \sigma_s}\bracket{p(\omega(\sigma, u))}\bigg|_{\sigma = 0}\\
   &= \nabla p(\omega(0, u)) \cdot \frac{\partial}{\partial \sigma_s}\omega(\sigma, u)\bigg|_{\sigma = 0}  = \bracket{ u_s \frac{\partial}{\partial {x_n}} - u_n\frac{\partial}{\partial{x_s}}}p(u)~.
\end{align*}
Now define, for each angle $\theta \in \R$, a rotation along the $(u_s , u_n)$ plane by $\theta$ degrees,
\begin{align*}
  R_s(\theta)&: \R^n \longrightarrow \R^n\\
  R_s(\theta)u & = (u_1 , ... ,u_{s-1} , u_s\cos(\theta) - u_n \sin(\theta) , u_{s+1} , ..., u_{n-1} , u_s\sin(\theta) + u_n \cos{\theta})~.
\end{align*}
Then $\bracket{\frac{\dd}{\dd \theta}R_s(\theta)u}_{s} = -(R_s(\theta)u)_n$ and $\bracket{\frac{\dd}{\dd \theta}R_s(\theta)u}_{n} = (R_s(\theta)u)_s$, and $\bracket{\frac{\dd}{\dd \theta}R_s(\theta)u}_{k} = 0$ for all $k \neq s, n$.
Hence, for each $u \in \R^n$, by the chain rule, we get
\begin{align*}
  \frac{\dd}{\dd \theta} p \bracket{R_s(\theta)u}  &= \nabla p \bracket{R_s(\theta)u} \cdot  \frac{\dd}{\dd \theta}R_s(\theta)u\\
  &= \frac{\partial}{\partial u_s}p \bracket{R_s(\theta)u} \bracket{\frac{\dd}{\dd \theta}R_s(\theta)u}_s + \frac{\partial}{\partial u_n}p \bracket{R_s(\theta)u} \bracket{\frac{\dd}{\dd \theta}R_s(\theta)u}_n\\
   &= -(R_s(\theta)u)_n \frac{\partial}{\partial u_s}p \bracket{R_s(\theta)u} + (R_s(\theta)u)_s \frac{\partial}{\partial u_n}p \bracket{R_s(\theta)u}\\
    &= 0~.
\end{align*}
Hence, for every $u \in \R^n$ and $1 \le s \le n-1$, the function $\theta \mapsto p \bracket{R_s(\theta)u}$ is constant in $\theta$. In particular, this means
\begin{align*}
  p \bracket{R_s(\theta)u} = p \bracket{R_s(0)u} = p(u)~\mbox{ , for all } \theta \in \R ~,~ 1 \le s \le n-1~.
\end{align*}

Given an arbitrary point $u = (u_1 ,..., u_n) \in \R^n \backslash \{ 0\}$, there is a sequence of angles $\theta_1 , ..., \theta_{n-1} \in \R$ such that
\begin{align*}
  R_{n-1}(\theta_{n-1})\cdots R_1(\theta_1)\frac{u}{\abs{u}} = e_n~.
\end{align*}
As a result, using that $p(u)$ is homogeneous of degree $d$, we can write
\begin{align*}
  p(u) = \abs{u}^d p \bracket{ \frac{u}{\abs{u}}} = \abs{u}^dp \bracket{R_{n-1}(\theta_{n-1})\cdots R_1(\theta_1)\frac{u}{\abs{u}}}  = p \bracket{ 0 , ... , 0 , \abs{u}}~. 
\end{align*}
Therefore, we have shown that $p: \R^n \to \R$ is a radial function. 

For the reverse containment, suppose $p(X) = \abs{X}^d$ for some even integer $d \ge 2$. Then for each $ u \in \R^n$ and $1 \le s \le n-1$, we have
\begin{align*}
    K^{e_s} p (u) = [\sigma_s]p(\omega(\sigma, u))  = \bracket{ u_s \frac{\partial}{\partial {x_n}} - u_n\frac{\partial}{\partial{x_s}}}p(u) = (u_s u_n - u_n u_s) d \abs{u}^{d-2} = 0~.
\end{align*}
Therefore, $K^{e_s} p$ is the zero polynomial for all $1 \le s \le n-1$, and consequently $L^{(1)}p = 0$.

\subsection{Proof of Lemma \ref{L2}} \label{lemma3}
Recall the definition of $L^{(2)}_s: \mathcal{P}^{n , d} \to \mathcal{P}^{n , (N+1)d}$,
\begin{align*}
  L^{(2)}_sp  = \sum_{j : \abs{\beta(j)} = 2} D_j^{e_s} K^{\beta(j)} p~,
\end{align*}
where $D_j^{e_s} \in \mathcal{P}^{n , Nd}$ is the polynomial
\begin{align*}
  D_j^{e_s}(X)= \sum_{i = 1}^{N} (-1)^{i+j} B_{i, e_s}(X) D_{j,i}(X)~,
\end{align*}
and $D_{j,i} \in \mathcal{P}^{n,(N-1)d}$ and the maps $K^{e_s}: \mathcal{P}^{n,d} \to \mathcal{P}^{n,d}$ are defined as in Lemma \ref{minorexp}.

It follows from Lemma \ref{L1} that the kernel of $L^{(1)}$ is trivial when $d$ is odd. Thus, it suffices to prove Lemma \ref{L2} in the case when $d$ is even. The result of Lemma \ref{L2} can be derived as a direct consequence of the following lemma.
\begin{lem}\label{L22}
Let $d\ge 2$ be even. Let $\mathcal{Q} \subseteq \mathcal{P}^{n , d}$ be a codimension $n$ linear subspace, equipped with a basis $\mathfrak{B}$. Suppose there is $c \in (\R \backslash \{0\})^n$ such that $\mathcal{Q}$ simultaneously satisfies the tangential and radial nondegeneracy conditions for $c$. Then, there is $1 \le s \le n-1$, depending only on $\mathcal{Q}$, such that $L^{(2)}_s (\abs{X}^d) \neq 0$ as a polynomial.
\end{lem}

\begin{proof}[Proof of Lemma \ref{L2}]
    Let $1 \le s \le n-1$ and define $p(X) := \abs{X}^d$. If $q\in \ker \bracket{L^{(1)}}$, then by Lemma \ref{L1}, $q = Cp$ for some $C  \in \R$. Since $p \neq 0$ as a polynomial and $\norm{\cdot}_{\ell^1(\mathcal{Q})}$ is a norm, we can write
    \begin{align*}
        \abs{C} = \frac{\norm{q}_{\ell^1(\mathcal{Q})}}{\norm{p}_{\ell^1(\mathcal{Q})}}~.
    \end{align*}
Since $\inorm{\cdot}_X$ is a norm,
\begin{align*}
  \inorm{L^{(2)}_s q(X)}_{X} = \abs{C} \inorm{L^{(2)}_sp(X)}_{X} = \inorm{L^{(2)}_sp(X)}_{X} (\norm{p}_{\ell^1(\mathcal{Q})})^{-1}  \norm{q}_{\ell^1(\mathcal{Q})}~.  
\end{align*}
Assuming Lemma \ref{L22}, there is an index $1 \le s_* \le n-1$, depending only on the construction of $L^{(2)}_s$ (and hence only on $n , d$ and the basis $\mathfrak{B} $ of $\mathcal{Q}$) such that $L^{(2)}_{s_*}p \neq 0$ as a polynomial. Then we have
\begin{align*}
    C_{s_*}: = \inorm{L^{(2)}_{s_*}p(X)}_{X} (\norm{p}_{\ell^1(\mathcal{Q})})^{-1}  >  0~,
\end{align*}
and so $\inorm{L^{(2)}_{s_*}p(X)}_{X}  = C_{s_*} \norm{p}_{\ell^1(\mathcal{Q})}$.
In particular, this constant $C_{s_*}$ is dependent only on $n , d$ and the basis $\mathfrak{B}$ of $\mathcal{Q}$.

\end{proof}

The rest of this section will be dedicated to the proof of Lemma \ref{L22}.
\begin{proof}[Proof of Lemma \ref{L22}]
Let $d\ge 2$ be even and $p(X) : = \abs{X}^d$. The strategy is to show that if $L^{(2)}_sp$ were the zero polynomial for all $1 \le s \le n-1$, then there would be a contradiction with the assumed nondegeneracy conditions.

Let $c \in (\R\backslash \{0\})^n$ be the constant which satisfies the tangential \eqref{cond1} and radial \eqref{cond2} nondegeneracy conditions in the assumption of Lemma \ref{L22}. We first compute an explicit expression for $L^{(2)}_sp(c)$. Recall the definition of $L^{(2)}_s$ at the start of this subsection.
Since $D_{j , i}(X)$ is the $ji$-th minor of the matrix
\begin{align*}
  M(X): = \bracket{B_{j, \beta(i)}\bracket{X}}_{1 \le i,j\le N}~,
\end{align*}
the polynomial
$D_j^{e_s}(X)$ is exactly the determinant of the matrix obtained by replacing the $j$-th row of $M(X)$ with the row vector $\bracket{B_{i, e_s}(X)}_{1 \le i \le N}$.
Define a column vector $(C^{\beta(j)}(X))_{1 \le j \le N}$ in the following way,
\begin{align*} 
  C^{\beta(j)}(X)   = \left\{\begin{array}{cl}
    K^{\beta(j)}p(X)  & \mbox{ if $\abs{\beta(j)} = 2$}  \\ [5pt]
    0 & \mbox{ if otherwise}
  \end{array}\right. \, ~.
\end{align*} 
Then $L^{(2)}_sp$ is exactly $-1$ times the determinant of the following $(N+1) \times (N+1)$ matrix
\[
\widetilde M_s(X)
: =
\begin{bmatrix}
0 & \bigl(B_{j,e_s}(X)\bigr)_{1\le j\le N} \\[0.4em]
C^{\beta(1)} \bracket{X} &  \\[-0.2em]
\vdots & M(X) \\
C^{\beta(N)}(X) &
\end{bmatrix}~.
\]
This can be obtained by performing a minor expansion along the first column.
Therefore, to show that there is some $1 \le s \le n-1$ such that $L^{(2)}_sp$ is not the zero polynomial, it suffices to show that $\det (\widetilde M_s(c)) \neq 0$ for some $c \in \R^n$, $1 \le s \le n-1$. We will show this for the constant $c \in (\R\backslash \{0\})^n$ which satisfies the tangential \eqref{cond1} and radial \eqref{cond2} nondegeneracy conditions.

By Lemma \ref{det}, for the chosen $c$ which satisfies the tangential \eqref{cond1} nondegeneracy condition, we have
\begin{align*}
    D(c) = \det(M(c)) \neq 0~.
\end{align*}
Hence, the inverse matrix $M^{-1}(c)$ is well-defined and we can use Schur's complement formula (see \cite[Page 364]{MR255573}) to write
\begin{align*}
  \det (\widetilde{M}_s(c)) = -\det(M(c))\det \bracket{\bigl(B_{j,e_s}(c)\bigr)_{1\le j\le N}  \cdot M^{-1}(c)\bigl(C^{\beta(j)}(c)\bigr)_{1\le j\le N} } ~.
\end{align*}
Recall that
$M(c)$ is the real-valued matrix representing the coordinate projection map $L(c)$ defined in \eqref{Lcdef}. Therefore, $M^{-1}(c)$ is the matrix representing the map
\begin{align*}
    L^{-1}(c) : \mathrm{span}_\R\{\sigma^{\beta(i)}~~|~~ 1\le i \le N \} \to   \mathrm{span}_\R \{ p_1,...,p_N\} = \mathcal{Q}
\end{align*}
with respect to the same bases. Meanwhile, by construction of the maps $K^\beta p(X) = [\sigma^\beta]p(\omega(\sigma, X))$, $(C^{\beta(j)}(c))_{1 \le j \le N} \in \R^N$ is the vector representing the nonzero polynomial
\begin{align*}
  q_c(\sigma) = \sum_{\abs{\beta} = 2}^{} [\sigma^\beta]\bracket{ \abs{\omega(\sigma, c )}^d} \sigma^\beta~ \in \mathrm{span}\{ \sigma^{\beta(i)}~~|~~ 1 \le i \le N\}~,
\end{align*}
in coordinates explicitly given by the basis $\{\sigma^{\beta(i)}\}_{1\le i \le N}$.
We expand $\abs{\omega(\sigma, c)}^d$ using the binomial theorem and compute explicitly
\begin{align*}
  q_c(\sigma) = \frac{d}{2} \abs{c}^{d-2}\bracket{\abs{\omega(\sigma, c)}^2 - \abs{c}^2}~. 
\end{align*}
Therefore, the vector
\begin{align*}
  (\lambda_j)_{1 \le j \le N} : =   M^{-1}(c)\bracket{C^{\beta(j)}(c)}_{1 \le j \le N} \in \R^N
\end{align*}
is the vector representing the polynomial $P_c : = L^{-1}(c)q_c \in \mathcal{Q}$ in coordinates given by the basis $\{p_1,..., p_N\}$ of $\mathcal{Q}$. We write $P_c = \sum_{ j = 1}^N \lambda_j p_j$. Then, since $B_{i, \beta}(X) = [\sigma^\beta]p_i(\omega(\sigma , X))$, we have
\begin{align*}
  L^{(2)}_sp (c)& = -D(c)\bracket{B_{j , e_s}(c)}_{1 \le j \le N} \cdot M^{-1}(c)\bracket{C^{\beta(j)}(c)}_{1 \le j \le N}\\ & = -D(c)\sum_{j = 1}^{N} B_{j , e_s }(c) \lambda_j \\
   &= - D(c)[\sigma^{e_s}] \bracket{\sum_{j = 1}^N \lambda_j p_j(\omega(\sigma , c )) } \\
  &=- D(c)[\sigma^{e_s}] P_c(\omega(\sigma , c ))~.
\end{align*}

Assume for contradiction that $L^{(2)}_sp (c) = 0$ for all $1 \le s \le n-1$. Then, by the statement from Lemma \ref{det} that $D(c) \neq 0$, we must have $[\sigma^\beta]P_c(\omega(\sigma , c )) = 0$ for all $\abs{\beta} =  1$. Since we have $L(c)P_c = q_c$ and $q_c$ contains no monomials in $\sigma$ of degree strictly larger than $2$, we must also have $[\sigma^\beta]P_c(\omega(\sigma , c )) = 0$ for all $\abs{\beta} \ge 3$ by construction of $L(c)$. Hence, $P_c(\omega(\sigma , c ))$ (considered as a polynomial in $\sigma$) must take the following form
\begin{align*}
P_c(\omega(\sigma, c)) = A(c) + \frac{d}{2} \abs{c}^{d-2}\bracket{\abs{\omega(\sigma, c)}^2 - \abs{c}^2}~,
\end{align*}
for some constant $A(c)$ that depends only on $c$.

We claim that
\begin{align}
  P_c(X) = \bracket{A(c) - \frac{d}{2} \abs{c}^d}\bracket{\frac{c \cdot X}{\abs{c}^2}}^{d}  + \frac{d}{2}|c|^{d-2} \bracket{\frac{c \cdot X}{\abs{c}^2}}^{d-2}\abs{X}^2~. \label{Pc}
\end{align} If this is true, then by the radial nondegeneracy condition \eqref{cond2}
\begin{align*}
    \mathcal{Q}\cap \mathrm{span}_\R \{(c\cdot X)^d~,~ |X|^2(c \cdot X)^{d-2}\} = \{0\}~,
\end{align*}
we must have $P_c= 0$ as a polynomial in $\R[X_1,...,X_n]$. However, this is not possible as the two polynomials of the sum in \eqref{Pc} are linearly independent in $\R[X_1,...,X_n]$, and we have a nonzero linear combination
\begin{align*}
    P_c(X) = a_1 \bracket{c \cdot X}^{d} + a_2 \bracket{c \cdot X}^{d-2}\abs{X}^2~,~ a_2 \neq 0~.
\end{align*}The assumption that $c \in (\R \backslash \{0\})^n$ in particular ensures $a_2 \neq 0$.
Hence, we have arrived at a contradiction, and we may conclude that $L^{(2)}_sp (c) \neq 0$ for some $1 \le s \le n-1$.

It remains to show \eqref{Pc}, which uses a similar argument as the proof of Lemma \ref{det}.
Consider again the hyperplane $H : = \{ x \in \R^n ~~|~~ c \cdot x = \abs{c}^2\}$. 
We have the parametrization $H = \cbracket{\omega(\sigma , c) ~~|~~ \sigma \in \R^{n-1}}$.
For each $x \in H$, we can choose $\sigma = (c^{(n)} - x^{(n)})/c_n$ so that $x = \omega(\sigma , c)$ and
\begin{align*}
  P_c(x) = P_c \bracket{ \omega(\sigma , c)}  &=  A(c) + \frac{d}{2} \abs{c}^{d-2}\bracket{\abs{x}^2 - \abs{c}^2}~.
\end{align*}

Since $P_c$ is homogeneous of degree $d$, on the set $\cbracket{x \in \R^n ~~|~~ c \cdot x  \neq 0}$, we can write
\begin{align*}
  x = \bracket{\frac{c\cdot x}{\abs{c}^2}} y~, ~ \mbox{ where } y := \frac{\abs{c}^2}{c\cdot x} x \in H~.
\end{align*}
Thus, if $c\cdot x \neq 0$, then by homogeneity of $P_c$, we have
\begin{align*}
  P_c(x) = \bracket{\frac{c \cdot x}{\abs{c}^2}}^{d} P_c(y) &=   \bracket{\frac{c \cdot x}{\abs{c}^2}}^{d}\bracket{ A(c)- \frac{d}{2} \abs{c}^d + \frac{d}{2} \abs{c}^{d-2}\abs{y}^2 } \\
  & = \bracket{\frac{c \cdot x}{\abs{c}^2}}^{d} \bracket{A(c) - \frac{d}{2} \abs{c}^d} + \bracket{\frac{c \cdot x}{\abs{c}^2}}^{d-2}\abs{x}^2|c|^{d-2}\frac{d}{2} ~.
\end{align*}
Again, it follows that
\begin{align*}
  P_c(X) - \bracket{\frac{c \cdot X}{\abs{c}^2}}^{d} \bracket{A(c) - \frac{d}{2} \abs{c}^d} - \bracket{\frac{c \cdot X}{\abs{c}^2}}^{d-2}\abs{X}^2\frac{d}{2} |c|^{d-2} \in \mathcal{P}^{n , d}
\end{align*}
is a polynomial which vanishes on the nonempty open set $\cbracket{x \in \R^n~~|~~ c \cdot x \neq 0}$. Hence, this polynomial must vanish identically.

\end{proof}

\section{Proof of Proposition \ref{derivationofmainthm}}\label{subspacesec}

In this section, we prove Proposition \ref{derivationofmainthm}, which implies the main result in Theorem \ref{mainresult2}. 

\subsection{Required Lemmas}To prove Proposition \ref{derivationofmainthm}, we need the following lemmas, which allow us to deduce injectivity of maps using estimates on the dimension of subspaces of $\mathcal{P}^{n , d}$. 
\begin{lem}\label{lemma:subspace} Let $n , d \ge 2$ be integers. Suppose $m > \binom{d+n-2}{n-2}$, so in particular $m\ge n$. Let $\mathcal{M} \subseteq \mathcal{P}^{n,d}$ be a $m$ dimensional linear subspace. Then, $\mathcal{M}$ admits an orthonormal basis $\{f_1,...,f_m\}$ with respect to the Fischer inner product on $\mathcal{P}^{n,d}$ defined in \eqref{innerprod}, such that there exists a set of indices $\{i_1, ... , i_n\} \subseteq \{1,...,m\}$ and an open dense set $U \subseteq \R^n$ such that all $c \in U$ satisfy
\begin{align*}
  \det \sbracket{  \begin{matrix}
    \nabla f_{i_1}(c)  \\
      \vdots \\
     \nabla f_{i_n}(c)
    \end{matrix}}\neq 0~.
\end{align*}
Consequently, the $m \times n$ matrix defined for $c \in \R^n$,
\begin{align}
  F(c) : = \sbracket{  \begin{matrix}
    \nabla f_1(c)  \\
      \vdots \\
     \nabla f_m(c)
    \end{matrix}}~\label{matrixrep2}
\end{align}
has rank $n$ for all $c \in U$.
\end{lem}

\begin{lem} \label{laplaciamlemma}Fix $c \in \R^n$, $n , d \ge 2$, and $m \ge n$.
  Let $f_1, ..., f_m \in \mathcal{P}^{n , d}$ be nonzero polynomials such that the $m \times n$ matrix $F(c)$ in \eqref{matrixrep2}
  has rank $n$, $f_i(c) \neq 0$ for all $1 \le i \le m$, and $\Delta f_1 (c)\neq 0$, where $\Delta = \sum_{i =  1}^{n} \partial_i^2$ is the Laplacian operator. Moreover, assume that for all multi-indices $1 \le i_1 < ... < i_n \le m$ such that the $n \times n$ matrix
  \begin{align*}
      A(c, i_1,...,i_n) : = (d-1)!\sbracket{  \begin{matrix}
    \nabla f_{i_1}(c)  \\
      \vdots \\
     \nabla f_{i_n}(c)
    \end{matrix}}
  \end{align*}
  has nonzero determinant, the $n \times 2$ matrix
  \begin{align*}
    B(c, i_1, ..., i_n) : = \sbracket{ \begin{matrix} d! f_{i_1} (c) &  (d-2)!\Delta f_{i_1}(c) \\ 
  \vdots & \vdots \\ 
  d!f_{i_n} (c) &  (d-2)!\Delta f_{i_n}(c) \end{matrix}}
  \end{align*}
has rank strictly less than $2$.
Then, there exists a nonzero constant $\lambda(c) \in \R \backslash \{0\}$ such that 
  \begin{align*}
      \Delta f_i(c) = \lambda(c) f_i(c) ~,~ \mbox{ for all } 1 \le i \le m~.
  \end{align*}
  
\end{lem}

\begin{proof}[Proof of Lemma \ref{lemma:subspace}]

Fix a monomial order on the multi-indices $\alpha \in \Z^{n}_{\ge 0}$, that is, a total order on $\Z^{n}_{\ge 0}$ such that if $\alpha > \beta$, then $\alpha + \gamma > \beta + \gamma$. By Gaussian elimination and the Gram-Schmidt process, we can construct an orthonormal basis $\{f_1 ,..., f_m\}$ of $\mathcal{M}$ such that each $f_i$ has a distinct leading monomial $X^{\alpha(i)}$ among the leading monomials of all the $f_j$'s with respect to the chosen monomial order. Let $\{\alpha(i)~~|~~ 1 \le i \le m\}$ be the set of leading multi-indices.

We claim that if $m > \binom{d+n-2}{n-2}$, then the set $\{\alpha(i)~~|~~ 1 \le i \le m\}$ contains a linearly independent subset of size $n$. Assuming that this is true for now, we proceed with the proof.
Without loss of generality, we reindex the multi-indices corresponding to the leading monomials so that $\{\alpha(i)~~|~~ 1 \le i \le n\}$ is a linearly independent set in $\R^n$. We will show
\begin{align}
  \det \sbracket{\begin{matrix}
    \nabla f_1(X)  \\
      \vdots \\
     \nabla f_n(X)
    \end{matrix}} \neq 0 \label{minor}
\end{align}
as a polynomial in $X$. To prove that a polynomial is nonzero, it suffices to find a monomial with a nonzero coefficient. Under the chosen monomial order, we write
\begin{align*}
    f_i(X) = c_iX^{\alpha(i)} + \sum_{\beta < \alpha(i)}c_{i,\beta} X^{\beta}~,~ \mbox{ for } c_i \neq 0~.
\end{align*}
Then, 
\begin{align*}
  \bracket{X_1 \cdots X_n}\det \bracket{\bracket{ \partial_jf_i(X)}_{1 \le i , j \le n}} = &  \det \bracket{\bracket{ X_j\partial_jf_i(X)}_{1 \le i , j \le n}}\\ = &  \sum_{\sigma \in S_n}^{}\mathrm{sgn}(\sigma) \prod_{i = 1}^{n} X_{\sigma(i)}\partial_{\sigma(i)} f_i(X) \\ =&\sum_{\sigma \in S_n}^{}\mathrm{sgn}(\sigma) \prod_{i = 1}^{n} \bracket{c_i \alpha(i)_{\sigma(i)}X^{\alpha(i)} + \sum_{\beta < \alpha(i)}^{} c_{i , \beta} \beta_{\sigma(i)} X^{\beta}}\\
    =& X^{\sum_{i = 1}^n \alpha(i) }\det \bracket{ \bracket{\alpha(i)_j}_{1 \le i , j \le n}} \prod_{i = 1}^n c_i + \mbox{ lower order terms}~.
\end{align*}
The leading term on the right hand side has nonzero coefficient since $\{\alpha(i)~~|~~ 1 \le i \le n\}$ is a linearly independent set in $\R^n$, so the polynomial on the left hand side is nonzero. Since the polynomial ring $\R[X]$ is an integral domain, the determinant in \eqref{minor} is a nonzero polynomial in $X$. The complement of the zero set of a nonzero polynomial is an open dense set $U \subseteq \R^n$. Therefore, we have
\begin{align*}
    \det \sbracket{\begin{matrix}
    \nabla f_1(c)  \\
      \vdots \\
     \nabla f_n(c)
    \end{matrix}} \neq 0 ~,~ \mbox{ for all } c \in U~.
\end{align*}

It remains to show that if $m > \binom{d+n-2}{n-2}$, then the set of multi-indices
\begin{align*}
    S: = \{\alpha(i)~~|~~ 1 \le i \le m\} \subseteq \{\alpha \in \Z_{\ge 0}^n~~|~~ |\alpha| = d\}
\end{align*}
contains a linearly independent subset of size $n$ in $\R^n$.   
We prove the contrapositive. Suppose $\mathrm{span}_\R S$
has dimension strictly less than $n$. Then there exists a nonzero vector $c \in \R^n$ such that $c \cdot \alpha = 0$ for all $\alpha \in S$. 
  For this $c$, there must exist coordinates $c_i , c_j$ such that $c_i \neq c_j$, since if not, then there is $\lambda \neq 0$ such that $c \cdot \alpha = \lambda\abs{\alpha} = \lambda d \neq 0$ for $\alpha \in S$. Without loss of generality, assume $c_{n-1} \neq c_n$. Then each $\alpha \in S$ satisfies the system of equations
    \begin{align*}
        \alpha_{n-1} + \alpha_n & = d - \sum_{ i =1}^{n-2} \alpha_i  \\ 
        c_{n-1}\alpha_{n-1} + c_n\alpha_n & = - \sum_{ i =1}^{n-2}c_i\alpha_i~.
    \end{align*}
    We bound the maximum number of multi-indices in $\{\alpha \in \Z_{\ge 0}^n~~|~~ |\alpha| = d\}$ for which this system holds.
    Treating $(\alpha_{n-1}, \alpha_n)$ as unknown variables, the determinant of the system is $c_{n} - c_{n-1} \neq 0$, so there is a unique solution $(\alpha_{n-1}, \alpha_n) \in \R^2$ for each choice of $(\alpha_1,..., \alpha_{n-2}) \in \R^{n-2}$ which determines the right hand side. Since $\alpha$ has to satisfy $|\alpha| = d$ and have integer entries, the possible choices of $(\alpha_1,..., \alpha_{n-2})$ range in the set
    \begin{align*}
       \cbracket{\gamma \in \Z_{\ge 0}^{n-2}~~|~~ 0 \le |\gamma| \le d}~~,
    \end{align*}
    which has cardinality $\binom{d+n-2}{n-2}$.
    Therefore, we have $ \abs{S} \le \binom{d+n-2}{n-2}$.

\end{proof}

\begin{proof}[Proof of Lemma \ref{laplaciamlemma}]
Fix $c \in \R^n$. Since $f_i$ is homogeneous of degree $d$, we have
\begin{align*}
  0 \neq d f_i(c) = c \cdot \nabla f_i(c)~.
\end{align*}
Therefore, all rows $\nabla f_i(c)$ of $F(c)$ are nonzero. For any indices $1 \le i_1 < ... < i_n \le m$, the matrix $A(c, i_1,...,i_n)$ is a submatrix of $(d-1)!F(c)$. 
Since $F(c)$ has rank $n$ and the first row $\nabla f_1(c)$ is nonzero, there are indices $2 \le i_2, ..., i_n \le m$ such that
\begin{align*}
    \det(A(c, 1, i_2,...,i_n)) \neq 0~.
\end{align*}
Without loss of generality, we reindex so that $i_j = j$ for $2 \le j \le n$. By assumption of the rank of $B(c, 1, 2, ... , n)$ being strictly less than $2$, we must have that the two columns are linearly dependent.
Since $\Delta f_1(c) \neq 0$ and $f_i(c) \neq 0$ for all $1 \le i \le m$, both columns of the matrix $B(c, 1, 2, ... , n)$ are nonzero, so there is a nonzero real number $\lambda(c) \in \R \backslash \{0\}$ such that
\begin{align}
    \lambda(c) f_{i}(c) = \Delta f_{i} (c)~,~ \mbox{ for all } 1 \le i \le n~. \label{lambdaeq}
\end{align}
We claim that \eqref{lambdaeq} actually holds for all $1 \le i \le m$.
Fix $n+1 \le k \le m$, and consider the $(n+1) \times n$ submatrix of $F(c)$, defined by
\begin{align*}
  \sbracket{\begin{matrix}
     \nabla f_1(c)  \\
      \vdots \\
      \nabla f_n(c) \\
      \nabla f_k(c)
    \end{matrix}}~.
  \end{align*}This submatrix has rank $n$ since it contains a nonvanishing $n \times n$ minor given by the first $n$ rows.
 By the rank condition of the $(n+1) \times n$ submatrix, there exists a subset $\{ j_2, ..., j_n\} \subseteq \{1 , ... , n\}$ such that
\begin{align*}
    \det(A(c , k, j_2 , ..., j_n)) \neq 0~.
\end{align*}
Again, since $f_i(c) \neq 0$ for all $1 \le i \le m$, the first column of the matrix $B(c , k , j_2, ..., j _n)$ is nonzero,
so there is a real number $\lambda^{(k)}(c) \in \R$, such that
\begin{align*}
  \lambda^{(k)}(c) f_i(c) = \Delta f_i (c)~, \mbox{ for all } i \in \{k , j_2, ..., j_n\}~.
\end{align*}
Recall that $j_2 , ..., j_n$ already satisfy \eqref{lambdaeq}, so we must have $\lambda^{(k)}(c) = \lambda(c)$ for every $n + 1 \le k \le m$, since
\begin{align*}
    \lambda^{(k)}(c) f_{j_2}(c) = \Delta f_{j_2} (c) = \lambda(c) f_{j_2}(c)~,~ \mbox{ and }~ f_{j_2}(c) \neq 0~.
\end{align*}

Therefore, there is a single nonzero real number $\lambda(c) \in \R \backslash \{0\}$ such that
\begin{align*}
  \lambda(c) f_j(c) = \Delta f_j(c)~, \mbox{ for all }~ 1 \le j \le m~.
\end{align*}
\end{proof}

\subsection{Proof of Proposition \ref{derivationofmainthm}, (i)}

It suffices to only prove the case in which $Y$ attains the maximum dimension proposed by the hypothesis of Proposition \ref{derivationofmainthm}, since we can always embed $Y$ into a larger subspace if needed.
Suppose $\dim(Y) = \binom{d+ n -2}{n-1} - 1$. Recall that $\dim(\mathcal{P}^{n,d}) = \binom{d+n-1}{n-1}$.
We may directly compute
\begin{align}
  \binom{d+n-1}{n-1} - \binom{d+n-2}{n-2} & = \binom{d+ n -2}{n-1}~.\label{computedimension}
\end{align}
Then, $Y^\perp$ has dimension
\begin{align*}
    m: = \binom{d+n-2}{n-2} + 1  ~.
\end{align*}
Take $\{f_1,..., f_m\}$ to be the orthonormal basis of $Y^\perp$ with respect to the Fischer inner product via setting $\mathcal{M} = Y^\perp$ in Lemma \ref{lemma:subspace}. 

The strategy is to show that there exist indices $\{i_1, ..., i_n\} \subseteq \{1 , ..., m\}$ and an open dense set $U \subseteq \R^n$ such that the following orthogonal projection map
\begin{align}
  I(c): \mathrm{span}_\R \cbracket{\bracket{ c \cdot X}^{d-1}X_j~~|~~ 1 \le j \le n } \longrightarrow \mathrm{span}_\R \{ f_{i_1},..., f_{i_n} \}~, \label{projection2}
\end{align}
is injective for all $c \in U$. Indeed, if this is true, we can define a codimension $n$ subspace
\begin{align*}
    \mathcal{Q} : = \bracket{ \mathrm{span}_\R \{ f_{i_1},..., f_{i_n} \}}^\perp \subseteq \mathcal{P}^{n , d}~.
\end{align*}
Then, $\mathcal{Q}$ contains $Y$ by construction of the basis $\{f_1,..., f_m\}$ of $Y^\perp$. To verify the tangential nondegeneracy condition, for any fixed
$c \in  U$, consider the intersection
\begin{align*}
   \mathcal{Q} \cap \mathrm{span}_\R \cbracket{\bracket{ c \cdot X}^{d-1}X_j~~|~~ 1 \le j \le n }~.
\end{align*}
If $p$ is contained in this intersection, then $I(c)p = 0$ since $p \in \mathcal{Q}$ by definition of an orthogonal projection. Injectivity of $I(c)$ in \eqref{projection2} then implies $p = 0$. This shows that $\mathcal{Q}$ satisfies the tangential nondegeneracy condition \eqref{cond1} for every $c \in U$.

It remains to prove that there exist indices $\{i_1, ..., i_n\}$ and an open dense set $U \subseteq \R^n$ such that $I(c)$ is injective for all $c \in U$. 
By Lemma \ref{lemma:subspace}, there exist rows $i_1,..., i_n$ and an open dense set $U \subseteq \R^n$ such that the $n \times n$ matrix
\begin{align}
  A(c, i_1,...,i_n) : = (d-1)!\sbracket{  \begin{matrix}
    \nabla f_{i_1}(c)  \\
      \vdots \\
     \nabla f_{i_n}(c)
    \end{matrix}} \label{Amatrix}
\end{align}
has nonzero determinant for all $c \in U$. Fix $c \in (\R \backslash \{0\})^n$.
If $f \in \mathcal{P}^{n , d}$, then expanding the polynomial $\bracket{ c \cdot X}^{d-1}X_j$ explicitly via the binomial theorem yields
\begin{align*}
    \inner{\bracket{ c \cdot X}^{d-1}X_j}{ f } = (d-1)!\partial_{j} f(c)~.
\end{align*}
Given a polynomial $f = \sum_{j = 1}^n \lambda_j \bracket{ c \cdot X}^{d-1}X_j$, under the inner product \eqref{innerprod}, we have
\begin{align}
  I(c)f = \sum_{k = 1}^n \inner{f}{f_{i_k}}f_{i_k} = \sum_{k = 1}^n \sum_{j = 1}^n \lambda_j \inner{\bracket{ c \cdot X}^{d-1}X_j}{f_{i_k}} f_{i_k} = (d-1)!\sum_{k = 1}^n \bracket{\sum_{j = 1}^n \lambda_j  \partial_j f_{i_k}(c)} f_{i_k}~. \label{representation}
\end{align}
Thus, the orthogonal projection map $I(c)$ in \eqref{projection2} is represented by the matrix $A(c, i_1,...,i_n)$ in \eqref{Amatrix} with respect to the bases explicitly given by the spans in the definition of $I(c)$.
Since $A(c, i_1,...,i_n)$ has nonzero determinant for all $c \in U$, the map $I(c)$ is injective for all $c \in U \cap (\R \backslash \{0\})^n$.

\subsection{Proof of Proposition \ref{derivationofmainthm}, (ii)}
Again, the goal is to construct a suitable $n$ dimensional subspace $\mathcal{Q}^\perp \subseteq Y^\perp$ using orthogonal projections.
Suppose $d\ge 4$ is even and $\dim(Y) = \binom{d+n-3}{n-1} - 1$. Then we may record that the $\binom{d+n-3}{n-1}$ dimensional subspace $\abs{X}^2\mathcal{P}^{n,d-2}$ is not fully contained in $Y$. 
We directly compute
\begin{align}
  \binom{d+n-1}{n-1} - \binom{d+n-3}{n-1} & = \frac{2d+n-2}{d+n-2}\binom{d+n-2}{n-2}\quad~, \label{computationdeven}
\end{align}
so that $Y^\perp$ has dimension
\begin{align*}
    m : = \frac{2d+n-2}{d+n-2}\binom{d+n-2}{n-2} +1 >  \binom{d+n-2}{n-2}~.
\end{align*}
Let $\{f_1,..., f_m\}$ be the orthonormal basis of $Y^\perp$ with respect to the Fischer inner product \eqref{innerprod} via setting $\mathcal{M} = Y^\perp$ in Lemma \ref{lemma:subspace}. For the same inner product, the vector space $\mathcal{P}^{n,d}$ admits an orthogonal decomposition (see \cite[Section 3.1.2]{MR290095})
\begin{align*}
  \mathcal{P}^{n,d} = \mathcal{H}^{n , d} \oplus   \abs{X}^2\mathcal{P}^{n,d-2}~,
\end{align*}
where $\mathcal{H}^{n , d}$ is the space of harmonic homogeneous polynomials of degree $d$ in $n$ variables, that is, the set of polynomials $h$ such that $\Delta h = 0$. Since $\abs{X}^2\mathcal{P}^{n,d-2}$ is not fully contained in $Y$, we must have $ Y^\perp \not \subseteq \mathcal{H}^{n , d}$.
Without loss of generality, we may assume that $f_1$ is a nonharmonic basis vector of $Y^\perp$.

If there is a $c \in (\R \backslash \{0\})^n$ and indices $\{i_1,..., i_n\} \subseteq \{1, ... , m\}$ such that both orthogonal projection maps
\begin{align}
    I(c)~:& ~\mathrm{span}_\R \cbracket{\bracket{ c \cdot X}^{d-1}X_j~~|~~ 1 \le j \le n } \longrightarrow \mathrm{span}_\R \{ f_{i_1},..., f_{i_n} \}~, \label{orthprojI}\\
    J(c) ~:&~ \mathrm{span}_\R \cbracket{\bracket{c \cdot X}^d , \bracket{c \cdot X}^{d-2}\abs{X}^2} \longrightarrow \mathrm{span}_\R \{ f_{i_1},..., f_{i_n} \}\label{orthproj}
\end{align}
are injective, then taking
\begin{align*}
    \mathcal{Q} :=(\mathrm{span}_\R \{ f_{i_1},..., f_{i_n} \})^\perp
\end{align*}
defines a codimension $n$ subspace of $\mathcal{P}^{n , d}$ containing $Y$ such that $\mathcal{Q}$ satisfies both nondegeneracy conditions \eqref{cond1} and \eqref{cond2} of Theorem \ref{mainresult} for that choice of $c$. We dedicate the rest of the proof in this subsection to finding such $c$ and indices $\{i_1,..., i_n\}$.

As shown by \eqref{representation} in part (i), for $c \in (\R \backslash \{0\})^n$ and indices $\{i_1, ..., i_n\}$, the projection map $I(c)$ in \eqref{orthprojI} is represented by a $n \times n$ matrix $A(c, i_1,..., i_n)$ defined as in \eqref{Amatrix}.
For $J(c)$, similarly, we may compute the following inner products for $f \in \mathcal{P}^{n,d}$ by expanding every term via the binomial theorem,
\begin{align*}
\inner{f}{(c\cdot X)^d} & = d! f(c) ~,~ \inner{f}{(c\cdot X)^{d-2}\abs{X}^2} = (d-2)!\Delta f(c) ~~,
\end{align*}
so that $J(c)$ is represented by the $n \times 2$ matrix
\begin{align}
  B(c, i_1,...,i_n) : = \sbracket{ \begin{matrix} \inner{f_{i_1}}{(c\cdot X)^d} &  \inner{f_{i_1}}{(c\cdot X)^{d-2}\abs{X}^2} \\ 
  \vdots & \vdots \\ 
  \inner{f_{i_n}}{(c\cdot X)^d} &  \inner{f_{i_n}}{(c\cdot X)^{d-2}\abs{X}^2} \end{matrix}} = \sbracket{ \begin{matrix} d! f_{i_1} (c) &  (d-2)!\Delta f_{i_1}(c) \\ 
  \vdots & \vdots \\ 
  d! f_{i_n} (c) &  (d-2)!\Delta f_{i_n}(c) \end{matrix}}~. \label{Bmatrix}
\end{align}
Therefore, to show that there are $c \in (\R \backslash \{0\})^n$ and $\{i_1,..., i_n\}$ such that both $I(c)$ and $J(c)$ are injective, it suffices to find $c \in (\R \backslash \{0\})^n$ and $\{i_1,..., i_n\}$ such that $\det(A(c, i_1,...,i_n)) \neq 0$ and $B(c, i_1,...,i_n)$ has rank $2$.

We argue by contradiction. Suppose for a contradiction that for all $c \in (\R \backslash \{0\})^n$ and all choices of $\{i_1,...,i_n\}$ such that $\det(A(c, i_1,...,i_n)) \neq 0$, the matrix $B(c, i_1,...,i_n)$ has rank strictly less than $2$. First, by Lemma \ref{lemma:subspace}, there is an open dense set $U \subseteq \R^n$ such that the matrix $F(c)$ in \eqref{matrixrep2} has rank $n$ for all $c \in U$, and recall that $f_1$ is a nonharmonic basis element of $Y^\perp$. Since the complement of the zero set of a nonzero polynomial is open and dense in $\R^n$, there is an open dense set $U'$ in $\R^n$ with $ U' \subseteq U$, such that 
\begin{align*}
    \Delta f_1(c)~, ~f_j(c) \neq 0 ~,~ \text{ for all } 1 \le j \le m~,~ c \in U'~.
\end{align*}

Fix $c \in U' $ for now. Then $c$ and the polynomials $f_1 ,... , f_m$ satisfy the hypothesis of Lemma \ref{laplaciamlemma}. By the Lemma, there exists a constant $\lambda(c) \in \R \backslash \{0\}$ such that 
  \begin{align}
      \Delta f_i(c) = \lambda(c) f_i(c) ~,~ \mbox{ for all } 1 \le i \le m~.\label{equalindentity}
  \end{align}
Let $f , g\in \{f_1,...,f_m\}$. By \eqref{equalindentity}, the following identity about the basis elements of $Y^\perp$ holds on the open and dense set $U' \subseteq \R^n$:
\begin{align}
    f\Delta g = g\Delta f~.~\label{equation}
\end{align}
Since a polynomial vanishing on a nonempty open set vanishes identically, the above equation holds as a polynomial identity. Fix a basis element $g$ of $Y^\perp$. Let $S$ be the set of all solutions $f$ to \eqref{equation} in $\mathcal{P}^{n , d}$. Then, $S$ is a linear subspace of $\mathcal{P}^{n , d}$ and we have
\begin{align*}
    Y^\perp \subseteq S~.
\end{align*}

We proceed to derive a contradiction by showing that $ \dim(S) < \dim(Y^\perp)$.
First, since $g$ in \eqref{equation} is nonzero, there is a hyperplane $H \subseteq \R^n$ such that $g|_H$ is not the zero function. Perform an orthogonal change of coordinates $B:\R^n \to \R^n$, $z = Bx$, which rotates the hyperplane $H$ to the hyperplane ${\{z \in \R^{n}~~|~~ z_n = 0\}}$, so that for each $f \in \mathcal{P}^{n , d}$,
\begin{align*}
  f(Z) = f(BX)~,~ f\circ B^T|_{\{Z_n = 0\}} = f|_H~.
\end{align*}
Since the change of coordinates is linear, $f(Z)$ remains a homogeneous polynomial of degree $d$ in $Z$.
In the new coordinates, every $f \in \mathcal{P}^{n , d}$ can be expressed in the following form
\begin{align*}
    f(Z) = \sum_{k = 0}^{d} Z_n^k f^{(k)} \bracket{ Z^{(n)}}~,~ f^{(k)} \in \mathcal{P}^{n-1 , d-k}~.
\end{align*}
Expressing $g$ in this form, we obtain that $g^{(0)} \neq 0$ as a polynomial in $\mathcal{P}^{n-1 , d}$, since $g(z)$ is not identically zero on ${\{z \in \R^{n}~~|~~ z_n = 0\}}$ by construction. Now consider the map
\begin{align*}
    \Phi: S \longrightarrow \mathcal{P}^{n-1 , d} \times \mathcal{P}^{n-1 , d-1}~,~ f \longmapsto (f^{(0)} , f^{(1)})~.
\end{align*}
We claim that $\Phi$ is injective on $S$. If this is true, then
\begin{align*}
  m = \dim(Y^\perp) \le \dim(S) & \le \dim(\mathcal{P}^{n-1 , d}) + \dim(\mathcal{P}^{n-1 , d-1})~
\end{align*}
gives an upper bound on $m$, but combining with \eqref{computationdeven}, we obtain that this upper bound is simultaneously less than $m$:
\begin{align*}
  \dim(\mathcal{P}^{n-1 , d}) + \dim(\mathcal{P}^{n-1 , d-1}) & = \binom{d+n-2}{n-2} + \binom{d+n-3}{n-2} \\& = \binom{d+n-1}{n-1} - \binom{d+n-3}{n-1} = m-1~.
\end{align*}
This is a contradiction.

It remains to show that $\Phi$ is injective. Let $\Delta_X$ denote the Laplacian with respect to the variable $X$ and $\Delta_Z$ denote the Laplacian with respect to the variable $Z$. Since the Laplacian is rotation invariant, we have
\begin{align*}
    B^T \circ \Delta_Z f = \Delta_X f\circ B^T~,
\end{align*}
so equation \eqref{equation} is unchanged under the orthogonal change of coordinates.
Suppose $f \in S$ satisfies $f^{(0)} = f^{(1)} = 0$ and $f \neq 0$.  Let $r \ge 2$ be the smallest $r$ such that $f^{(r)} \neq 0$. Then we can compute $g\Delta_Z f$:
\begin{align*}
  \Delta_Z f(Z) & =  \Delta_Z \sum_{k = r}^{d} Z_n^k f^{(k)} \bracket{ Z^{(n)}} = r(r-1)Z_n^{r-2}f^{(r)} \bracket{ Z^{(n)}} + O(Z_n^{r-1})~,\\
  g(Z) \Delta_Z f \bracket{ Z} & =  r(r-1) g^{(0)} \bracket{ Z^{(n)} }Z_n^{r-2}f^{(r)} \bracket{ Z^{(n)}}+ O(Z_n^{r-1})~.
\end{align*}
On the other hand, since all monomials in $Z_n$ in $f$ have degree at least $r $ and $g^{(0)} \neq 0$, the coefficient of $Z_n^{r-2}$ in $f\Delta_Z g$ must be zero for $f$ to satisfy \eqref{equation}. Hence,
\begin{align*}
    r(r-1) g^{(0)} \bracket{ Z^{(n)} }f^{(r)} \bracket{ Z^{(n)}} = 0~.
\end{align*}
Since $g^{(0)} \neq 0$, we conclude that $f^{(r)} = 0$ as a polynomial, which is a contradiction to the definition of $r$. Hence, $\Phi$ must be injective on $S$.

In summary, we conclude that there is a choice of $c \in (\R \backslash \{0\})^n$ and indices $\{i_1, ..., i_n\}$ such that $\mathcal{Q} =(\mathrm{span}_\R \{ f_{i_1},..., f_{i_n} \})^\perp$ satisfies both nondegeneracy conditions \eqref{cond1} and \eqref{cond2} for $c$ and contains $Y$.

\subsection{Proof of Proposition \ref{derivationofmainthm}, (iii)} Assume $d = 2$, $\dim(Y) = n-1$, and $|X|^2 \notin Y$. Notice that
\begin{align*}
  n - 1 = \binom{n}{n-1} -1 = \binom{d+ n - 2}{n-1} - 1~,
\end{align*}
So by the same computation as in \eqref{computedimension}, we have $\dim(Y^\perp) > \binom{d+ n - 2}{n-2}$. By Lemma \ref{lemma:subspace}, fix an orthonormal basis $\{f_1, \dots, f_m\}$ of $Y^\perp$ such that there exists an open dense set $U \subseteq \R^n$ such that for all $c \in U$, the matrix $F(c)$ in \eqref{matrixrep2} has rank $n$.

Again, it suffices to find $c \in (\R \backslash \{0\})^n$ and indices $i_1, ..., i_n$ such that the matrix $A(c,i_1, ..., i_n)$ in \eqref{Amatrix} has nonzero determinant, while the $B(c,i_1, ..., i_n)$ in \eqref{Bmatrix} has rank $2$. Once this is established, we can define
\begin{align*}
    \mathcal{Q} : = (\mathrm{span}_\R \{f_{i_1} , ... ,f_{i_n}\})^\perp~.
\end{align*}
By contradiction, suppose that for each $c \in (\R \backslash \{0\})^n$ and each $i_1, ..., i_n$ such that $\det(A(c , i_1, ..., i_n)) \neq 0$,
the rank of $B(c, i_1, ..., i_n)$ is strictly less than $2$. By the same argument as Proposition \ref{derivationofmainthm}, part (ii), we get that for any basis elements $f , g$ of $Y^\perp$, the following identity holds:
\begin{align*}
    f\Delta g = g\Delta f~. 
\end{align*}
Fix a basis element $g \in \mathcal{P}^{n , d}$ of $Y^\perp$, and 
let $S$ be the set of solutions to the above equation in $\mathcal{P}^{ n , 2}$, so that $Y^\perp \subseteq S$. 

We proceed to derive a contradiction by showing that there is a basis element $g$ of $Y^\perp$ such that $\dim(S) < \dim(Y^\perp)$.
Since we have the orthogonal decomposition (see \cite[Section 3.1.2]{MR290095})
\begin{align*}
    \mathcal{P}^{n , 2} = \mathcal{H}^{n , 2} \oplus  \mathrm{span}_\R \{ |X|^2 \} ~,
\end{align*}
and we have assumed $|X|^2 \notin Y$,
there is at least one element in $Y^\perp$ that is not harmonic. Take $g$ to be a nonharmonic basis element of $Y^\perp$, so that $\Delta g \in \R \backslash \{0\}$ (this is a real number since $d = 2$). 
Since $d = 2$, we have $\Delta f_j \in \R$ for all $1 \le j \le m$. It follows that
\begin{align*}
  f_j = \frac{\Delta f_j}{\Delta g} g~,~ \mbox{ for all } 1 \le j \le m~,
\end{align*}
which implies $\dim(Y^\perp) = 1$. This is a contradiction to the assumption
\begin{align*}
   \dim(Y^\perp) > \binom{d+ n - 2}{n-2} =  \binom{n}{n-2} \ge 1~.
\end{align*}

\begin{rem} \label{conditions}
  Fix the inner product \eqref{innerprod} on $\mathcal{P}^{n,d}$. Let $\mathcal{X}, \mathcal{M} \subseteq \mathcal{P}^{n , d}$ be linear subspaces, and assume that $\mathcal{M}$ is equipped with a basis (not necessarily orthonormal) $\{g_1,..., g _m\}$, and $\mathcal{X}^\perp$ is equipped with an orthonormal basis $\{f_1,..., f_s\}$.
  In the proof of Proposition \ref{derivationofmainthm}, we used the observation that $\mathcal{X} \cap \mathcal{M} = \{0\}$ if the orthogonal projection of $\mathcal{M}$ onto $\mathcal{X}^\perp$,
  \begin{align*}
      \pi : \mathcal{M} & \longrightarrow \mathcal{X}^\perp\\
      p = \sum_j \lambda_j g_j & \longmapsto \sum_{i} \bracket{\sum_j \lambda_j  \langle g_j, f_i \rangle} f_i = \sum_i \langle p, f_i \rangle f_i
  \end{align*}
is injective.
When $\mathcal{Q}$ is a codimension $n$ subspace of $\mathcal{P}^{n,d}$, its orthogonal complement has dimension $n$, and we can choose an orthonormal basis $\{f_1,... ,f_n\}$. When $c \in (\R \backslash \{0\})^n$ and
\begin{align*}
    \mathcal{M} = \mathrm{span}_\R \cbracket{ X_j(c \cdot X)^{d-1}~~|~~ 1 \le j \le n}~,
\end{align*}
under the fixed inner product \eqref{innerprod},
the matrix which represents the projection map $\pi$ (with respect to the bases $\{X_j(c \cdot X)^{d-1}\}$ and $\{f_i\}$) is
\begin{align*}
    (d-1)!\sbracket{  \begin{matrix}
    \nabla f_{1}(c)  \\
      \vdots \\
     \nabla f_{n}(c)
    \end{matrix}}~.
\end{align*}
We named condition \eqref{cond1} tangential because this representation involves the gradients of the basis polynomials of $\mathcal{Q}^\perp$ at the point $c$. Similarly, when
\begin{align*}
    \mathcal{M} = \mathrm{span}_\R \cbracket{ (c \cdot X)^d~,~ \abs{X}^{2}(c \cdot X)^{d - 2}}~,
\end{align*}
the matrix which represents the projection map $\pi$ (with respect to the bases $\{(c \cdot X)^d, \abs{X}^{2}(c \cdot X)^{d - 2}\}$ and $\{f_i\}$) is
\begin{align*}
    \sbracket{ \begin{matrix} d! f_{1} (c) &  (d-2)!\Delta f_{1}(c) \\ 
  \vdots & \vdots \\ 
  d! f_{n} (c) &  (d-2)!\Delta f_{n}(c) \end{matrix}}~.
\end{align*}
We named condition \eqref{cond2} radial because this representation involves the Laplacians of the basis polynomials of $\mathcal{Q}^\perp$ at the point $c$.

\end{rem}

\section{Proof of Proposition \ref{monomials}} \label{monomialsec}

\begin{proof}[Proof of Proposition \ref{monomials}]

Fix the inner product on $\mathcal{P}^{n , d}$ defined as in \eqref{innerprod}. Then, the monomials $X^{\alpha}$, $\abs{\alpha} = d$, form an orthogonal basis for $\mathcal{P}^{n,d}$. For a codimension $n$ subspace $\mathcal{Q}$ equipped with a monomial basis, we have
\begin{align*}
    \mathcal{Q}^\perp = \mathrm{span}_\R\{ X^{\theta(i)}~~|~~ 1 \le i \le n\}~,
\end{align*}
where $X^{\theta(i)}$, $1 \le i \le n$ are the monomials in $\mathcal{P}^{n,d}$ that are omitted from the basis of $\mathcal{Q}$.

Assume that the set $T := \{\theta(i)~~|~~ 1 \le i \le n\}$ of the omitted multi-indices is linearly independent in $\R^n$. We will show that for all $c \in (\R \backslash \{0\})^n$, the subspace $\mathcal{Q}$ satisfies the tangential nondegeneracy condition \eqref{cond1}.
Consider the orthogonal projection map
\begin{align}
  \pi: \mathrm{span}_\R\{ (c \cdot X)^{d-1}X_1 ,..., (c \cdot X)^{d-1}X_n\} &\longrightarrow \mathcal{Q}^\perp = \mathrm{span}_\R\{ X^{\theta(i)}~~|~~ 1 \le i \le n\}~ \nonumber \\ 
  p &\longmapsto \sum_{i = 1}^{n}[X^{\theta(i)}](p(X)) X^{\theta(i)}~.\label{C}
\end{align}
Here, $[X^{\theta(i)}](p(X))$ denotes the coefficient of the monomial $X^{\theta(i)}$ in the polynomial $p(X)$.
Take any $c\in (\R\backslash \{0\})^n$.
If $\pi$ is injective, then $\mathcal{Q}$ satisfies the tangential nondegeneracy condition \eqref{cond1} for $c$. Indeed, if
\begin{align*}
    p \in \mathcal{Q} \cap \mathrm{span}_\R\{ (c \cdot X)^{d-1}X_1 ,..., (c \cdot X)^{d-1}X_n\}
\end{align*}
then $\pi(p) = 0$ since $p \in \mathcal{Q}$, and so $p = 0$ by injectivity.

To verify injectivity,
consider the $n \times n$ real-valued matrix $C = (C_{i,j})_{1 \le i,j \le n}$ where
\begin{align*}
  C_{i , j} = [X^{\theta(i)}](c \cdot X)^{d-1}X_j = c^{\theta(i) - e_j} \binom{d}{\theta(i)} \frac{1}{d} \theta(i)_j ~.
\end{align*}
Then $C$ represents $\pi$ in coordinates given by the basis specified by the spans in \eqref{C}. In particular, we have 
\begin{align*}
  \det(C) = d^{-n} c^{\sum_{i = 1}^{n} \theta(i)} \prod_{i = 1}^n \binom{d}{\theta(i)} c^{- \sum_{j = 1}^{n} e_j} \det \bracket{ \bracket{\theta(i)_j}_{1\le i,j \le n} }~.
\end{align*}
This determinant is nonzero,
since $T = \{\theta(i)~~|~~ 1 \le i \le n\}$ is a linearly independent set in $\R^n$. Hence, we have shown that for any $c \in (\R\backslash \{0\})^n$, the map $\pi$ is injective.

Next, in the case when $d\ge 2$ is even, we want to find $c \in (\R\backslash \{0\})^n$ such that the monomial span $\mathcal{Q}$ satisfies condition \eqref{cond2}. Here, we will use the hypothesis that $T$ is a linearly independent set in $\R^n$, and that there is $\theta \in T$ such that $\theta_j \ge 2$, for some $1 \le j \le n$.
Fix an arbitrary $c \in (\R\backslash \{0\})^n$ for now and consider
\begin{align*}
  p \in \mathcal{Q} \cap \mathrm{span}_\R  \cbracket{ (c \cdot X)^d~,~ \abs{X}^{2}(c \cdot X)^{d - 2}} ~.
\end{align*}
We write $p(X) = a (c \cdot X)^d + b \abs{X}^{2}(c \cdot X)^{d - 2}$ for $a,b \in \R$. The goal is to show that there is a suitable $c \in (\R \backslash \{0\})^n$ that forces $a = b = 0$, which implies $p=0$ as a polynomial. Since $p \in \mathcal{Q}$, we must have for all $\theta \in T$,
\begin{align}
  0 = [X^{\theta}]p(X) = \frac{(d-2)!}{\theta!}c^\theta \bracket{ad(d-1) + b \sum_{ j =1}^{n} \theta_j (\theta_j -1) c_j^{-2} }~. \label{coeff}
\end{align}
Hence, the constants $a,b$ must satisfy the equation
\begin{align}
  ad(d-1) + b \sum_{ j =1}^{n} \theta_j (\theta_j -1) c_j^{-2} = 0~\mbox{ , for every } \theta \in T~. \label{intersection}
\end{align}

We claim that there exists $c \in (\R\backslash \{0\})^n$ such that the function $\theta \mapsto \sum_{ j =1}^{n} \theta_j (\theta_j -1) c_j^{-2}$ is nonconstant on $T$. For this $c$, the only way in which \eqref{intersection} is satisfied uniformly over $\theta \in T$ is when $a = b = 0$. 
Define $ q (\theta) : = (\theta_1(\theta_1 - 1) ,..., \theta_n (\theta_n - 1))$, acting on $\theta \in \R^n$. First, there must exist $\theta , \eta \in T $ such that $q (\theta) \neq q(\eta)$. By contradiction, suppose not. Then, $q (\theta)$ is constant as $\theta$ varies in $T$. By assumption, $T$ is not entirely square-free, so there is an index $1\le j \le n$ such that $\theta_j(\theta_j - 1) =:  h > 0$ for all $\theta \in T$. Using the quadratic formula to solve for $\theta_j$, we get that 
\begin{align*}
  \theta_j = \frac{1 \pm \sqrt{1 + 4h}}{2}~\mbox{ for all } \theta \in T~.
\end{align*}
There is only one permissible solution for $\theta_j \ge 0$, which means $\theta_j$ is constant as $\theta$ varies in $T$. Moreover, we have $\abs{\theta} = d$ for each $\theta \in T$. So the $n \times n$ matrix
\begin{align*}
    \sbracket{ \begin{matrix}
      \theta(1)_1 & \theta(1)_2 & \cdots & \theta(1)_n \\
      \vdots & \vdots & \ddots & \vdots \\
      \theta(n)_1 & \theta(n)_2 & \cdots & \theta(n)_n  
    \end{matrix}}
\end{align*}
has linearly dependent columns, 
contradicting linear independence of the set of multi-indices $T = \{\theta(i)\}$. Thus, we conclude that there are $\theta , \eta \in T$ such that $q(\theta) \neq q(\eta)$.

Next, take $\theta , \eta \in T$ such that $q(\theta)_j \neq q(\eta)_j$, for some index $1\le j \le n$. Define $c \in (\R\backslash \{0\})^n$ so that $c_i = 1$ for all $i \neq j$, and $c_j = \epsilon^{1/2}$ for an $\epsilon > 0$ of our choice. Then
\begin{align*}
  \sum_{k = 1}^{n} (\theta_k (\theta_k - 1) - \eta_k(\eta_k - 1)) c_k^{-2} = (\theta_j (\theta_j - 1) - \eta_j(\eta_j - 1))\frac{1}{\epsilon} + \sum_{k \neq j}^{} (\theta_k (\theta_k - 1) - \eta_k(\eta_k - 1))~.
\end{align*}
The second term on the right-hand side is bounded in absolute value by a constant depending only on $n$ and $d$. We take $\epsilon> 0$ to be small enough so that the entire sum is nonzero. In this way, we have constructed $c \in (\R\backslash \{0\})^n$ so that the function $\theta \mapsto \sum_{ j =1}^{n} \theta_j (\theta_j -1) c_j^{-2}$ is nonconstant on $T$. Therefore, the only way in which \eqref{intersection} is satisfied uniformly over $\theta \in T$ is when $a = b = 0$. This completes the proof that
\begin{align*}
  \mathcal{Q} \cap \mathrm{span}_\R  \cbracket{ (c \cdot X)^d~,~ \abs{X}^{2}(c \cdot X)^{d - 2}} = \{0\}~.
\end{align*}

\end{proof}
\bibliographystyle{alpha} 
\bibliography{bibliography}

\end{document}